\documentclass{amsart}
\usepackage{amsmath,amssymb,amsthm,mathrsfs,mathtools,tikz-cd,hyperref}
\usepackage[capitalise,noabbrev]{cleveref}
\usepackage{dsfont}  
\usepackage{enumerate}
\usepackage{esint}   
\numberwithin{equation}{section}

\hypersetup{colorlinks=true, 	linkcolor=blue!75!black, 	citecolor=red!50!black, 	urlcolor=green!50!black, 	linktoc=all }

\usepackage{colortbl}

\def\NN{\mathbb{N}}
\def\RR{\mathbb{R}}
\def\ZZ{\mathbb{Z}}
\def\TT{\mathbb{T}}

\newtheorem{theorem}{Theorem}[section]
\newtheorem{proposition}[theorem]{Proposition}
\newtheorem{corollary}[theorem]{Corollary}

\newtheorem{lemma}[theorem]{Lemma}
\theoremstyle{definition}
\newtheorem{definition}[theorem]{Definition}
\newtheorem{problem}[theorem]{Problem}
\theoremstyle{remark}
\newtheorem{remark}[theorem]{Remark}

\newcommand{\norm}[1]{\left\|#1\right\|}
\newcommand{\abs}[1]{\left\vert#1\right\vert}
\newcommand{\T}[1]{\widetilde{#1}}
\newcommand{\R}[1]{\mathring{#1}}
\DeclareMathOperator\Div{div}
\DeclareMathOperator\curl{curl}
\DeclareMathOperator\Id{Id}
\DeclareMathOperator\dist{dist}
\DeclareMathOperator\supp{supp}
\DeclareMathOperator\tr{tr}
\def\gradsim{\nabla_{\text{sym}}}
\newcommand{\olsi}[1]{\,\overline{\!{#1}}} 

\newcommand{\Rbar}[1]{\R{\olsi{R}}_q\strut^{\hspace{-3.5pt}(#1)}}
\newcommand{\vtildei}[1]{\T{v}_i^{\hspace{0.5pt}#1}}
\newcommand{\ptildei}[1]{\T{p}_i^{\hspace{0.5pt}#1}}
\newcommand{\ztildei}[1]{\T{z}_i^{\hspace{0.5pt}#1}}
\newcommand{\Rtildei}[1]{\R{\T{R}}_i\strut^{\hspace{-3.3pt}#1}}

\def\dim{n}
\def\Proy{P}
\def\normal{\nu}
\def\hodge{\star}

\usepackage{subfiles}

\begin{document}

\title[An Extension theorem for weak solutions of the 3d Euler equations]{An extension theorem for weak solutions\\ of the 3d incompressible Euler equations\\ and applications to singular flows}

\author[A.~Enciso, J.~Pe\~nafiel-Tom\'as and D.~Peralta-Salas]{Alberto Enciso, Javier Pe\~nafiel-Tom\'as and Daniel Peralta-Salas}

\address{
\newline
\textbf{{\small Alberto Enciso}} \vspace{0.1cm}
\vspace{0.1cm}
\newline \indent Instituto de Ciencias Matem\'aticas, Consejo Superior de Investigaciones Cient\'ificas, 28049 Madrid, Spain}
\email{aenciso@icmat.es}

\address{
\vspace{-0.25cm}
\newline
\textbf{{\small Javier Pe\~nafiel-Tom\'as}} 
\vspace{0.1cm}
\newline \indent Instituto de Ciencias Matem\'aticas, Consejo Superior de Investigaciones Cient\'ificas, 28049 Madrid, Spain}
\email{javier.penafiel@icmat.es}

\address{
\vspace{-0.25cm}
\newline
\textbf{{\small Daniel Peralta-Salas}} \vspace{0.1cm}
\vspace{0.1cm}
\newline \indent Instituto de Ciencias Matem\'aticas, Consejo Superior de Investigaciones Cient\'ificas, 28049 Madrid, Spain}
\email{dperalta@icmat.es}

\begin{abstract}
We prove an extension theorem for local solutions of the 3d incompressible Euler equations. More precisely, we show that if a smooth vector field satisfies the Euler equations in a spacetime region $\Omega\times(0,T)$, one can choose an admissible weak solution on $\mathbb R^3\times (0,T)$ of class $C^\beta$ for any $\beta<1/3$ such that both fields coincide on~$\Omega\times (0,T)$. Moreover, one controls the spatial support of the global solution. Our proof makes use of a new extension theorem for local subsolutions of the incompressible Euler equations and a $C^{1/3}$ convex integration scheme implemented in the context of weak solutions with compact support in space. We present two nontrivial applications of these ideas. First, we construct infinitely many admissible weak solutions of class $C^\beta_{\text{loc}}$ with the same vortex sheet initial data, which coincide with it at each time~$t$ outside a turbulent region of width $O(t)$. Second, given any smooth solution~$v$ of the Euler equation on $\TT^3\times(0,T)$ and any open set~$U\subset\TT^3$, we construct admissible weak solutions which coincide with $v$ outside~$U$ and are uniformly close to it everywhere at time~$0$, yet blow up dramatically on a subset of~$U\times (0,T)$ of full Hausdorff dimension. These solutions are of class~$C^\beta$ outside their singular set.
\end{abstract}

\maketitle

\setcounter{tocdepth}{1}
\tableofcontents
\section{Introduction}


Convex integration methods, introduced by Nash~\cite{Nash} in the context of the $C^1$ isometric embedding problem and subsequently refined by Gromov in his work on flexible geometric PDEs and by M\"uller and \v{S}ver\'ak~\cite{MullerSverak} in their theory of differential inclusions, have experimented an extraordinary development in connection with the study of weak solutions of the incompressible Euler equations. This system reads as
\[
\partial_tv+\Div(v\otimes v) +\nabla p=0\,,\qquad \Div v=0\,,
\]
where the time-dependent vector field $v$ is the velocity of the fluid and the scalar function $p$ is the hydrodynamic pressure. One typically considers the Euler equations either on the whole space~$\RR^3$, or on the torus $\TT^3:=(\RR/\ZZ)^3$, or on a bounded domain $\Omega\subset\RR^3$ with smooth boundary (where additional complications may arise).

The motivation to consider weak solutions in this setting is twofold. First, the 3d Euler equations are expected to dynamically produce singularities from smooth initial conditions~\cite{Hou,JGS}. Second, weak solutions are necessary to describe some of the phenomena that appear in turbulence, such as the energy dissipation in non-smooth Euler flows famously conjectured by Onsager in 1949~\cite{Onsager49}. Roughly speaking, Onsager's conjecture asserts that weak solutions that are H\"older continuous in space with exponent greater than $1/3$ must conserve energy, while for any smaller exponent  there should be weak solutions that do not.

The rigidity part of Onsager's conjecture was proved by Constantin, E and Titi~\cite{CET} after a partial result of Eyink~\cite{Eyink}. The endpoint case was addressed  in~\cite{CCFS}. Concerning the flexible part of the conjecture, following the construction of $L^2$~solutions with compact support in space and time due to Scheffer~\cite{scheffer1993} and Shnirelman~\cite{shnirelman97},  a systematic approach was introduced in the seminal work of De Lellis and Sz\'ekelyhidi, who introduced $L^\infty$-convex integration and $C^0$-Nash iteration schemes in this setting~\cite{Bounded,Continuous}. After a series of significant intermediate results~\cite{Daneri, Quinto, Isett}, the flexible part of Onsager's conjecture was finally established by Isett~\cite{Isett}, and further refined by Buckmaster, De Lellis, Sz\'ekelyhidi and Vicol~\cite{Onsager_final} to construct solutions for which the kinetic energy is strictly decreasing. In addition to the classical H\"older-based approach, the so-called intermittent $L^p$-based flavor of convex integration, introduced by Buckmaster and Vicol~\cite{VicolNS} to prove the non-uniqueness of weak solutions of the 3d Navier--Stokes equations, has also attracted much attention, as it can effectively capture new aspects of Kolmogorov's theory of turbulence. For detailed expositions of these and other results on various models in fluid mechanics, we refer the reader to the surveys~\cite{ VicolBAMS, LellisBAMS17} and the papers~\cite{ Colombo, VicolBook,   SQG, CCF, CheskiShy14, Cheskidov23, ChesShv, Faraco, Faraco2,NovackVicol, ShvJAMS}.

A key property of the solutions that one constructs using convex integration techniques is their {\em flexibility}\/. This refers to the fact that, at a certain regularity level, the equations are no longer predictive: there exist infinitely many solutions, in stark contrast to what happens in the case of smooth solutions. Three possible formulations of this property are as follows; as discussed in~\cite[Remark 1.3]{NovackVicol}, once one of them has been established within a certain functional framework, it is usually straightforward to pass to another formulation using techniques that are now standard. 

Restricting to the case of~$\TT^3$ for concreteness, let us denote by $\mathcal V\subset L^2(\TT^3)$ some suitable function space, which in our case will be some H\"older space $C^\beta(\TT^3)$. Three standard ways of stating the flexibility of weak solutions in this regularity class are as follows:
\begin{enumerate}
    \item {\em Solutions of compact time support:}\/ Given any positive constants $E,T$, there exists a weak solution $v\in C([-T,T], \mathcal V)$ whose time support is contained in $(-T,T)$ and such that $\|v(0)\|_{L^2(\TT^3)}>E$.

    \item {\em Solutions with fixed energy profile:}\/ Given any smooth positive function $e:[0,T]\to(0,\infty)$, there exists a weak solution $v\in C([0,T],\mathcal V)$ such that $\|v(t)\|_{L^2(\TT^3)}=e(t)$ for all $t\in[0,T]$.

    \item {\em Arbitrary initial and final states:}\/ Given any divergence-free vector fields $v_{\mathrm{start}},v_{\mathrm{end}}\in \mathcal V$ with the same mean, any $T>0$ and any $\epsilon>0$, there exists a weak solution $v\in C([0,T],\mathcal V)$ such that
    \begin{equation}\label{E.introepsilon}
        \|v(0)-v_{\mathrm{start}}\|_{L^2(\TT^3)}+\|v(T)-v_{\mathrm{end}}\|_{L^2(\TT^3)}<\epsilon \,.
    \end{equation}
(If $v_{\mathrm{start}},v_{\mathrm{end}}$ are smooth, one can take $\epsilon=0$ by gluing in time.)
\end{enumerate}

\subsection{Main result}

Our objective in this paper is to prove an extension theorem for local solutions of the 3d incompressible Euler equations. Roughly speaking, we prove that if a smooth vector field satisfies the Euler equations in a spacetime region $\Omega\times(0,T)$ (so it is a ``local'' solution of Euler), one can choose a weak solution on $\mathbb R^3\times (0,\infty)$ of class $C^\beta$ for any $\beta<1/3$ (which is the sharp H\"older regularity) such that both fields coincide on~$\Omega\times (0,T)$. Moreover, one controls the spatial support of the ``global solution'' which extends the local one.

This property is very different from the approximation theorems that one can prove for smooth solutions of various classes of linear PDEs~\cite{aprox1,aprox2,aprox3,aprox4}, and also from the fact (often known as h-principle) that weak solutions of certain regularity can approximate, in Sobolev spaces of negative index, any given subsolution of the Euler equations. 

Before presenting this result, let us recall the definition of weak solution (which, as we will be dealing with continuous functions exclusively, is just the distributional one).
	More precisely, given some $T\in(0,\infty]$ and some open set $\Omega\subseteq\RR^3$ with smooth boundary,  we will say that a vector field $v\in C(\Omega\times [0,T),\RR^3)$  is a {\em weak solution}\/ of the Euler equations on $\Omega\times (0,T)$ if 
	\[\int_0^T\int_{\Omega}(\partial_t\varphi \cdot v+\nabla\varphi:(v\otimes v))dxdt=0\]
	for all divergence-free $\varphi\in C^\infty_c(\Omega\times(0,T),\RR^3)$, and $\Div v=0$ in the sense of distributions. 

The main result of this paper can then be stated as follows:

\begin{theorem}
	\label{teorema gordo}
	Fix some $T>0$ and a bounded open set $\Omega\subset \RR^3$ with smooth boundary and with a finite number of connected components. Assume that $v_0\in C^{\infty}(\overline{\Omega}\times[0,T],\RR^3)$ is a solution of the Euler equations on the spacetime region $\Omega\times(0,T)$. Then, for any $0<\beta<1/3$, there exists an admissible weak solution $v \in C^\beta(\RR^3\times[0,T])$ of the Euler equations such that $v|_{\overline{\Omega}\times[0,T]}=v_0$ if and only if 
	\begin{equation}
		\int_{\Sigma}v_0\cdot \normal=\int_{\Sigma}\left[(a\cdot x)\partial_t v_0+(a\cdot v_0)v_0+p_0 a\right]\cdot \normal= 0  \label{condiciones compatibilidad intro}
	\end{equation}
 for all $a\in\RR^3$, all $t\in[0,T]$ and all connected components $\Sigma$ of $\partial \Omega$. These conditions are automatically satisfied if $\partial \Omega$ is connected. Furthermore, there exists $e_0>0$ such that we may prescribe any energy profile $e\in C^\infty([0,T],[e_0,+\infty))$, that is, $\|v(t)\|_{L^2(\RR^3)}=e(t)$. In addition, given any open set $\Omega'\supset \overline\Omega$, one can in fact assume that the spatial support of~$v$ is contained in this region. 
\end{theorem}

\begin{remark}\label{R.teoremagordo}
In fact, one can obtain a global weak solution $v\in C^\beta(\RR^3\times[0,+\infty))$ such that $v|_{\RR^3\times [0,T]}$ satisfies the claimed properties. Its temporal support may be assumed to be contained in $[0,T']$ for any $T'>T$. We cannot then prescribe the energy profile for $t>T$, but we can still chose~$v$ so that it remains admissible, i.e.,
	\[\int_{\RR^3}\abs{v(x,t)}^2dx\leq \int_{\RR^3}\abs{v(x,0)}^2dx \qquad \forall t\in[0,+\infty).\]
\end{remark}


\begin{remark}
	\label{remark presión}
	It can be proved~\cite{pressure} that the pressure function $p:=-\Delta^{-1} \Div\Div(v\otimes v)$ associated to this weak solution is in
 $ L^\infty_tC^{2\beta}_x\cap C^{2\beta-\delta}_{xt}$ for any $\delta>0$. 
\end{remark}

Before moving on to discuss some applications, let us provide some intuition about the compatibility conditions (\ref{condiciones compatibilidad intro}). When $\partial\Omega$ is connected, it is easy to see that any smooth Euler flow~$v_0$ on~$\Omega$ satisfies this condition. Indeed, these two conditions are respectively obtained by integrating over the domain~$\Omega$ the incompressibility condition $\Div v_0=0$ and the projected Euler equation $a\cdot(\partial_t v_0+\Div(v_0\otimes v_0)+\nabla p_0)=0$. If $v_0$ is the restriction to~$\Omega$ of a global Euler flow, one can refine the argument to show that these conditions must hold on each connected component~$\Sigma$ of the boundary~$\partial\Omega$, and not every field satisfying the Euler equations on~$\Omega$ will satisfy them. Details are given in \cref{pegado subsoluciones}.





\subsection{Applications}

We shall next present two applications of the above extension result to the analysis of weak solutions of the 3d Euler equations. These applications do not follow directly from our main theorem, but they use it in an essential way.

Specifically, for these applications we consider subsolutions that are not smooth up to the endpoints of the interval $(0,T)$, which implies a lack of uniform-in-time bounds. Thus the scheme does not work as is because the available bounds are not uniform, but we will show in Section~\ref{section intervalo abierto} that one can tweak the construction in many interesting situations.

The first application we consider concerns the case of the standard vortex sheet $u_0$, which we can define as the periodic extension to $\mathbb T^3$ of:
\[u_0(x):=\begin{cases}+e_1\quad \text{if }x_3\in\left[0,\frac{1}{4}\right]\cup \left[\frac{3}{4},1\right],\\[3pt] -e_1\quad \text{if }x_3\in\left(\frac{1}{4},\frac{3}{4}\right).\end{cases}\]
It follows from the classical local existence results and from the weak-strong uniqueness property~\cite{Brenier, weak_strong} that wild initial data must be somewhat irregular. However, until the publication of~\cite{vortex_sheet} it was not known how irregular they must be. In that paper it was proved that the vortex sheet $u_0$ is a wild initial data but the constructed solutions are only in $L^\infty$. Results for non-flat vortex sheets have been recently established in~\cite{Mengual}.

One can use a suitable modification of our main theorem to extend this result to solutions of class $C^\beta_\text{loc}$:
\begin{theorem}
	\label{teorema vortex sheet}
	Let $0<\beta<1/3$ and let $T>0$. There exist infinitely many admissible weak solutions of the Euler equations $v \in C^\beta_\text{loc}(\mathbb{T}^3\times(0,T))$ with initial datum $u_0$. For all $t\in(0,T)$, $v(x,t)$ coincides with~$u_0(x)$ outside a ``turbulent'' zone of size $O(t)$.
\end{theorem}

The second application we will present concerns the existence of a wealth of reasonably well behaved solutions that blow up on a set of maximal Hausdorff dimension. To make this precise, let us say that a point~$(x_0,t_0)$ in spacetime is in the {\em singular set}\/ of~$v$, which we will denote by~$\mathscr{S}^\infty_v$, if $v\not\in L^\infty((t_0-\delta,t_0+\delta)\times B)$ for any ball $B\ni x_0$ and any $\delta>0$. More generally, the {\em $q$-singular set}\/ of~$v$, $\mathscr S_v^q$, consists of the spacetime points $(x_0,t_0)$ such that $v\not\in L^\infty((t_0-\delta,t_0+\delta), L^q( B))$ for any ball~$B$ and any $\delta>0$ as above. Clearly $\mathscr S_v^{q}\subset \mathscr S_v^{q'}$ if $q<q'$ and $\mathscr S_v^q$ is a closed set. 

We are now ready to state the result. Basically, the theorem says that, given any smooth solution~$v_0$ on $\Omega\times(0,T)$ and any open set $U\subset\Omega$, there is an admissible weak solution~$v$ which coincides with $v_0$ outside~$U$ and which is uniformly close to~$v_0$ at time~$0$, yet blows up dramatically on a subset of~$U\times (0,T)$ of full dimension.
Interestingly,  smooth stationary Euler flows with compact support~\cite{Gavrilov,Vicolcptsupp} are very useful as building blocks in the construction of these solutions.

\begin{theorem}
	\label{teorema blowup}
	Consider some $0<\beta<1/3$ and some $q>2$. Let $T>0$ and let $\Omega$ be $\TT^3$ or an open subset of $\RR^3$. Fix some open set~$U$ whose closure is contained in~$\Omega$. Let $v_0$  be a smooth solution of the Euler equations in $\Omega\times(0,T)$. For any $\varepsilon>0$ there exists a weak solution $v\in L^2(\Omega\times(0,T))$ of the Euler equations whose $q$-singular set $\mathscr{S}_v^q$ is contained in $ U\times(0,T]$ and has Hausdorff dimension $4$. Furthermore, $v$~coincides with~$v_0$ on $(\Omega\backslash U)\times[0,T]$ and  satisfies
 \[\norm{v(\cdot,0)-v_0(\cdot,0)}_{C^0(\Omega)}<\varepsilon.\]
 Moreover, $v\in C^\beta_\text{loc}((\Omega\times[0,T])\backslash \mathscr{S}_v^q)$ and the energy profile $\int_\Omega \abs{v}^2\,dx$ can be chosen to be nonincreasing. 
\end{theorem}

\subsection{Strategy of the proof}

We prove Theorem~\ref{teorema gordo} in two stages: first we extend the field to $\RR^3\times[0,T]$ as a smooth subsolution (see Definition~\ref{D.subsolution} in the main text), and then we use a Nash iteration to perturb it into a weak solution. These stages are interrelated in that  tools and ideas  that we  develop to manipulate subsolutions also play a fundamental role in our convex integration scheme.

Concerning the extension of a local smooth solution of the Euler equations as a subsolution, the key result we prove is the following. In view of future applications of this result, which will appear elsewhere, we are stating these results for the Euler equations in any spatial dimension $n\geq 2$. 

\begin{theorem}
	\label{extender subs con soporte compacto}
	Let $\Omega_0\subset \mathbb{R}^\dim,\;\dim\geq 2$ be a bounded open set with smooth boundary and finitely many connected components and let $I\subset \mathbb{R}$ be a closed and bounded interval. Let $(v_0,p_0,\R{R}_0)\in C^\infty(\overline{\Omega}_0\times I)$ be a subsolution in $\Omega_0\times I$. Let $\Omega$ be a neighborhood of $\overline{\Omega}_0$. There exists a subsolution $(v,p,\R{R})\in C^\infty(\RR^\dim\times I)$ that extends $(v_0,p_0,\R{R}_0)$ and such that $\supp (v,p,\R{R})(\cdot,t)\subset \Omega$ for all $t\in I$ if and only if for each connected component $\Sigma$ of $\partial \Omega_0$ and all times $t\in I$ we have
	\begin{align*}
		&\int_{\Sigma}v_0\cdot \normal=\int_{\Sigma}\left[(a\cdot x)\partial_t v_0+(a\cdot v_0)v_0+p_0 a-a^t\mathring{R}_0\right]\cdot \normal= 0  \qquad \forall a\in\RR^\dim. 
	\end{align*}
	These conditions are automatically satisfied if $\partial \Omega_0$ is connected.
\end{theorem}

Regarding the convex integration scheme, we start off with the strategy from~\cite{Onsager_final}, which we implement in the context of solutions with compact support. The main issue we have to address is that, as we want the resulting solution to coincide with $v_0$ in~$\Omega\times[0,T]$, we must ensure that the scheme does not modify the subsolution in that region.

The result of our construction is:
\begin{theorem}
	\label{teorema integración convexa}
	Fix some $T>0$ and let $\Omega\subset \RR^3$ be a bounded open set with smooth boundary and with a finite number of connected components. Let $(v_0,p_0,\R{R}_0)\in C^\infty(\RR^3\times[0,T])$ be a subsolution of the Euler equations. Suppose that $\supp\R{R}_0\subset \overline{\Omega}\times[0,T]$. Let $0<\beta<1/3$ and let $e\in C^\infty([0,T],(0,\infty))$ be an energy profile such that 
	\begin{equation}
		\label{condición energía 1}
		e(t)>\int_\Omega\abs{v_0(x,t)}^2dx+6\|\R{R}_0\|_{L^\infty}\abs{\Omega}
	\end{equation}
 for all $0\leq t\leq T$.
	Then, there exists a weak solution of the Euler equations, $v\in C_c^\beta(\RR^3\times[0,\infty))$,  such that $v=v_0$ in $(\RR^3\backslash \Omega)\times[0,T]$ and
$$	
 \int_{\Omega}\abs{v(x,t)}^2dx=e(t)
$$
for all~$t\in[0,T]$.
\end{theorem}

\subsection{Organization of the paper}
In \cref{section matrices} we develop a set of tools to handle the construction, extension and gluing of subsolutions that will be used throughout the paper; in particular we prove Theorem~\ref{extender subs con soporte compacto}. In \cref{section overview} we present the iterative process used to prove \cref{teorema integración convexa}, which is carried out in a number of stages. The technical details of each stage of the construction are discussed in detail in Sections~\ref{proof of prop} to~\ref{section perturbation}. The very short \cref{section small} shows how to pass from Theorems~\ref{extender subs con soporte compacto} and~\ref{teorema integración convexa} to \cref{teorema gordo} and Remark~\ref{R.teoremagordo}. The modification of this scheme to account for the lack of uniform-in-time bounds is carried out in \cref{section intervalo abierto}. The applications concerning vortex sheets and blowup, cf. \cref{teorema vortex sheet} and \cref{teorema blowup}, are discussed in Sections~\ref{section vortex sheet} and~\ref{section blowup}, respectively. The paper concludes with two appendices, one about H\"older norms and another with some auxility estimates.


\section{Extension of subsolutions}\label{section matrices}

The goal of this section is to prove the extension theorem of smooth subsolutions stated in Theorem~\ref{extender subs con soporte compacto}. This is a key ingredient to prove our main theorem on the extension of weak solutions of the Euler equations. In Subsection~\ref{SS.GS} we sketch the strategy to prove Theorem~\ref{extender subs con soporte compacto}. Some instrumental tools from Hodge theory are presented in Subsection~\ref{SS.BT}, and in Subsection~\ref{SS.div} we show how to construct compactly supported solutions to the key divergence equation. Finally, the proof of Theorem~\ref{extender subs con soporte compacto} is completed in Subsection~\ref{SS.sub}.

In addition to constructing the desired solutions, we must estimate their derivatives. We refer to \cref{A:Holder} for the definition of the Hölder norms used. Specifically, we warn the reader that when dealing with time-dependent functions, we consider the supremum in time of the corresponding norms in space. Nevertheless, obtaining bounds on the derivatives of the solutions to certain differential equations is not enough for our construction. As we will see, we also need to control their $C^0$ norm. This can be achieved if we work with Besov spaces, which are defined in \cref{A:Holder}.

Throughout this section, we denote the space of $n\times n$ symmetric matrices as $\mathcal{S}^n$ and the space of $n\times n$ skew-symmetric matrices as $\mathcal{A}^n$. Unless otherwise stated, the dimension is $n\geq 2$. We define the divergence of a matrix $M\in C^\infty(\mathbb{R}^n\times\mathbb{R},\mathbb{R}^{n\times n})$ as the vector field whose coordinates are given by \[(\Div \, M)_i\coloneqq\sum_{j=1}^n\partial_j M_{ij},\]
where the derivatives are taken only with respect to the spatial variables. More generally, partial derivatives with Latin subscripts denote partial derivatives in the spatial coordinates, whereas temporal partial derivatives are always denoted by $\partial_t$.

We will repeatedly use Einstein's summation convention: when an index appears twice in an expression, it is implicitly summed over its range. Indices that appear only once in an expression are free indices and are not summed over.

Let us now recall the definition of subsolution of the Euler equations:

\begin{definition}\label{D.subsolution}
Let $V\subset \mathbb{R}^n\times \mathbb{R}$ be an open set. We will say that a triplet $(v,p,\R{R})\in C^\infty(V, \mathbb{R}^n\times \mathbb{R}\times \mathcal{S}^n)$ is a \emph{subsolution} of the Euler equations if
	\begin{equation}
		\label{defin weak sol}
		\begin{cases} \partial_t v+v\cdot\nabla v+\nabla p=\Div \R{R}, \\ \Div v=0.\end{cases}
	\end{equation}
	The symmetric matrix $\R{R}$ is known as the \emph{Reynolds stress} and it measures the deviation from being a solution of the Euler equations. It is customary to also impose that \begin{equation}
	    \tr \R{R}=0.
	\end{equation}
 All along this article, $\circ$ above a symmetric matrix will indicate that it is trace-free.
\end{definition}

Finally, let us fix some notation that will be used all along this section. We introduce the following norms in the space of $n\times n$ matrices:
	\begin{equation}
		\label{def norma matrices}
		\abs{M}\coloneqq \max_{\zeta\in\mathbb{S}^{n-1}}\abs{M\zeta}, \qquad \abs{M}_2\coloneqq \left(\sum_{i,j=1}^n M_{ij}^2\right)^{1/2}.
	\end{equation}
	Unless otherwise stated, we will always use the operator norm $\abs{\cdot}$. However, in some parts of the article we will exploit the elementary property that $\abs{\cdot}_2^2$ depends smoothly on the matrix entries. Note that $\abs{M}_2^2=\text{tr}(M^tM)$, which is invariant under orthogonal transformations. Hence, in the case of a symmetric matrix $S\in \mathcal{S}^n$, we have
	\begin{equation}
		\label{inequality norm matrices}
		\abs{S}\leq \abs{S}_2,
	\end{equation}
	which can be easily deduced by using an orthonormal basis of eigenvectors.

\subsection{General strategy}\label{SS.GS}
Our techniques for extending subsolutions and performing convex integration in the nonperiodic setting rely on obtaining compactly supported solutions to the (matrix) divergence equation when the source term is compactly supported. Let us illustrate the key ideas behind our method with the following toy problem: 

\begin{problem}
	Given $\rho\in C^\infty_c(\RR^3)$ such that $\int\rho=0$, find $v\in C^\infty_c(\RR^3,\RR^3)$ such that $\Div v=\rho$.
\end{problem}

It is easy to see that $v_0=\nabla\,\Delta^{-1}\rho$ solves our equation, but in general it is not compactly supported. To fix this, let $B$ be a ball containing the support of $\rho$, so that $v_0$ is divergence-free outside $B$. In addition, it follows from the divergence theorem that
\[0=\int_B \rho=\int_B \Div v_0=\int_{\partial B}v_0\cdot \normal.\]
This ensures that in $\RR^3\backslash B$ the divergence-free field $v_0$ can be written as $v_0=\curl w_0$ for some smooth field $w_0$. We extend $w_0$ to a smooth field $w\in C^\infty(\RR^3,\RR^3)$ and we define
\[v\coloneqq v_0-\curl w.\]
Since $w$ extends $w_0$, we see that $v$ vanishes outside $B$. Furthermore, $\Div v=\rho$ because $\Div\curl\equiv 0$. Therefore, $v$ is the sought field, which is clearly not unique.
	
Our approach to solving the divergence equation in the matrix case is the same: the potential-theoretic solution of the equation is not compactly supported. However, far from the support of the source  our matrix will be the image of certain differential operator $\mathscr{L}$ applied to a smooth potential, which is in the kernel of the divergence. We will extend the potential to the whole space and then subtract it from the potential-theoretic solution, obtaining a compactly supported solution. 

Just like in the vector case, we will have to impose certain integrability conditions on the source term for this to be possible. As we will see, these conditions are related to the classical conservation laws of linear and angular momentum in the Euler equations.

A totally different method to construct compactly supported solutions to the divergence equation in the (symmetric) matrix case was developed by Isett and Oh in~\cite{Isett}. Their theorem is stated in a very different setting and adapting it to what we need would require certain work. On the other hand, it will be relatively easy to deduce our result as a consequence of our analysis of the operator $\mathscr{L}$ introduced below, which is necessary for our result on the extension of subsolutions. Hence, we have preferred to take this path, which we believe is simpler (partly because it has a nice interpretation in terms of the elementary operations of vector calculus).
	
\subsection{Basic tools}\label{SS.BT}
The tools that we will need come from the Hodge decomposition theorem for manifolds with boundary. A good reference for this topic is~\cite{schwarz}. Nevertheless, we do not need the full generality of these results, as we will work in bounded domains of $\mathbb{R}^\dim$. Let us summarize the notation and main definitions that we will need.

Let $\Omega\subset \mathbb{R}^\dim$ be a bounded domain with smooth boundary. We denote by $\Lambda^k$ the vector space of skew-symmetric $k$-forms over $\mathbb{R}^\dim$, for $0\leq k\leq \dim$. In this setting, differential $k$-forms are maps $\omega\in C^\infty(\overline{\Omega},\Lambda^k)$. They form vector spaces in which we have two differential operators: the exterior derivative $d:C^\infty(\overline{\Omega},\Lambda^k)\to C^\infty(\overline{\Omega},\Lambda^{k+1})$ and the codifferential $\delta:C^\infty(\overline{\Omega},\Lambda^{k+1})\to C^\infty(\overline{\Omega},\Lambda^k)$. The Euclidean product induces a natural scalar product $(\cdot,\cdot )$ in $C^\infty(\overline{\Omega},\Lambda^k)$. The tangential part of a differential form is $\mathbf{t}\omega\coloneqq j^*\omega$, where $j:\partial \Omega\hookrightarrow\overline{\Omega}$ is the natural inclusion, and $j^*$ is the pushforward. We define the Dirichlet harmonic $k$-forms as:
\[\mathcal{H}^k_D(\overline{\Omega})\coloneqq\{\omega\in C^\infty(\overline{\Omega},\Lambda^k):d\omega=0,\delta\omega=0,\mathbf{t}\omega=0\}.\]
Finally, to obtain quantitative estimates we will need to work in Hölder spaces. We refer to Appendix~\ref{A:Holder} for the definition of these norms. We also recommend to take a look at the appendix to check our convention of H\"older norms when the field is time-dependent.

With this notation, the first basic lemma that we shall use is:
\begin{lemma}
	\label{schwarz no t}
	Let $\Omega\subset \mathbb{R}^\dim$ be a bounded domain with smooth boundary and let $\rho\in C^\infty(\overline{\Omega},\Lambda^k)$. The boundary value problem
	\[(P)\begin{cases}\delta \omega=\rho,\\d\omega=0,\\ \mathbf{t}\omega=0,\\ ( \omega, \lambda) =0 \quad \forall \lambda\in \mathcal{H}^{k+1}_D(\overline{\Omega})\end{cases}\]
	is solvable if and only if
	\[\delta \rho=0\qquad \text{and}\qquad \int_{C_{\dim-k}}\hodge \rho=0\quad \forall\; (\dim-k)\text{-cycle }C_{\dim-k}.\]
	In that case, the solution is unique and we have\[\norm{\omega}_{C^{N+1+\alpha}(\Omega)}\leq C\norm{\rho}_{C^{N+\alpha}(\Omega)}\]
	for any $N\geq 0$, $\alpha\in (0,1)$ and certain constants $C\equiv C(N,\alpha,\Omega)$.
\end{lemma}
For a proof, see \cite[Theorem 7.2]{Dacorogna} and \cite[Corollary 3.2.4]{schwarz}. In \cite[Theorem 3.2.5]{schwarz} we can find the general case of a Riemannian manifold with boundary, but it does not include estimates for Hölder norms, only for Sobolev norms. We recall that, as usual, $\hodge$ is the Hodge star operator acting on differential forms and an $(\dim-k)$-cycle is an $(\dim-k)$-chain (in the sense of algebraic topology) whose boundary is zero.

We are mainly interested in the problem $\delta \omega=\rho$, but we have to add the other conditions to select a single solution. This allows us to obtain a time-dependent version of Lemma~\ref{schwarz no t}:
\begin{lemma}
	\label{schwarz}
	Let $\Omega\subset \mathbb{R}^\dim$ be a bounded domain with smooth boundary and let $I\subset \mathbb{R}$ be a closed and bounded interval. Let $\rho\in C^\infty(\overline{\Omega}\times I,\Lambda^k)$. There exists a differential form $\omega\in C^\infty(\overline{\Omega}\times I,\Lambda^{k+1})$ solving the boundary value problem $(P)$ at each $t\in I$ if and only if
	\[\delta \rho=0\qquad \text{and}\qquad \int_{C_{\dim-k}}\hodge \rho=0\quad \forall\; (\dim-k)\text{-cycle }C_{\dim-k}, \forall t\in I.\]
	In that case, the solution is unique and we have\[\norm{\omega}_{N+1+\alpha}\leq C\norm{\rho}_{N+\alpha}\]
	for any $N\geq 0$, $\alpha\in (0,1)$ and certain constants $C\equiv C(N,\alpha,\Omega)$.
\end{lemma}
\begin{proof}
	Given a time-dependent differential form, we denote by a subscript the differential form at a given time. By Lemma~\ref{schwarz no t}, the necessity of the conditions is clear. To prove that they are also sufficient, let us suppose that $\rho_t$ satisfies the conditions of Lemma~\ref{schwarz no t} at all times $t\in I$. Hence, applying Lemma~\ref{schwarz no t} at each time, we see that there exists a time-dependent $(k+1)$-form $\omega$ solving $(P)$ at each $t\in I$. The question is whether $\omega$ depends smoothly on $t$.
	
	Since $\partial_t \rho$ also satisfies the hypotheses of Lemma~\ref{schwarz no t} at all times $t\in I$, there exists a $(k+1)$-form $\widetilde{\omega}$ solving $(P)$ with data $\partial_t \rho$. For a fixed $t_0\in I$ and $h\neq 0$ small we see that $h^{-1}(\omega_{t_0+h}-\omega_{t_0})-\widetilde{\omega}_{t_0}$ is the unique solution of $(P)$ with data $h^{-1}(\rho_{t_0+h}-\rho_{t_0})-(\partial_t \rho)_{t_0}$. Therefore,
	\[\norm{\frac{\omega_{t_0+h}-\omega_{t_0}}{h}-\widetilde{\omega}_{t_0}}_{k+1+\alpha}\leq C\norm{\frac{\rho_{t_0+h}-\rho_{t_0}}{h}-(\partial_t \rho)_{t_0}}_{k+\alpha}\longrightarrow 0 \quad \text{as }h\longrightarrow 0.\]
	We deduce that $\widetilde{\omega}$ is the partial derivative with respect to time of $\omega$. Iterating this argument, we conclude that $\omega$ depends smoothly on time. The claimed estimates are easily obtained by taking the supremum on $t\in I$ in the estimates of Lemma~\ref{schwarz no t}.
\end{proof}

\begin{remark} 
	\label{solo homología}
	If $\delta \rho=0$, the integral of $\hodge \rho$ on an $(\dim-k)$-cycle depends only on the homology class of the cycle. Indeed, if $C$ and $C'$ are two $(\dim-k)$-cycles in $\overline{\Omega}$ that are the boundary of an $(\dim-k+1)$-chain, by Stokes' theorem we have
	\[\int_{C'}\hodge \rho-\int_{C}\hodge \rho=\int_{\partial \mathcal{N}}\hodge \rho=\int_{\mathcal{N}}d\hodge \rho=(-1)^k\int_\mathcal{N}\hodge \delta \rho=0.\]
\end{remark}

The machinery of differential geometry is quite powerful, but we are interested in the simpler setting of bounded domains $\Omega\subset\mathbb{R}^\dim$. Taking advantage of the canonical basis of Euclidean space, we may forget about differential forms and work with simpler objects. Indeed, there is a natural correspondence between 1-forms and vector fields and between 2-forms and skew-symmetric matrices $\mathcal{A}^\dim$:
\begin{align*}
&C^\infty\left(\overline{\Omega},\Lambda^1\right)\to C^\infty\left(\overline{\Omega},\RR^\dim\right), \qquad \sum_{i,j=1}^na_{i}\,dx_i\hspace{44pt}\mapsto\;\, \sum_{i=1}^na_{i}\hspace{1pt}e_i, \\
&C^\infty\left(\overline{\Omega},\Lambda^2\right)\to C^\infty\left(\overline{\Omega},\mathcal{A}^\dim\right), \qquad \frac12\sum_{i,j=1}^na_{ij}\,dx_i\wedge dx_j\;\,\mapsto\;\, \sum_{i,j=1}^na_{ij}e_i\otimes e_j.
\end{align*}
Here we have used that $dx_i\wedge dx_j=-dx_j\wedge dx_i$. Using the canonical base of $1$-forms, the action of the codifferential can be summarized as
\[\delta(f\hspace{1pt}dx_i)=\partial_if, \hspace{50pt} \delta(f\hspace{1pt}dx_i\wedge dx_j)=-\partial_jf\hspace{1pt}dx_i,\]
where $f$ is any smooth function and $1\leq i<j\leq \dim$. One can then check that the following diagram commutes:
\begin{center}
	\begin{tikzcd}
		C^\infty\left(\overline{\Omega},\Lambda^2\right) \arrow[d,leftrightarrow] \arrow[r, "\delta"] & C^\infty\left(\overline{\Omega},\Lambda^1\right) \arrow[r, "\delta"] \arrow[d,leftrightarrow]&C^\infty\left(\overline{\Omega},\Lambda^0\right) \arrow[d,equal]\\
		C^\infty\left(\overline{\Omega},\mathcal{A}^\dim\right) \arrow[r, "div"]& C^\infty\left(\overline{\Omega},\mathbb{R}^\dim\right) \arrow[r, "div"] & C^\infty\left(\overline{\Omega}\right)
	\end{tikzcd}
\end{center}
This allows us to write everything in terms of matrices and to simplify the notation. Using this correspondence and Remark~\ref{solo homología}, we may formulate a particular case of Lemma~\ref{schwarz} as follows:
\begin{lemma}
	\label{invertir matrices}
	Let $\Omega\subset \RR^\dim$ be a bounded domain with smooth boundary and let $I\subset \mathbb{R}$ be a closed and bounded interval. Let $v\in C^\infty(\overline{\Omega}\times I,\RR^\dim)$. The following are equivalent:
	\begin{enumerate}
		\item there exists $A\in C^\infty(\overline{\Omega}\times I,\mathcal{A}^\dim)$ such that $\Div A=v$,
		\item $\Div v=0$ and $\int_\Sigma v\cdot \nu=0$ for any connected component $\Sigma$ of $\partial \Omega$ and any fixed $t\in I$,
		\item $\int_{C_{\dim-1}} v\cdot \nu=0$ for any $(\dim-1)$-cycle $C_{\dim-1}$ and any $t\in I$.
	\end{enumerate}
	In that case, $A\in C^\infty(\overline{\Omega}\times I,\mathcal{A}^\dim)$ may be chosen so that 
	\[\norm{A}_{N+1+\alpha}\leq C\norm{v}_{N+\alpha}\]
	for any $N\geq 0$, $\alpha\in (0,1)$ and certain constants $C\equiv C(N,\alpha,\Omega)$.
\end{lemma}
The proof is straightforward taking into account that the codifferential $\delta$ becomes the operator $\Div$ and the integral on cycles becomes the flux of the corresponding vector field across a closed surface. By \cref{solo homología}, the integral only depends on the homology class, so we can choose to integrate on the connected component of $\partial\Omega$ belonging to each class.

\subsection{The divergence equation}\label{SS.div}
After collecting some basic tools from Hodge theory in the previous subsection, we will now show how to obtain compactly supported solutions to the divergence equation. We begin by introducing some potential-theoretic solutions, which we will later modify in order to fix the support. Let us consider the following differential operator that maps smooth vector fields (with bounded derivatives) to $C^\infty(\RR^\dim,\mathcal{S}^\dim)$: 
\begin{equation}
\label{def mathscr R}
	(\mathscr{R}f)_{ij}\coloneqq \Delta^{-1}(\partial_i f_j+\partial_j f_i)-\delta_{ij}\Delta^{-1}\Div f.
\end{equation}
Here $\Delta^{-1}$ refers to the potential-theoretic solution of the Poisson equation, that is, the (spatial) convolution of the source term  with
the fundamental solution of the Laplace equation in $\mathbb R^n$. We remind the reader that partial derivatives with latin subscripts denote partial derivatives in the spatial coordinates, whereas temporal partial derivatives are always denoted by $\partial_t$.

A direct calculation shows that $\Div\mathscr{R}f=f$. 
We notice that $\mathscr{R}$ is not trace-free. This is not an issue in our proofs, because our constructions with potentials do not preserve being trace-free, so we will usually absorb the trace into the pressure at the end.

Let us now derive a very useful identity. Let $\Omega\subset \RR^\dim$ be a bounded open set with smooth boundary. Let $v\in C^\infty(\overline{\Omega},\RR^\dim)$ and $S\in C^\infty(\overline{\Omega},\mathcal{S}^\dim)$. Integrating by parts, we have
\begin{equation}
	\label{green id}
	\int_\Omega v\cdot \Div S+\int_\Omega(\gradsim v):S=\int_{\partial \Omega}v^t S \,\nu,
\end{equation}
where $\normal$ is the unitary normal vector associated to the exterior orientation and the operator $\gradsim$ is given by:
\[\gradsim:C^\infty(\Omega,\mathbb{R}^\dim)\to C^\infty(\Omega,\mathcal{S}^\dim), v\mapsto \frac{1}{2}(\nabla v+\nabla v^t).\] Its kernel are the so-called Killing vector fields. It is a finite-dimensional vector space that plays an important role in Riemannian geometry. It is well known (see \cite[page 52]{Petersen}) that in $\mathbb{R}^\dim$ a basis of this vector space is given by
\begin{equation}
	\label{definition base Killing}
	\mathcal{B}\coloneqq \{e_1, \dots, e_\dim, \xi_{12}, \dots , \xi_{(\dim-1)\dim},\}
\end{equation}
where
\begin{equation}
	\label{definition Killing momento angular}
	\xi_{ij}\coloneqq x_i \hspace{1pt}e_j-x_j \hspace{1pt} e_i, \qquad 1\leq i<j\leq \dim.
\end{equation}

Next, we introduce two vector spaces and a differential operator that will be very important in our construction:
\begin{definition} Let $\Omega\subset \mathbb{R}^\dim$ be a bounded domain an let $I\subset \mathbb{R}$ be a closed and bounded interval. We define two vector spaces:
	\begin{align*}
		\mathcal{P}(\overline{\Omega}\times I)&\coloneq \left\{A\in C^\infty(\overline{\Omega}\times I,\RR^{\dim^4}):A^{ik}_{jl}=-A^{ki}_{jl},\, A^{ik}_{jl}=-A^{ik}_{lj}\right\},\\
		\mathcal{G}(\overline{\Omega}\times I)&\coloneq \left\{S\in C^\infty(\overline{\Omega}\times I,\mathcal{S}^\dim)\,:\,\Div S=0,\; \int_\Sigma \xi^t S\, \normal=0\hspace{10pt} \begin{matrix}\forall \xi \in\text{ker}\,\gradsim, \\ \forall \Sigma \text{ comp. of }\partial \Omega, \\ \forall t\in I.\end{matrix}\,\right\}
	\end{align*}
	and we consider the differential operator
	\[\mathscr{L}:\mathcal{P}(\overline{\Omega}\times I)\to C^\infty(\overline{\Omega}\times I,\RR^{\dim^2}), \;\; \left[\mathscr{L}(A)\right]_{ij}=\frac{1}{2}\sum_{k,l}\partial_{kl}\left(A^{ik}_{jl}+A^{jk}_{il}\right).\]
\end{definition}
This operator already appeared in the context of convex integration in the original article by De Lellis and Sz\'ekelyhidi~\cite{Bounded}, who noticed that the image of the operator $\mathscr{L}$ is contained in the space of divergence-free matrices.

For our purposes, this operator can be regarded as a matrix analog of the $\curl$ operator in arbitrary dimension. In order to perform the construction sketched in Subsection~\ref{SS.GS}, the next step is to understand how to invert this operator (under the appropriate boundary conditions).  The following lemma is the key of our approach to solve the divergence equation:

\begin{lemma}
	\label{key lemma}
	Let $\Omega\subset \RR^\dim$ be a bounded domain with smooth boundary and let $I\subset \mathbb{R}$ be a closed and bounded interval. Then $\mathcal{G}(\overline{\Omega}\times I)$ is the image of the differential operator $\mathscr{L}$. Furthermore, given $S\in\mathcal{G}(\overline{\Omega}\times I)$, there exists $A\in \mathcal{P}(\overline{\Omega}\times I)$ such that $S=\mathscr{L}(A)$ and for any $\alpha \in (0,1)$ we have
	\[\norm{A}_{N+2+\alpha}\leq C\norm{S}_{N+\alpha}\]
	for all $N\geq 0$ and certain constants $C\equiv C(N,\alpha,\Omega)$.
\end{lemma}
\begin{proof}
	First, we prove that the image of $\mathscr{L}$ is contained in $\mathcal{G}(\overline{\Omega}\times I)$. Fix an arbitrary $A\in\mathcal{P}(\overline{\Omega}\times I)$. It is clear from the definition that $\mathscr{L}(A)$ is symmetric. The fact that it is divergence-free follows from the skew-symmetric properties of $A$:
	\begin{align*}
		[\Div \mathscr{L}(A)]_i&=\frac{1}{2}\sum_{j,k,l}\partial_{jkl}\left(A^{ik}_{jl}+A^{jk}_{il}\right)=\\
		&=\sum_k\frac{1}{2}\,\partial_k\left(\sum_{j,l}\partial_{jl}A^{ik}_{jl}\right)+\sum_l\frac{1}{2}\,\partial_l\left(\sum_{j,k}\partial_{jk}A^{jk}_{il}\right)=0.
	\end{align*}
	Next, we fix and arbitrary Killing vector field $\xi$ and a connected component $\Sigma$ of $\partial \Omega$. We choose a smooth cut-off function $\varphi$ that vanishes in a neighborhood of $\Sigma$ and it is identically 1 in a neighborhood of the other connected components of $\partial \Omega$. By the choice of $\varphi$ we have
	\[\int_{\partial \Omega} \xi^t\mathscr{L}(\varphi A)\,\normal=\int_{\partial \Omega}\xi^t\mathscr{L}(A)\,\normal-\int_\Sigma \xi^t\mathscr{L}(A)\,\normal.\]
	By our previous discussion, $\mathscr{L}(A)$ and $\mathscr{L}(\varphi A)$ are symmetric and divergence-free. In addition, $\xi$ is a Killing vector, so $\gradsim w=0$. Thus, from (\ref{green id}) we deduce that the term on the left-hand side of the previous equation vanishes and so does the first term on the right-hand side. Therefore, we see that $\int_\Sigma \xi^t\mathscr{L}(A)\,n=0$ and, since $A$, $\xi$ and $\Sigma$ are arbitrary, we conclude that the image of $\mathscr{L}$ is contained in $\mathcal{G}(\overline{\Omega}\times I)$. 
	
	Now we will prove the other inclusion and the stated estimate. We fix an arbitrary $S\in \mathcal{G}(\overline{\Omega}\times I)$. By definition, when choosing the canonical basis of $\mathbb{R}^\dim$ as Killing vectors , we obtain
	\[\int_\Sigma S_{ij}\normal_j=0\qquad \forall \Sigma \text{ connected component of }\partial \Omega\]
	for any $i=1, \dots, m$. If we fix $i$, we may apply Lemma~\ref{invertir matrices} to conclude that there exists $B^i\equiv B^i_{jl}\in C^\infty(\overline{\Omega}\times I,\mathcal{A}^\dim)$ such that $\partial_l B^i_{jl}=S_{ij}$. 
	
	Next we fix a Killing field $\xi$ of the form $\xi_i=R_{ik}x_k$, where $R\in \mathcal{A}^\dim$. For $j=1, \dots m$ we compute
	\[\partial_l(\xi_i B^i_{jl})=\xi_i\partial_l B^i_{jl}+(\partial_l \xi_i)B^i_{jl}=\xi_i S_{ij}+R_{il} B^i_{jl}.\]
	Note that, since $S\in \mathcal{G}(\overline{\Omega}\times I)$,
	\[\int_\Sigma \xi_i S_{ij} \normal_j=0 \qquad \forall \Sigma \text{ connected component of }\partial \Omega, \quad \forall t\in I.\]
	Regarding the left-hand side term, we define the forms
	\begin{align*}
		\omega &\coloneq \sum_{i,j,l} \partial_{l}(\xi_i B^i_{jl})\,dx_j,\\
		\eta &\coloneq \sum_{j<l, i} \xi_i B^i_{jl}\,dx_j\wedge dx_l.
	\end{align*}
	Since $B^i_{jl}$ is skew-symmetric in the lower indices, we see that $\delta \eta=\omega$. Using the properties of the Hodge star operator and the codifferential, we have $\hodge\delta \eta=d\hodge\eta$. These forms allow us to rewrite the integral on an $(\dim-1)$-cycle at any $t\in I$ as:
	\[\int_{C_{\dim-1}}\partial_l(\xi_iB^i_{jl})\,\normal_j=\int_{C_{\dim-1}}\omega(n)\,\tilde{\mu}=\int_{C_{\dim-1}}\hodge \,\omega=\int_{C_{\dim-1}}d(\hodge\, \eta)=0,\]
	where $\tilde{\mu}$ is the measure induced by the standard measure in $\mathbb{R}^\dim$. We have used Stokes' theorem and the fact that $(\dim-1)$-cycles have no boundary. We conclude that for any $(\dim-1)$-cycle
	\[\int_{C_{\dim-1}} R_{il}B^i_{jl}\normal_j=0.\]
	Choosing $R=e_{i_0}\otimes e_{l_0}-e_{l_0}\otimes e_{i_0}$, that is, choosing $\xi$ as $\xi_{l_0i_0}$, we see that for any $i,l=1, \dots m$ we have:
	\[\int_{C_{\dim-1}} (B^i_{jl}-B^l_{ji})\,\normal_j=0 \qquad \text{for any }(\dim-1)\text{-cycle }C_{\dim-1} \text{ and all }t\in I.\]
	Applying again Lemma~\ref{invertir matrices}, we obtain $A^{ik}_{jl}$ skew-symmetric in $j, k$ such that
	\[\partial_k A^{ik}_{jl}=B^i_{jl}-B^l_{ji}.\]
	Therefore,
	\begin{align*}
		\frac{1}{2}\partial_{kl}(A^{ik}_{jl}+A^{jk}_{il})&=\frac{1}{2}\partial_l\left[\partial_k A^{ik}_{jl}+\partial_k A^{jk}_{il}\right]=\frac{1}{2}\partial_l\left[(B^i_{jl}-B^l_{ji})+(B^j_{il}-B^l_{ij})\right]\\&=\frac{1}{2}(\partial_l B^i_{jl} +\partial_l B^j_{il})=S_{ij},
	\end{align*}
	where we have used that $B$ is skew-symmetric in the lower indices and the symmetry of $S$: $\partial_l B^i_{jl}=S_{ij}=S_{ji}=\partial_l B^j_{il}$. In summary, we have found $A\in C^\infty\left(\overline{\Omega}\times I,\mathbb{R}^{\dim^4}\right)$ such that:
	\begin{enumerate}[(i)]
		\item $A^{ik}_{jl}=-A^{ij}_{kl}$,
		\item $\partial_k A^{ik}_{jl}=-\partial_k A^{lk}_{ji}$,
		\item $\frac{1}{2}\partial_{kl} \left(A^{ik}_{jl}+A^{jk}_{il}\right)=S_{ij}$.
	\end{enumerate}
	We define
	\[\widetilde{A}^{ik}_{jl}\coloneq \frac{1}{2}\left(A^{il}_{jk}-A^{kl}_{ji}\right).\]
	It is clear that $\widetilde{A}^{ik}_{jl}$ is skew-symmetric in $i,k$. In addition, it is skew-symmetric in $j,l$ by (i). Furthermore,
	\[\partial_l \widetilde{A}^{ik}_{jl}=\frac{1}{2}\left(\partial_l A^{il}_{jk}-\partial_l A^{kl}_{ji}\right)\stackrel{\text{(ii)}}{=}\partial_l A^{il}_{jk}.\]
	Hence, using (iii) we conclude
	\[\frac{1}{2}\partial_{kl} \left(\widetilde{A}^{ik}_{jl}+\widetilde{A}^{jk}_{il}\right)=\frac{1}{2}\partial_{kl} \left(A^{ik}_{jl}+A^{jk}_{il}\right)=S_{ij}.\]
	Therefore, $\widetilde{A}\in\mathcal{P}(\overline{\Omega}\times I)$ and $\mathscr{L}(\widetilde{A})=S$, as we wanted. The estimates for $\widetilde{A}$ follow from applying twice the estimates from Lemma~\ref{invertir matrices}.
\end{proof}

Finally, we are ready to prove the main result of this subsection, which establishes the existence of compactly supported solutions to the divergence equation. In a different setting, a related class of compactly supported solutions to the divergence equation were constructed by Isett and Oh~\cite[Theorem 10.1]{Isett}. Our approach is based on the operator $\mathscr{L}$, which will be essential for the extension of subsolutions in \cref{extender fuera}. We observe that the compatibility conditions~\eqref{integrability condition divergence} in Lemma~\ref{invertir divergencia matrices} are precisely the conditions~(202) in~\cite{Isett}.

We recall that the Besov norms that we use are defined in \cref{A:Holder}. We need to work with these norms because they are necessary to derive estimates for the $C^\alpha$ norm of the resulting matrix. This will be essential in the proof of \cref{teorema integración convexa}.

\begin{lemma}
	\label{invertir divergencia matrices}
	Let $\Omega\subset \RR^\dim$ be a bounded domain with smooth boundary and let $I\subset \mathbb{R}$ be a closed and bounded interval. Let $f\in C^\infty(\RR^\dim\times I,\RR^\dim)$  such that $\supp f(\cdot,t)\subset \Omega$ for all $t\in I$. Then, there exists $S\in C^\infty(\RR^\dim\times I,\mathcal{S}^\dim)$ such that $\Div S=f$ and $\supp S(\cdot,t)\subset \Omega$ for all $t\in I$ if and only if
	\begin{equation}
		\label{integrability condition divergence}
		\int_\Omega f\cdot \xi=0\qquad \forall \xi\in \ker \gradsim, \;\forall t\in I.
	\end{equation}
	In that case, we may choose $S$ so that for all $N\geq0$ and any $\alpha\in (0,1)$ we have
	\[\norm{S}_{N+\alpha}\leq C\norm{f}_{B^{N-1+\alpha}_{\infty,\infty}}\]
	for certain constants $C=C(\Omega,N,\alpha)$.
\end{lemma}
\begin{proof}
	First of all, we show that the integrability condition (\ref{integrability condition divergence}) is necessary. Let us suppose that such an $S$ exists. We fix a ball $B\supset \overline{\Omega}$ and use the identity (\ref{green id}) to obtain
	\[0=\int_{\partial B}\xi^tS \normal=\int_{B}\xi\cdot \Div S=\int_{\Omega}\xi\cdot f \qquad \forall \xi\in\ker\gradsim,\;\forall t\in I.\]
	
	Let us show that condition (\ref{integrability condition divergence}) is also sufficient. The field $S_0\in C^\infty(\RR^\dim\times I,\mathcal{S}^\dim)$ given by $S_0\coloneqq \mathscr{R}f$ solves the equation $\Div S_0=f$, where $\mathscr{R}$ was defined in (\ref{def mathscr R}). It is easy to check that it satisfies the estimates
	\begin{equation}
		\label{estimates S0}
		\norm{S_0}_{N+\alpha}\leq C\norm{f}_{B^{N-1+\alpha}_{\infty,\infty}}
	\end{equation}
	for certain constants $C=C(N,\alpha)$ because $\mathscr{R}$ is an operator of order $-1$. However, it is not compactly supported, in general. We must modify it far from the support of $f$.
	
	We begin by studying the boundary conditions. Let $\Sigma_i$ be a connected component of $\partial \Omega$ and let $U_i$ be the domain bounded by it. We claim that
	\begin{equation}
		\label{integral source Ui}
		\int_{U_i}\xi\cdot f=0 \qquad \forall\xi\in\ker\gradsim,\;\forall t\in I.
	\end{equation}
	Indeed, since $\Omega$ is a bounded domain, $U_i$ must be either the complement of the unbounded connected component of $\RR^3\backslash \overline{\Omega}$ or one of the bounded connected components of $\RR^3\backslash \overline{\Omega}$ (if there are any). In the first case, (\ref{integral source Ui}) follows from the integrability condition (\ref{integrability condition divergence}) because $f(\cdot,t)\subset \Omega\subset U_i$. In the second case, (\ref{integral source Ui}) is trivial because $f(\cdot,t)$ vanishes on $U_i\subset \RR^3\backslash \overline{\Omega}$.
	Thus, applying the identity (\ref{green id}) to each $U_i$ we obtain
	\begin{equation}
		\label{boundary conditions S0}
		\int_{\Sigma_i}\xi^t S_0\normal=0 \qquad \forall\xi\in\ker\gradsim,\;\forall t\in I.
	\end{equation}
		
	Next, note that for sufficiently small $r>0$ the boundary of the open set
	\[G\coloneqq \{x\in \Omega:\dist(x,\partial \Omega)<r\}\]
	has twice as many connected components as $\partial \Omega$. Furthermore, the boundary of each connected component $G_i$ of $G$ consists of exactly two hypersurfaces, which we denote as $\Sigma_i$ and $\Sigma'_i$, and we have $\Sigma_i\subset \partial \Omega$. By further reducing $r>0$, we may assume that $f(\cdot,t)$ vanishes on $G$ at all times $t\in I$. Then, it follows from (\ref{boundary conditions S0}) and the identity (\ref{green id}) that
	\[\int_{\Sigma'_i}\xi^tS\normal=-\int_{\Sigma_i}\xi^tS\normal+\int_{\partial G_i}\xi^tS\normal=-\int_{\Sigma_i}\xi^tS\normal+\int_{G_i}\xi\cdot f=0\]
	for any $\xi\in \ker\gradsim$. Next, we fix a ball $B\supset \overline{\Omega}$ and we consider the domain \[U\coloneqq (B\backslash \Omega)\cup G.\] We see that \[\partial U=\partial B\cup\bigcup_i \Sigma_i'.\]
	Again, it follows from (\ref{green id}) and the integrability condition (\ref{integrability condition divergence}) that
	\[\int_{\partial B}\xi^tS_0\normal=\int_B\xi\cdot f=\int_{\Omega}\xi\cdot f=0 \qquad \forall\xi\in\ker\gradsim,\;\forall t\in I.\]
	We conclude that $S_0$ is divergence-free on $U$ and in each connected component $\Sigma$ of $\partial U$ we have
	\[\int_\Sigma \xi^tS_0\normal=0\qquad \forall\xi\in\ker\gradsim,\;\forall t\in I,\]
	that is, $S_0\in \mathcal{G}(\overline{U}\times I)$. By \cref{key lemma} there exists $A_0\in\mathcal{P}(\overline{U}\times I)$ such that $S_0(x,t)=\mathscr{L}(A_0)(x,t)$ for all $x\in \overline{G}$ and $t\in I$. Furthermore, for any $N\geq 0$ and $\alpha\in (0,1)$ we have
	\[\norm{A_0}_{N+2+\alpha}\leq C(U,N,\alpha)\norm{S_0}_{N+\alpha}\leq C(U,N,\alpha)\norm{f}_{B^{N-1+\alpha}_{\infty,\infty}}.\]
	The constants depend on $U$, which depends not only on the geometry of $\Omega$ but also on the minimum distance between the support of $f(\cdot,t)$ and $\partial \Omega$ through the parameter $r$. However, since $U$ tends to $B\backslash \overline{\Omega}$ as $r\to0$, the constants remain uniformly bounded, so they ultimately depend only on $\Omega$. For this, smoothness of $\partial\Omega$ is essential, as it allows us to choose parametrizations of $\overline{U}$ converging to parametrizations of $\overline{B}\backslash \Omega$ in any Hölder norm as $r\to 0$.
	
	Applying \cref{extension holder functions} and antisymmetrizing, we see that there exists a map $A\in C^\infty(\RR^\dim\times I,\RR^{\dim^4})$ such that $A^{ik}_{jl}=-A^{ki}_{jl},\, A^{ik}_{jl}=-A^{ik}_{lj}$ that extends $A_0$ outside $\overline{U}\times I$.
	Furthermore, for any $N\geq 0$ and $\alpha\in(0,1)$, we have
	\begin{equation}
		\label{estimates A}
		\norm{A}_{N+2+\alpha}\leq C(U,N)\norm{A_0}_{N+2+\alpha}\leq C(U,\Omega,N,\alpha)\norm{f}_{B^{N-1+\alpha}_{\infty,\infty}}.
	\end{equation}
    Again, since $U$ tends to $B\backslash\overline{\Omega}$ as $r\to 0$ in a suitable manner, the constants ultimately depend only on $\Omega$, $N$ and $\alpha$, since they will be uniformly bounded on $r\in(0,1)$.
	
	Finally, for $x\in B$ and $t\in I$ we define
	\[S\coloneqq S_0-\mathscr{L}(A).\]
	Since the image of $\mathscr{L}$ is contained in the kernel of the divergence, we see that $\Div S=f$. By construction $A$ extends $A_0$, so $\mathscr{L}(A)=\mathscr{L}(A_0)=S_0$ on $\overline{U}\times I$. Therefore, $S(\cdot,t)$ vanishes in $U$, so we may extend it by 0 to $\RR^\dim\times I$.
	
	In conclusion, we have constructed $S\in C^\infty(\RR^\dim\times I,\mathcal{S}^\dim)$ such that $\Div S=f$ and $\supp S(\cdot,t)\subset \Omega$ for all $t\in I$. Furthermore, the desired estimate follows from (\ref{estimates S0}) and (\ref{estimates A}) because $\mathscr{L}$ is a second order differential operator.
\end{proof}

\subsection{Subsolutions and proof of Theorem~\ref{extender subs con soporte compacto}}\label{SS.sub}
In this subsection we use \cref{key lemma} and \cref{invertir divergencia matrices} to glue and extend subsolutions, which will yield the proof of Theorem~\ref{extender subs con soporte compacto}. It should be apparent by now that controlling the $L^2$-product with the Killing fields is very important in these constructions. It is not difficult to construct $f\in C^\infty_c(\RR^\dim,\RR^\dim)$ with the desired $L^2$-product with the Killing fields. However, when working with subsolutions we will also need that $f$ be divergence-free. In addition, in our constructions we will work in domains of a certain form. Our approach is based on the following:
\begin{lemma}
	\label{controlar Killing}
Let $\Omega\subset \RR^\dim$ be a bounded domain with smooth boundary and let $I\subset \mathbb{R}$ be a closed and bounded interval. Let $r>0$ and let $L_{ij}\in C^\infty(I)$ for $1\leq i<j\leq\dim$. There exists a divergence-free field $w\in C^\infty(\RR^\dim\times I,\RR^\dim)$ such that the support of $w(\cdot,t)$ is contained in $\{x\in\RR^\dim:0<\dist(x,\Omega)<r\}$  and
\begin{align*}
	&\int a\cdot w\,dx=0, \qquad \int \xi_{ij}\cdot w\,dx=L_{ij}(t)
\end{align*}
for all $t\in I$ and $1\leq i<j\leq\dim$, where $\xi_{ij}$ is given by (\ref{definition Killing momento angular}). Furthermore, for any $N\geq 0$ we have
\begin{align*}
	&\norm{w}_{N}\leq C(N,n)\abs{\Omega}^{-1}\, r^{-(N+1)}\max_{ij,\, t\in I}\abs{L_{ij}(t)}, \\
	&\norm{\partial_t w}_{N}\leq C(N,n)\abs{\Omega}^{-1}\, r^{-(N+1)}\max_{ij,\, t\in I}\abs{L'_{ij}(t)}.
\end{align*}
\end{lemma}
\begin{proof}
We will construct our field as $w=\Div A$ for some $A\in C^\infty_c(\RR^\dim\times I,\mathcal{A}^\dim)$ that we will choose later. Since $A$ is compactly supported, it follows from the divergence theorem that $\int a\cdot w=0$ for any $a\in \RR^\dim$. Furthermore, for any $1\leq i<j\leq \dim$ we have
\begin{equation}
	\label{angular momentum vector potential}
	\int \xi_{ij}\cdot w=\int (\xi_{ij})_k\partial_l A_{kl}=-\int \partial_l(\xi_{ij})_k A_{kl}=-\int (A_{ji}-A_{ij})=2\int A_{ij}
\end{equation}
because $\partial_l(\xi_{ij})_k=\delta_{il}\delta_{jk}-\delta_{jl}\delta_{ik}$. Here we have denoted by $(\xi_{ij})_k$ the $k$-th component of the vector $\xi_{ij}$ and we have used Einstein's summation convention when summing over $k$ and $l$. By \cref{cutoff} we may choose a nonnegative cutoff function $\varphi\in C^\infty_c(\Omega+B(0,r))$ that is identically 1 in a neighborhood of $\Omega$ and such that 
\[\norm{\varphi}_N\leq C(N,n)\,r^{-N}.\]
We define
\[A(x,t)\coloneqq\sum_{1\leq i<j\leq n}L_{ij}(t)\,(e_i\otimes e_j-e_j\otimes e_i)\left(2\int \varphi\right)^{-1}\varphi(x).\]
Since $\varphi$ is constant in a neighborhood of $\Omega$ and its support is contained in $\Omega+B(0,r)$, we see that the support of $w(\cdot,t)$ is contained in $\{x\in\RR^\dim:0<\dist(x,\Omega)<r\}$. By construction
\[2A_{ij}(x,t)=L_{ij}(t)\left(\int\varphi\right)^{-1}\varphi(x),\]
so it follows from \cref{angular momentum vector potential} that $\int \xi_{ij}\cdot w=L_{ij}$. Finally, the claimed estimates follow at once from the bounds for $\varphi$ and the fact that $\int \varphi\geq \abs{\Omega}$.
\end{proof}

Now we have all the ingredients that we need to glue subsolutions in space. The following lemma is the key tool in this section. It will be used not only in the proof of~\cref{extender subs con soporte compacto}, but also in the convex integration scheme. We use skew-symmetric matrices instead of potential
vectors because the lemma is stated in any dimension $n\geq 2$.
\begin{lemma}
\label{pegado subsoluciones}
Let $T>0$ and let $\Omega_1\Subset \Omega_2\subset \RR^\dim$ be bounded domains with smooth boundary. Let $(v_i,p_i,\mathring{R}_i)\in C^\infty(\Omega_2\times[0,T])$ be subsolutions for $i=1,2$. Let $r>0$ be sufficiently small. There exists a subsolution $(v,p,\mathring{R})\in C^\infty(\Omega_2\times[0,T])$ such that
\begin{equation}
    \label{subsolucion en cada trozo}
    (v,p,\mathring{R})(x,t)=\begin{cases}(v_1,p_1,\mathring{R}_1)(x,t) \qquad x\in \overline{\Omega}_1, \\ (v_2,p_2,\mathring{R}_2)(x,t)\qquad \dist(x,\Omega_1)\geq r\end{cases}
\end{equation}
if and only if for each connected component $\Sigma$ of $\partial \Omega_1$, and all times $t\in[0,T]$, we have
\begin{align}
&\int_{\Sigma}v_1\cdot \normal=\int_{\Sigma}v_2\cdot \normal, \label{comp1}\\
&\int_{\Sigma}\left[(a\cdot x)\partial_t v_1+(a\cdot v_1)v_1+p_1 a-a^t\mathring{R}_1\right]\cdot \normal= \nonumber \\&\qquad \quad =\int_{\Sigma}\left[(a\cdot x)\partial_t v_2+(a\cdot v_2)v_2+p_2 a-a^t\mathring{R}_2\right]\cdot \normal  \qquad \forall a\in\RR^\dim. \label{comp2}
\end{align}
Suppose that, in addition, we have $v_1=\Div A_1$ and $v_2=\Div A_2$ for some potentials $A_i\in C^\infty(\Omega_2\times I,\mathcal{A}^\dim)$. Then, there exists $A\in C^\infty(\Omega_2\times I,\mathcal{A}^\dim)$ such that $v=\Div A$ and $A(x,t)=A_2(x,t)$ if $\dist(x,\Omega_1)\geq r$. 
\end{lemma}
\begin{remark}
	\label{condiciones triviales}
	The compatibility conditions (\ref{comp1}) and (\ref{comp2}) are automatically satisfied if $\partial \Omega_1$ is connected or $\Omega_2=\RR^\dim$. This will be explained in the proof of the lemma.
\end{remark}
\begin{remark}
	The subsolution $(v_1,p_1,\mathring{R}_1)$ need not be defined in all of $\Omega_2$ and the subsolution $(v_2,p_2,\mathring{R}_2)$ need not be defined in $\Omega_1$. We have assumed this to simplify slightly the statement of the lemma. 
\end{remark}
\begin{proof}
First of all, note that a subsolution $(v_0,p_0,\mathring{R}_0)$ in a bounded domain $G$ with smooth boundary satisfies
\begin{align}
    &0=\int_G \Div v_0=\int_{\partial G}v_0\cdot\normal, \label{integral subs 1} \\
    &0=\int_G a\cdot\left[\partial_tv_0+\Div\left(v_0\otimes v_0+p_0\Id-\mathring{R}_0\right)\right]\label{integral subs 2} \\&\hspace{40pt}=\int_{\partial G}\left[(a\cdot x)\partial_t v_0+(a\cdot v_0)v_0+p_0 a-a^t\mathring{R}_0\right]\cdot \normal \qquad \forall a\in\RR^\dim, \nonumber
\end{align}
where we have used the divergence theorem, identity (\ref{green id}) and the fact that $\Div[(a\cdot x)\partial_tv_0]=a\cdot \partial_t v_0$ because $\partial_t v_0$ is divergence-free.

From these equations we readily deduce that the compatibility conditions (\ref{comp1}) and (\ref{comp2}) are automatically satisfied if $\partial\Omega_1$ is connected, as both integrals vanish for each field. In the case $\Omega_2=\RR^\dim$, we apply \cref{integral subs 1,integral subs 2} to the domain  bounded by each connected component of $\partial\Omega_1$. We conclude that both integrals vanish for each field in each connected component of $\partial \Omega_1$.

Next, we check that the conditions are necessary; we study (\ref{comp1}) because the expressions are shorter, but the argument for (\ref{comp2}) is exactly the same. First, if $\partial \Omega_1$ is connected, it readily follows from (\ref{integral subs 1}) that (\ref{comp1}) must be satisfied. Hence, we focus on bounded domains $\Omega_1$ whose boundary is not connected. In that case, $\RR^\dim \backslash\Omega_1$ must have at least one bounded connected component. Given a bounded connected component of $\RR^\dim\backslash \Omega_1$, we define $G$ to be its intersection with $\Omega_2$. Then, $\partial G$ is composed of a connected component $\Sigma$ of $\partial \Omega_1$ and (possibly) some connected components $\Sigma'_1, \dots, \Sigma'_m$ of $\partial \Omega_2$. Since $v$ equals $v_1$ on $\Sigma$ and $v_2$ on the other connected components of $\partial G$, it follows from (\ref{integral subs 1}) that:
\[0=\int_{\partial G} v\cdot \normal=\int_{\Sigma}v_1\cdot\normal+\sum_{i=1}^\dim\int_{\Sigma'_i}v_2\cdot\normal.\]
On the other hand, applying (\ref{integral subs 1}) to $v_2$ on $G$, we have
\[-\int_\Sigma v_2\cdot\normal=\sum_{i=1}^\dim\int_{\Sigma'_i}v_2\cdot\normal,\]
which, combined with the previous equation, yields
\[\int_\Sigma v_1\cdot\normal=\int_\Sigma v_2\cdot\normal.\]
Since this applies to any bounded connected component of $\RR^\dim\backslash \Omega_1$, we can combine it with (\ref{integral subs 1}) with $G=\Omega_1$ to obtain an analogous identity for the remaining connected component of $\partial\Omega_1$, that is, the boundary of the unbounded connected component of $\RR^\dim\backslash \Omega_1$.  We conclude (\ref{comp1}).

Let us now prove that the compatibility conditions (\ref{comp1}) and (\ref{comp2}) are also sufficient. Let $r>0$ be small enough so that $\{x\in \Omega_2:\dist(x,\Omega_1)=r\}$ is diffeomorphic to $\partial \Omega_1$. We define $U\coloneqq\{x\in \Omega_2:0<\dist(x,\Omega_1)<r\}$. Then, the condition~(\ref{comp1}) ensures that there exists $A_{12}\in C^\infty(\overline{U}\times[0,T],\mathcal{A}^\dim)$ such that $v_2-v_1=\Div A_{12}$ in $U\times[0,T]$. Indeed, let $U_i$ be a connected component of $U$ and let $\Sigma_i$ and $\Sigma'_i$ be the connected components of $\partial U_i$, where $\Sigma_i\subset \partial\Omega_1$. Using~(\ref{comp1}) and the fact that $v_2-v_1$ is divergence-free:
\[\int_{\Sigma'_i}(v_2-v_1)\cdot\normal=\int_{\partial U}(v_2-v_1)\cdot\normal-\int_{\Sigma_i}(v_2-v_1)\cdot\normal=0.\]
Hence, the flux of $v_2-v_1$ through each connected component of $U$ vanishes, so by \cref{invertir matrices} there exists $A_{12}\in C^\infty(\overline{U}\times[0,T],\mathcal{A}^\dim)$ such that $v_2-v_1=\Div A_{12}$.

Next, using \cref{cutoff} we choose a cutoff function $\varphi\in C^\infty_c(\Omega_1+B(0,r))$ that equals 1 in a neighborhood of $\Omega_1$. We define
\begin{align*}
    v&\coloneqq \varphi \,v_1+(1-\varphi)v_2+w_c+w_L\equiv \varphi \,v_1+(1-\varphi)v_2+w,\\
    \widetilde{p}&\coloneqq \varphi\, p_1+(1-\varphi)p_2,
\end{align*}
where $w_c\coloneqq A_{12}\cdot \nabla \varphi$ so that $\varphi\, v_1+(1-\varphi)v_2+w_c$ is divergence-free. The additional correction $w_L$ is a divergence-free field supported within $U$ that will be defined later. Its purpose is to cancel the angular momentum so that the gluing can be performed in the interior of $U$. After a tedious computation we obtain
\begin{equation}
\label{ec pegado}
\partial_t v+\Div (v\otimes v)+\nabla \widetilde p=\Div\left(\varphi\,\mathring{R}_1+(1-\varphi)\mathring{R}_2+S_1\right)+\partial_t w_L+M\cdot\nabla\varphi,
\end{equation}
where
\begin{align}
	\label{def S1 pegado}
S_1&\coloneqq -\varphi(1-\varphi)(v_1-v_2)\otimes(v_1-v_2)+w\otimes\left(v-\frac{1}{2}w\right)+\left(v-\frac{1}{2}w\right)\otimes w,\\
 \label{def M pegado}
M&\coloneqq \partial_t A_{12}+v_1\otimes v_1-v_2\otimes v_2+(p_1-p_2)\Id-\mathring{R}_1+\mathring{R}_2.
\end{align}
Let $\widetilde{\rho}\coloneqq M\cdot \nabla\varphi$ and $\rho=\widetilde{\rho}+\partial_t w_L$. Our goal is to find $S_2\in C^\infty(\Omega_2\times[0,T],\mathcal{S}^\dim)$ supported on $U$ for all $t\in[0,T]$ and such that $\Div S_2=\rho$. Thus, we may set $R=\varphi\,\mathring{R}_1+(1-\varphi)\mathring{R}_2+S_1+S_2$ and absorb the trace into the pressure, obtaining $\mathring{R}$ and the final pressure $p$. To do so, first we must check that $\rho$ satisfies the compatibility conditions (\ref{integrability condition divergence}).

Note that $\widetilde{\rho}=\Div(\varphi M)$ because $\Div M=0$, since $(v_i,p_i,\mathring{R}_i)$ are subsolutions. Hence, by the divergence theorem for any $a\in\RR^\dim$ we have
\begin{equation}
\label{compatibilidad rho}
\int_U a\cdot \widetilde{\rho}=\int_{\partial U}a^t(\varphi M)\normal=\int_{\partial \Omega_1} a^tM\normal.
\end{equation}
Note that
\begin{align}
\int_{\partial\Omega_1}a^t(\partial_tA_{12})\normal&=-\int_{\partial\Omega_1}\normal^t(\partial_t A_{12})\nabla(a\cdot x)=\int_{\partial\Omega_1}(a\cdot x)\Div(\partial_t A_{12})\cdot \normal \nonumber \\&=\int_{\partial\Omega_1}(a\cdot x)(\partial_t v_1-\partial_t v_2)\cdot \normal. \label{integrar partes A}
\end{align}
Therefore, combining \cref{comp2,compatibilidad rho,integrar partes A} we conclude that $\int_U a\cdot \widetilde{\rho}=0$ for all $a\in\RR^\dim$. Since $\partial_t w_L$ is divergence-free and its support is contained in $U$, the same holds for $\partial_t w_L$, so $\int_U a\cdot\rho=0$ for all $a\in\RR^\dim$.

Next, we study the product with non-constant Killing fields. For each pair $1\leq i<j\leq \dim$ we define
\[l_{ij}(t)\coloneqq \int_U\xi_{ij}\cdot \widetilde{\rho}(x,t)\,dx,\]
where $\xi_{ij}$ are the elements of the basis of Killing fields defined in (\ref{definition Killing momento angular}). It will be useful later on to write the coefficients as:
\begin{equation}
    \label{expresion buena l_ij}
    l_{ij}=\int_U\xi_{ji}\cdot \Div(\varphi M)=\int_{\partial\Omega_1}\xi_{ij}^tM\normal,
\end{equation}
where we have used the fact that Killing fields are divergence-free as well as the values of $\varphi$ on $\partial \Omega_1$ and $\partial\Omega_2$. We define
\begin{equation}
	\label{def Lij}
	L_{ij}(t)\coloneqq -\int_0^tl_{ij}(s)\,ds.
\end{equation}
We then define the correction $w_L$ to be the divergence-free field obtained by applying \cref{controlar Killing} to the domain $\Omega_1$ with coefficients $L_{ij}\in C^\infty([0,T])$. Thus, we have
\[\int \xi_{ij}\cdot\partial_tw_L=\frac{d}{dt}\int\xi_{ij}\cdot w_L=L_{ij}'=-l_{ij}.\] 
We conclude that $\int_U\xi\cdot\rho=0$ for any Killing field $\xi$, as we wanted. Therefore, by \cref{invertir divergencia matrices} there exists $S_2\in C^\infty(\Omega_2\times[0,T],\mathcal{S}^\dim)$ supported on $U$ for all $t\in[0,T]$ and such that $\Div S_2=\rho$. We define the final pressure and the Reynolds stress as
\begin{align*}
	p&\coloneqq \T{p}-\frac{1}{\dim}\tr\left(S_1+S_2\right)= \varphi\, p_1+(1-\varphi)p_2-\frac{1}{\dim}\tr\left(S_1+S_2\right), \\
	\R{R}&\coloneqq\varphi\,\mathring{R}_1+(1-\varphi)\mathring{R}_2+S_1+S_2-\frac{1}{\dim}\tr\left(S_1+S_2\right)\Id.
\end{align*}
It follows from \cref{ec pegado} that the resulting triplet $(v,p,\mathring{R})$ is a subsolution and it satisfies (\ref{subsolucion en cada trozo}) because $S_1$ and $S_2$ are supported in $U$ for all $t\in[0,T]$.

Finally, let us consider that the velocity fields are given by $v_i=\Div A_i$. In that case, we may simply take $A_{12}=A_2-A_1$ instead of constructing a suitable potential using \cref{invertir matrices}. We see that
\[v-w_L=\varphi v_1+(1-\varphi)v_2+(A_2-A_1)\cdot\nabla \varphi=\Div(\varphi A_1+(1-\varphi)A_2).\]
Inspecting \cref{controlar Killing} leads us to define
\begin{equation}
	\label{A pegado}
	A\coloneqq \varphi A_1+(1-\varphi)A_2+\sum_{1\leq i<j\leq n}L_{ij}(t)\,(e_i\otimes e_j-e_j\otimes e_i)\left(2\int \varphi\right)^{-1}\varphi(x).
\end{equation}
Hence, we have $v=\Div A$ and we see that $A$ equals $A_2$ in $\{\dist(x,\Omega_1)\geq r\}$ because $\varphi$ vanishes in a neighborhood of this set.
\end{proof}

In the convex integration scheme we will need estimates of the glued subsolution. For the sake of clarity, we keep them separate in a different lemma:
\begin{lemma}
\label{estimaciones pegado subsolucones}
	Let $\alpha\in (0,1)$. In the conditions of \cref{pegado subsoluciones} and using the notation of its proof, the new subsolution satisfies:
	\begin{align}
		\norm{v-(\varphi v_1+(1-\varphi)v_2)}_N&\lesssim Tr^{-(N+1)}\norm{M}_{0;\hspace{0.3pt}U}+\sum_{k=0}^Nr^{-(k+1)}\norm{A_{12}}_{N-k;\hspace{0.3pt}U}, \\
		\norm{\partial_t(v-\varphi v_1-(1-\varphi)v_2)}_N&\lesssim r^{-(N+1)}\norm{M}_{0;\hspace{0.3pt}U}+\sum_{k=0}^Nr^{-(k+1)}\norm{\partial_t A_{12}}_{N-k;\hspace{0.3pt}U},\\
		\|\R{R}-\varphi\,\R{R}_1-(1-\varphi)\R{R}_2\|_0&\lesssim r^{-\alpha}\norm{M}_{0;\hspace{0.3pt}U}+\norm{v_1-v_2}_{0;\hspace{0.3pt}U}^2\\&\hspace{10pt}+(\norm{v_1}_{0;\hspace{0.3pt}U}+\norm{v_2}_{0;\hspace{0.3pt}U}+\norm{w}_{0;\hspace{0.3pt}U})\norm{w}_{0;\hspace{0.3pt}U}, \nonumber
	\end{align}
	In addition, if $v_1=\Div A_1$ and $v_2=\Div A_2$, the potential $A$ satisfies
	\begin{align}
		\norm{A-(\varphi A_1+(1-\varphi)A_2)}_N&\lesssim Tr^{-N}\norm{M}_{0;\hspace{0.3pt}U}, \\
		\norm{\partial_t(A-\varphi A_1-(1-\varphi)A_2)}_N&\lesssim r^{-N}\norm{M}_{0;\hspace{0.3pt}U}.
	\end{align}
	The implicit constants in these inequalities depend on $\Omega_1$, $N$ and $\alpha$.
\end{lemma}
\begin{proof}
	We begin by estimating $w_c=A_{12}\cdot\nabla\varphi$. Since $\varphi$ satisfies $\norm{\varphi}_N\lesssim r^{-N}$ and it is independent of time, it is clear that
	\[\norm{w_c}_N\lesssim\sum_{k=0}^Nr^{-(k+1)}\norm{A_{12}}_{N-k;\hspace{0.3pt}U}, \qquad \norm{\partial_t w_c}_N\lesssim\sum_{k=0}^Nr^{-(k+1)}\norm{\partial_t A_{12}}_{N-k;\hspace{0.3pt}U}\] 
	because the support of $\nabla \varphi$ is contained in $U$. Regarding $w_L$, it follows from (\ref{expresion buena l_ij}) that $\abs{l_{ij}}\lesssim \norm{M}_{0;\hspace{0.3pt}U}$, so $\abs{L_{ij}}\lesssim T\norm{M}_{0;\hspace{0.3pt}U}$. Hence, by \cref{controlar Killing} we have the bounds
	\[\norm{w_L}_N\lesssim Tr^{-(N+1)}\norm{M}_{0;\hspace{0.3pt}U}, \qquad \norm{\partial_t w_L}_N\lesssim r^{-(N+1)}\norm{M}_{0;\hspace{0.3pt}U}.\]
	The claimed estimates for $v-(\varphi v_1+(1-\varphi)v_2)=w_c+w_L$ follow at once. Let us now focus on the Reynolds stress. Using the assumption $\norm{w}_0\leq \norm{v_1}_0+\norm{v_2}_0$, we deduce from the definition (\ref{def S1 pegado}) that
	\[\norm{S_1}_0\lesssim \norm{v_1-v_2}_{0;\hspace{0.3pt}U}^2+(\norm{v_1}_{0;\hspace{0.3pt}U}+\norm{v_2}_0)\norm{w}_{0;\hspace{0.3pt}U}.\]
	Concerning $S_2$, let us first estimate $\rho$:
	\[\norm{\rho}_0\leq \norm{M\cdot\nabla\varphi}_{0;\hspace{0.3pt}U}+\norm{\partial_tw_L}_0\lesssim r^{-1}\norm{M}_{0;\hspace{0.3pt}U}.\]
	Since the support of $\rho(\cdot,t)$ is contained in $\{x\in\RR^\dim:0<\dist(x,\Omega_1)<r\}$, we may apply \cref{frecuencias bajas de dominio pequeño}, obtaining
	\[\norm{\rho}_{B^{-1+\alpha}_{\infty,\infty}}\lesssim r^{1-\alpha}\norm{\rho}_0\lesssim r^{-\alpha}\norm{M}_{0;\hspace{0.3pt}U}.\]
	Hence, it follows from the estimates in \cref{invertir divergencia matrices} that
	\[\norm{S_2}_0\lesssim \norm{\rho}_{B^{-1+\alpha}_{\infty,\infty}}\lesssim r^{-\alpha}\norm{M}_{0;\hspace{0.3pt}U}.\]
	Since
	\[\R{R}-(\varphi\,\R{R}_1+(1-\varphi)\R{R}_2)=S_1+S_2-\frac{1}{\dim}\tr(S_1+S_2)\Id,\]
	the claimed bound follows.
	
	Finally, let us estimate $A$ in the case that the velocities are given by $v_i=\Div A_2$. By (\ref{A pegado}) we have
	\[A-(\varphi A_1-(1-\varphi)A_2)=\sum_{1\leq i<j\leq n}L_{ij}(t)\,(e_i\otimes e_j-e_j\otimes e_i)\left(2\int \varphi\right)^{-1}\varphi(x).\]
	The claimed bounds then follow from the estimates derived in \cref{controlar Killing}.
\end{proof}

\cref{pegado subsoluciones} is almost what we want, but it can be made a bit sharper. In particular, in \cref{extender subs con soporte compacto} we do not want to assume that the subsolution $(v_0,p_0,\R{R}_0)$ is defined in a neighborhood of $\overline{\Omega}_0$. Fortunately, it turns out that all subsolutions can be extended, at least a little bit. Our operator $\mathscr{L}$ is essential for this:
\begin{lemma}
	\label{extender fuera}
	Let $\Omega_0\subset \mathbb{R}^\dim$ be a bounded domain with smooth boundary and let $I\subset \mathbb{R}$ be a closed and bounded interval. Let $(v_0,p_0,\R{R}_0)\in C^\infty(\overline{\Omega}_0\times I)$ be a subsolution in $\Omega_0\times I$. Let $\Omega$ be a sufficiently small open neighborhood of $\overline{\Omega}_0$. Then, there exists a subsolution $(v,p,\R{R})\in C^\infty(\Omega\times I)$ that extends $(v_0,p_0,\R{R}_0)$.
\end{lemma}

\begin{remark} The open neighborhood $\Omega$ need not be very small. It only needs to be bounded and such that each connected component of $\mathbb{R}^\dim\backslash \overline{\Omega}_0$ has nonempty intersection with $\mathbb{R}^\dim\backslash \overline{\Omega}$.
\end{remark}

\begin{proof}
	We begin by constructing the velocity field $v$. We choose $\rho\in C^\infty_c(\mathbb{R}^\dim\times I,\mathbb{R})$ such that $\supp \rho(\cdot,t)$ is contained in $\mathbb{R}^\dim\backslash \overline{\Omega}$ for all $t\in I$ and such that
	\[\int_{G}\rho(x,t)\,dx=\int_{\partial G} v_0\cdot \normal\qquad \forall t\in I\]
	for each bounded connected component $G$ of $\mathbb{R}^\dim\backslash \overline{\Omega}_0$, whose boundary we have oriented with the outer normal with respect to $G$. This can be done if $\Omega$ is a sufficiently small neighborhood of $\overline{\Omega}_0$ so that the intersection of $G$ with $\mathbb{R}^\dim\backslash \overline{\Omega}$ is nonempty. 
	
	Next, let $\T{v}\coloneqq \nabla \,\Delta^{-1}\rho$ so that $\widetilde{v}\in C^\infty(\mathbb{R}^\dim\times I,\mathbb{R}^\dim)$ and $\Div \widetilde{v}=\rho$ for all $t\in I$. By the divergence theorem we have
	\[\int_{\partial G} \widetilde{v}\cdot \nu=\int_G \rho=\int_{\partial G} v_0\cdot \normal\]
	in each bounded connected component $G$ of $\mathbb{R}^\dim\backslash \overline{\Omega}_0$ and for all $t\in I_0$. In addition, $\widetilde{v}-v_0$ is divergence-free in $\Omega_0\times I$ because $\rho$ vanishes in this set by construction. In particular, by the divergence theorem we have $\int_{\partial \Omega_0}(\widetilde{v}-v_0)\cdot \normal=0$, from which we conclude
	\[\int_\Sigma (\widetilde{v}-v_0)\cdot \nu=0 \qquad \forall t\in I_0\]
	for all connected components $\Sigma$ of $\partial \Omega_0$.  Therefore, by Lemma~\ref{invertir matrices} there exists $A\in C^\infty(\overline{\Omega}_0\times I,\mathcal{A}^n)$ such that $\Div A=\widetilde{v}-v_0$ in $\Omega_0\times I$. We choose a smooth extension $\widetilde{A}\in C^\infty(\mathbb{R}^\dim\times I,\mathbb{R}^{\dim\times \dim})$ and then we take the skew-symmetric part, so that $\widetilde{A}\in C^\infty(\mathbb{R}^\dim\times \mathbb{R},\mathcal{A}^\dim)$. We define:
	\[v\coloneqq \widetilde{v}-\Div\widetilde{A}\in C^\infty(\mathbb{R}^\dim\times I,\mathbb{R}^\dim).\]
	Since the support of $\rho(\cdot,t)$ is contained in $\mathbb{R}^\dim\backslash \overline{\Omega}$ and the second term is divergence-free, $v$ is divergence-free in $\Omega\times I$. In addition, our choice of $\widetilde{A}$ ensures that the restriction of $v$ to $\overline{\Omega}_0\times I$ is $v_0$, as we wanted.
	
	Next, we will extend $S_0\coloneqq v_0\otimes v_0+p_0\Id-\R{R}_0$ in a similar manner. First, we choose $f\in C^\infty_c(\mathbb{R}^\dim\times I,\mathbb{R}^n)$ such that $\supp f(\cdot,t)$ is contained in $\mathbb{R}^\dim\backslash \overline{\Omega}$ for all $t\in I$ and such that
	\begin{equation}
		\label{condition f}
		\int_{G_i} \xi\cdot f=\int_{G_i} \xi\cdot\partial_t v+\int_{\partial G_i} \xi^tU_0\normal\qquad \forall \xi\in \ker \gradsim,\;\forall t\in I
	\end{equation}
	for all bounded connected components $G_i$ of $\mathbb{R}^n\backslash \overline{\Omega}_0$, whose boundary we have oriented with the outer normal with respect to $G_i$. To find such an $f$, we fix a nonnegative radial function $\psi\in C^\infty_c(B(0,1))$ and a ball $\overline{B}(x_i,r_i)\subset G_i$. Due to symmetry, for any two elements $w_j\neq w_k$ of the basis $\mathcal{B}$ defined in (\ref{definition base Killing}) we have
	\[\int \psi\left(r_i^{-1}(x-x_i)\right)\,w_j(x)\cdot w_k(x)\,dx=0.\]
	Therefore, it suffices to choose
	\[f(x)\coloneqq \sum_i\sum_{j=1}^{\dim(\dim+1)/2}c_{ij}\,\psi\left(r_i^{-1}(x-x_i)\right)w_j(x)\]
	for the appropriate coefficients $c_{ij}$. Then, noticing that $\partial_t v$ is a smooth vector field with bounded derivatives, we define $\T{S}\coloneqq \mathscr{R}(-\partial_t v+f)$ so that $\T{S}\in C^\infty(\RR^\dim\times I,\mathcal{S}^n)$, and 	\begin{equation}
		\label{ec S extender fuera}
		\Div \T{S}=-\partial_t v+f.
	\end{equation}
	Since $v$ extends $v_0$ and $(v_0,p_0,\R{R}_0)$ is a subsolution in $\Omega_0\times I$, we have
	\begin{equation}
		\label{difference S S0}
		\Div(\T{S}-S_0)(x,t)=0\qquad \forall(x,t)\in \Omega_0\times I
	\end{equation}
	because $f$ vanishes in that set. In addition, using (\ref{green id}) it follows from (\ref{ec S extender fuera}) and (\ref{condition f}) that
	\[\int_{\partial G} \xi^t (\T{S}-S_0)\normal=0 \qquad \forall \xi\in\ker \gradsim\; \forall t\in I\]
	for all bounded connected components $G$ of $\mathbb{R}^m\backslash\overline{\Omega}_0$. Due to (\ref{green id}) and (\ref{difference S S0}), this integral also vanishes for the remaining connected component of $\partial\Omega_0$, that is, the connected component that separates $\Omega_0$ and the unbounded connected component of $\RR^\dim\backslash \Omega_0$.
	
	We conclude that $\T{S}-S_0$ is in $\mathcal{G}(\overline{\Omega}_0\times I)$. Therefore, by Lemma~\ref{key lemma} there exists $\T{E}\in \mathcal{P}(\overline{\Omega}_0\times I)$ such that $\mathscr{L}(\T{E})=\T{S}-S_0$ in $\Omega_0\times I$. We choose a smooth extension $E\in C^\infty(\mathbb{R}^\dim\times I,\mathbb{R}^{\dim^4})$ and then we make the appropriate antisymmetrization. We define
	\[S\coloneqq \T{S}-\mathscr{L}(E)\in C^\infty(\RR^\dim\times I,\mathcal{S}^\dim)\]
	By construction of $E$, we see that $S$ extends $S_0$. Additionally, we have
	\[\Div S=\Div \T{S}=-\partial_t v+f\]
	because the image of $\mathscr{L}$ is contained in the kernel of the divergence. We define
	\begin{align*}
		\R{R}&\coloneqq v\otimes v-S-\frac{1}{\dim}\tr(v\otimes v-S)\Id, \\
		p&\coloneqq -\frac{1}{\dim}\tr(v\otimes v-S).
	\end{align*}
	Since $f(\cdot,t)$ vanishes in $\Omega$ for all $t\in I$, we conclude that $(v,p,\R{R})$ is a subsolution in $\Omega\times I$ that extends $(v_0,p_0,\R{R}_0)$.
\end{proof}

Combining \cref{pegado subsoluciones} and \cref{extender fuera} we can finally prove \cref{extender subs con soporte compacto}:
\begin{proof}[Proof of~\cref{extender subs con soporte compacto}]
Working with each connected component of $\Omega_0$, we may assume that both $\Omega_0$ and $\Omega$ are connected (i.e., domains).	Then, reducing $\Omega$ if necessary, by \cref{extender fuera} we may assume that $(v_0,p_0,\R{R}_0)$ is a subsolution in $\Omega\times I$. The result then follows by applying \cref{pegado subsoluciones} with the domains $\Omega_0, \Omega$ and subsolutions $(v_0,p_0,\R{R}_0)$ and $(0,0,0)$.
\end{proof}

\section{Proof of Theorem~\ref{teorema integración convexa}}\label{section overview}

The construction of a weak solution to the Euler equations stated in Theorem~\ref{teorema integración convexa} consists in an iterative argument which is presented in Subsection~\ref{ss.iter}, cf. Proposition~\ref{prop iteracion integracion convexa}. This proposition together with Lemma~\ref{lema integración convexa una iteración} allow us to prove the theorem in Subsection~\ref{ss.proofteo}. We want to remark that most of this article, i.e., Sections~\ref{proof of prop} to~\ref{section perturbation}, is devoted to prove Proposition~\ref{prop iteracion integracion convexa}, which is the key result for our convex integration scheme.

\cref{ss.iter} and \cref{ss.proofteo} follow the general outline of Sections~2.1 and~2.2 in \cite{Onsager_final} but with two important differences: here the initial subsolution will be nontrivial, i.e., different from $(0,0,0)$, and the perturbations will be supported in a subset of $\Omega$ instead of the whole space. Regarding the first issue, we use \cref{lema integración convexa una iteración} to help start the iterative process. Concerning the second point, we introduce suitable sets related to the distance to $\partial\Omega$ and the size of $\mathring{R}_0$. The perturbations will be localized to these sets, which is summarized in an additional inductive hypothesis, (\ref{inductive 1}).

\subsection{The iterative process}\label{ss.iter}
Let us assume all along this subsection that the subsolution $(v_0,p_0)(\cdot,t)$ is compactly supported for each time $t\in[0,T]$. We will construct the desired weak solution of the Euler equations as the limit of a sequence of subsolutions, that is, at a given step $q\geq0$ we have $(v_q,p_q,\mathring{R}_q)\in C^\infty(\RR^3\times[0,T])$ solving the Euler-Reynolds system:
\begin{equation}
\begin{cases}
\partial_t v_q+\Div(v_q\otimes v_q)+\nabla p_q=\Div\mathring{R}_q,\\
\Div v_q=0,
\end{cases}
\end{equation}
to which we add the constraint that 
\begin{equation}
\tr\mathring{R}_q=0.
\end{equation}
The matrix $\mathring{R}_q$ measures the deviation from being a solution of the Euler equations. The goal of the process is to make $\mathring{R}_q$ vanish at the limit $q\to +\infty$, so that the limit field is a weak solution of the Euler equations.

Assume we are given the initial subsolution $(v_0,p_0,\mathring{R}_0)\in C^\infty(\RR^3\times[0,T])$. Let us then show how to construct the rest of the terms iteratively. To construct the subsolution at step $q$ from the one in step $q-1$, we will add an oscillatory perturbation with frequency $\lambda_q$. Meanwhile, the size of the Reynolds stress will be measured by an amplitude $\delta_q$. These parameters are given by
\begin{align}
	\lambda_q&=2\pi \lceil a^{b^q}\rceil, \label{def lambdaq}\\
	\delta_q&=\lambda_q^{-2\beta}, \label{def deltaq}
\end{align}
where $\lceil x\rceil$ denotes the ceiling, that is, the smallest integer $n\geq x$. The parameters $a,b>1$ are very large and very close to 1, respectively. They will be chosen depending on the exponent $0<\beta<1/3$ that appears in \cref{teorema integración convexa}, on $\Omega$ and on the initial subsolution. We introduce another parameter $\alpha>0$ that will be very small. The necessary size of all of the parameters will be discovered in the proof.

Throughout the process we will also try to achieve a given energy profile $e\in C^\infty([0,T])$, which must satisfy the inequality~(\ref{condición energía 1}). We will also assume
\begin{equation}
	\label{condición energía 2}
	\sup_{t\in[0,T]}\abs{\frac{d}{dt}e(t)}\leq 1.
\end{equation}
We will see that this can be assumed without losing generality.

Unlike the construction on the torus in~\cite{Onsager_final}, it is essential that we only perturb the field in the region where the Reynolds stress is nonzero. Hence, we have to pay special attention to the support of the fields. 

Since the map $(v_0,p_0,\mathring{R}_0)(\cdot,t)$ is assumed to be compactly supported at each time $t\in[0,T]$, we shall see that with a suitable rescaling we may assume that its support and $\Omega$ are contained in $(0,1)^3$. This is useful because sometimes it will be convenient to consider that we are working with periodic boundary conditions (that is, in $\mathbb{T}^3$) to reuse the results in~\cite{Onsager_final}. On the other hand, $\mathring{R}_0(\cdot,t)$ is supported in a potentially smaller domain $\overline{\Omega}$. In our construction we must ensure that we do not perturb the subsolution outside of this set.

It will be convenient to do an additional rescaling in our problem. In the rescaled problem the initial subsolution will depend on $a$, but we assume that nevertheless there exists a sequence $\{y_N\}_{N=0}^\infty$ independent of the parameters such that
\begin{align}
	\norm{v_0}_N+\norm{\partial_t v_0}_N&\leq y_N, \label{est v no cambia con el rescalado}\\
	\norm{p_0}_N&\leq y_N, \\
	\|\R{R}_0\|_N+\|\partial_t \R{R}_0\|_N&\leq y_N. \label{est R no cambia con el rescalado}
\end{align}

Since the initial Reynolds stress $\mathring{R}_0$ and its derivatives vanish at $\partial \Omega\times[0,T]$, for any $k\in\NN$ there exists a constant $C_k$ such that for any $x\in \Omega$ we have
\begin{equation}
	\vert\mathring{R}_0(x,t)\vert\leq C_k\dist(x,\partial \Omega)^k.
    \label{R inicial}
\end{equation}
The constants $C_k$ are independent of $a$ by (\ref{est R no cambia con el rescalado}). We define
\begin{equation}
    \label{def d_q}
    d_q\coloneqq\left(\frac{\delta_{q+2}\lambda_{q+1}^{-6\alpha}}{4C_{10}}\right)^{1/10}.
\end{equation}
Hence, we have
\begin{equation}
	\label{tamaño R0 cerca frontera}
    \vert\mathring{R}_0(x,t)\vert\leq \frac{1}{4}\delta_{q+2}\lambda_{q+1}^{-6\alpha} \qquad \forall x\in \Omega,\; \dist(x,\partial\Omega)\leq d_q.
\end{equation} 
At step $q$ the perturbation will be localized in a central region
\begin{equation}
    A_q\coloneqq \{x\in\Omega:\dist(x,\partial\Omega)\geq d_q\}
    \label{def A_q primera}
\end{equation}
so that $(v_q,p_q,\mathring{R}_q)$ equals the initial subsolution in $(\RR^3\backslash A_q)\times[0,T]$. Note that $d_q\to 0$ as $q\to\infty$ because so does $\delta_{q+2}$. Therefore, in the limit the perturbation covers all of the region where $\mathring{R}_0$ is nonzero. However, the velocity is not modified outside this set.

As we have mentioned, the error introduced in the gluing step of~\cite{Onsager_final} is spread throughout the whole space. To avoid this, we introduce an additional gluing \textit{in space}, which we will explain in more detail in the following sections.

The complete list of inductive estimates is the following:
\begin{align}
	(v_q,p_q,\mathring{R}_q)&=(v_0,p_0,\mathring{R}_0) \qquad \text{outside } A_q\times[0,T], \label{inductive 1}\\
	\|\mathring{R}_q\|_0&\leq \delta_{q+1}\lambda_q^{-6\alpha}, \label{inductive 2}\\ 
	\norm{v_q}_1&\leq M\delta_q^{1/2}\lambda_q, \label{inductive 3}\\
	\norm{v_q}_0&\leq 1-\delta_q^{1/2}, \label{inductive v0}\\
	\delta_{q+1}\lambda_q^{-\alpha}&\leq e(t)-\int_\Omega\abs{v_q}^2dx\leq \delta_{q+1}, \label{inductive 5}
\end{align}
where $M$ is a geometric constant that depends on $\Omega$ and is fixed throughout the iterative process.

\begin{remark}
\label{varias componentes}
If $\Omega$ has several connected components $\Omega^j$, we may fix an energy profile $e^j$ in each of them. In that case, (\ref{inductive 5}) would have to be replaced by
\[\delta_{q+1}\lambda_q^{-\alpha}\leq e^j(t)-\int_{\Omega^j}\abs{v_q}^2dx\leq \delta_{q+1}.\]
Since the construction does not differ much, for simplicity we will assume that $\Omega$ is connected.
\end{remark}

The following proposition is the key result to prove Theorem~\ref{teorema integración convexa}, because it establishes the existence of the iterative scheme in the convex integration process.

\begin{proposition}
	\label{prop iteracion integracion convexa}
	Let $T>0$ and let $\Omega\subset (0,1)^3\subset \RR^3$ be an open set with smooth
boundary and with a finite number of connected components. Let $(v_0,p_0,\mathring{R}_0)\in C^\infty(\RR^3\times[0,T])$ be a subsolution whose support is contained in $(0,1)^3\times[0,T]$ and such that $\supp\R{R}_0\subset \overline{\Omega}\times[0,T]$. Furthermore, assume that (\ref{est v no cambia con el rescalado})$-$(\ref{est R no cambia con el rescalado}) are satisfied for some sequence of positive numbers $\{y_N\}_{N=0}^\infty$. There exists a constant $M$ depending only on $\Omega$ with the following property: 
	\newline Assume $0<\beta<1/3$ and
	\begin{equation}
		\label{condicion b}
		1<b<\min\left\{\frac{1-\beta}{2\beta},\;\frac{11}{10}\right\}.
	\end{equation}
	Then there exists an $\alpha_0$ depending on $\beta$ and $b$ such that for any $0<\alpha<\alpha_0$ there is an $a_0$ depending on $\beta$, $b$, $\alpha$, $\Omega$ and $\{y_N\}_{N=0}^\infty$ such that for any $a\geq a_0$ the following holds: 
	\newline Given a strictly positive energy profile satisfying  (\ref{condición energía 2}) and a subsolution $(v_q,p_q,\mathring{R}_q)$ satisfying (\ref{inductive 1})$-$(\ref{inductive 5}), there exists a subsolution $(v_{q+1},p_{q+1},\mathring{R}_{q+1})$ satisfying the same equations (\ref{inductive 1})$-$(\ref{inductive 5}) with $q$ replaced by $q+1$. Furthermore, we have the estimate
	\begin{equation}
		\label{estimación cambio iteración}
		\norm{v_{q+1}-v_q}_0+\frac{1}{\lambda_{q+1}}\norm{v_{q+1}-v_q}_1\leq M\delta_{q+1}^{1/2}.
	\end{equation}
\end{proposition}

We wish to iterate this result to construct a sequence of subsolutions whose limit will be the desired weak solution. However, in order to start the process, the first term in the sequence must satisfy the inductive hypotheses (\ref{inductive 1})$-$(\ref{inductive 5}). Since we do not assume any bounds on $(v_0,p_0,\R{R}_0)$, these hypotheses will not be satisfied by the initial subsolution, in general. Although a time dilation would almost solve the problem, we need the following lemma to fully prepare the initial subsolution:
\begin{lemma}
	\label{lema integración convexa una iteración}
	Let $T>0$ and let $\Omega\subset (0,1)^3\subset \RR^3$ be an open set with smooth
boundary and with a finite number of connected components. Let $(v_0,p_0,\mathring{R}_0)\in C^\infty(\RR^3\times[0,T])$ be a subsolution whose support is contained in $(0,1)^3\times[0,T]$ and such that $\supp\R{R}_0\subset \overline{\Omega}\times[0,T]$. Let $\lambda>0$ be a sufficiently large constant. There exists a subsolution $(v,p,\mathring{R})\in C^\infty(\RR^3\times[0,T])$ such that for any $N\geq 0$ we have
	\[\norm{v}_N\lesssim \lambda^N, \hspace{50pt} \|\R{R}\|_0\leq\lambda^{-1/2} \]
	where the implicit constants are independent of $\lambda$. In addition, the energy satisfies
	\begin{equation}
		\int_\Omega \abs{v_0}^2dx<\int_\Omega \abs{v}^2dx<\int_\Omega \abs{v_0}^2dx+6\|\R{R}_0\|_0\abs{\Omega}. \label{energía una iteración}
	\end{equation}
	Furthermore, $(v,p,\R{R})$ equals the initial subsolution outside the set
	\[A_*\coloneqq \left\{x\in \Omega: \dist(x,\partial \Omega)>\lambda^{-1/12}\right\}\times[0,T].\]
\end{lemma}
While necessary, this result is nothing new and one could easily obtain it by combining~\cite{Continuous} and \cref{invertir divergencia matrices}, or using the ideas in~\cite{Isett}. For completeness, we sketch its proof at the end of Section~\ref{section perturbation}, considering a simplified version of the preceding construction.

\subsection{Proof of \cref{teorema integración convexa}}\label{ss.proofteo}
We first prove the theorem under the assumption that the subsolution $(v_0,p_0,\mathring{R}_0)(\cdot,t)$ is compactly supported for each time $t\in[0,T]$. This assumption will be relaxed in next subsection to a condition on the support of $\mathring{R}_0$. We fix $0<\beta<1/3$, we choose $b$ satisfying (\ref{condicion b}) and $\alpha$ smaller than the threshold given by \cref{prop iteracion integracion convexa}. Next, we use the scale invariance of the Euler equations and subsolutions
\[v_0(x,t)\mapsto v_0(\rho x,\rho t), \quad p_0(x,t)\mapsto p_0(\rho x,\rho t), \quad \R{R}_0(x,t)\mapsto \R{R}_0(\rho x,\rho t)\]
to assume that $\Omega\subset (0,1)^3$ and $\supp (v_0,p_0,\mathring{R}_0)\subset (0,1)^3\times[0,T]$. The desired energy profile must also be modified: $e(t) \mapsto \rho^{-3}e(t)$. Note that this preserves (\ref{condición energía 1}). For convenience, we denote the rescaled subsolution like the original. It then suffices to construct the desired solution in this case and then reverse the change of variables. 

Next, we use \cref{lema integración convexa una iteración} with $\lambda=\lambda_1^{12\alpha}$ to obtain a subsolution $(v_1,p_1,\R{R}_1)$ that equals the initial subsolution outside the set \[\{x\in \Omega:d(x,\partial \Omega)>\lambda_1^{-\alpha}\}\times[0,T]\] and satisfies the estimates
\[\norm{v_1}_N\leq c_N\lambda_1^{12N\alpha}, \hspace{50pt} \|\R{R}_1\|_0\leq \lambda_1^{-6\alpha},\]
where the constants $c_N$ are independent of $\lambda_1$ but they will depend on $\Omega$ and the initial subsolution. In addition, by (\ref{condición energía 1}) we have
\[e(t)>\int_\Omega \abs{v_0}^2dx+6\|\R{R}_0\|_0\abs{\Omega}>\int_\Omega\abs{v_1}^2dx.\]

Next, we use another scale invariance of the Euler equations (and the definition of subsolution):
\[v(x,t)\mapsto \Gamma \hspace{0.5pt}  v(x,\Gamma t), \quad p(x,t)\mapsto \Gamma^2 p( x,\Gamma t), \quad \R{R}(x,t)\mapsto \Gamma^2\R{R}( x,\Gamma t).\]
We choose
\[\Gamma\coloneqq \delta_2^{1/2}\max\left\{1,\supp_t\left(e(t)-\int_\Omega \abs{v_1(x,t)}^2dx\right)\right\}^{-1/2}\]
and we begin to work in this rescaled setting, which we will
indicate with a superscript $r$. We are thus working in the interval $[0,\T{T}]$, where $\T{T}\coloneqq \Gamma^{-1}T$, and we try to prescribe the energy profile $\tilde{e}\coloneqq \Gamma^2 e(t)$. By construction of $\Gamma$ we have 
\[\sup_t\left(\tilde{e}(t)-\int_\Omega \abs{v_1^r(x,t)}^2dx\right)\leq \delta_2\]
and
\[\inf_t\left(\tilde{e}(t)-\int_\Omega \abs{v_1^r(x,t)}^2dx\right)= \Gamma^2\inf_t\left(e(t)-\int_\Omega \abs{v_1(x,t)}^2dx\right).\]
It follows from the definition of $\Gamma$ that if $a$ is sufficiently large we have
\[\lambda_1^{\alpha}\left(\frac{\Gamma^2}{\delta_2}\right)\inf_t\left(e(t)-\int_\Omega \abs{v_0(x,t)}^2dx\right)\geq 1.\]
so (\ref{inductive 5}) holds. On the other hand,
\[\sup_t\abs{\tilde{e}'(t)}\leq \Gamma^{3/2}\sup_t\abs{e'(t)}\leq 1\]
because $\Gamma$ becomes arbitrarily small by increasing $a$.

Next, we observe that $(v_0^r,p_0^r,\R{R}{}_0^r)$ still satisfies (\ref{est v no cambia con el rescalado})$-$(\ref{est R no cambia con el rescalado}) with the same sequence $\{y_N\}_{N=0}^\infty$. Regarding $(v_1^r,p_1^r,\R{R}{}_1^r)$, it follows from the definition of the rescaling that 
\[\|\R{R}{}_1^r\|_0\leq \delta_2\lambda_1^{-6\alpha}.\]
On the other hand, since the constants $c_N$ are independent of $\lambda_1$, for sufficiently large $a$ we have
\begin{align*}
    \norm{v_1^r}_0&\leq\delta_{2}^{1/2}\norm{v_1}_0\leq \delta_{2}^{1/2}c_0\leq 1-\delta_1^{1/2},\\
    \norm{v_1^r}_1&\leq\delta_{2}^{1/2}\norm{v_1}_1\leq \delta_{2}^{1/2}c_1\lambda_1^{12\alpha}\leq M\delta_1^{1/2}\lambda_1.
\end{align*}
Finally, $(v_1^r,p_1^r,\R{R}{}_1^r)=(v_0^r,p_0^r,\R{R}{}_0^r)$ outside 
\[\{x\in \Omega:d(x,\partial \Omega)>\lambda_1^{-\alpha}\}\times[0,\T{T}].\]
Let us consider the sets $A_q$ defined in (\ref{def A_q primera}). We see that $(v_1^r,p_1^r,\R{R}{}_1^r)=(v_0^r,p_0^r,\R{R}{}_0^r)$ outside $A_1\times[0,\T{T}]$ for sufficiently small $\alpha$ and sufficiently large $a$. Remember that $\delta_q$ and $\lambda_q$ depend on $a$ through the expressions (\ref{def deltaq}) and (\ref{def lambdaq}), respectively.

From now on we assume that we are working in this rescaled problem and we omit the superscript $r$. Once we obtain the desired weak solution in this setting, to obtain the solution to the original problem it suffices to undo the scaling. By the previous discussion, the energy profile satisfies the inductive hypotheses (\ref{condición energía 2}) and the subsolution $(v_1,p_1,\R{R}_1)$ satisfies (\ref{inductive 1})$-$(\ref{inductive 5}). In addition, the initial subsolution $(v_0,p_0,\R{R}_0)$ satisfies (\ref{est v no cambia con el rescalado})$-$(\ref{est R no cambia con el rescalado}) for some sequence $\{y_N\}_{N=0}^\infty$ that does not depend on $a$. Therefore, we may apply \cref{prop iteracion integracion convexa} iteratively, obtaining a sequence of compactly supported smooth subsolutions $\{(v_q,p_q,\R{R}_q)\}_{q=1}^\infty$. 

It follows from (\ref{estimación cambio iteración}) that $v_q$ converges uniformly to some continuous map $v$. On the other hand, note that the pressure $p_q$ is the only compactly supported solution of 
\[\Delta p_q=\Div\Div(-v_q\otimes v_q+\R{R}_q).\]
Therefore, $p_q$ also converges to some pressure $p\in L^r(\RR^3)$ for any $1\leq r<\infty$. Since $\R{R}_q$ converges uniformly to $0$, we conclude that the pair $(v,p)$ is a weak solution of the Euler equations.

Furthermore, using (\ref{estimación cambio iteración}) we obtain
\begin{align*}
	\sum_{q=1}^\infty \norm{v_{q+1}-v_q}_{\beta'}&\leq \sum_{q=1}^\infty C(\beta',\beta) \norm{v_{q+1}-v_q}_{0}^{1-\beta'}\norm{v_{q+1}-v_q}_{\beta}^{\beta'}\\
	&\leq C(\beta',\beta) \sum_{q=1}^\infty (M\delta_{q+1}^{1/2})^{1-\beta'}\left(M\delta_{q+1}^{1/2}\lambda_q \right)^{\beta'}\\
	&\leq M \,C(\beta',\beta) \sum_{q=1}^\infty \lambda_q^{\beta'-\beta},
\end{align*}
so $\{v_q\}_{q=1}^\infty$ is uniformly bounded in $C^0_t C^{\beta'}_x$ for all $\beta'<\beta$. To recover the time regularity, see \cite{Onsager_final}.

Next, note that $(v_q,p_q,\R{R}_q)=(v_0,p_0,0)$ in $(\RR^3\backslash \Omega)\times[0,T]$ for all $q$ by (\ref{inductive 1}) and the definition of $A_q$. Hence, we have $(v,p)=(v_0,p_0)$ in $(\RR^3\backslash \Omega)\times[0,T]$. 

Finally, it follows from (\ref{inductive 5}) and the fact that $\delta_{q+1}\to 0$ as $q\to \infty$ that $\int_\Omega \abs{v(x,t)}^2dx=e(t)$, as we wanted. This completes the proof of the theorem for the case that the initial subsolution is compactly supported for all time.

\subsection{Dropping the compact support condition}\label{SS.drop}

Once we have proved \cref{teorema integración convexa} for the case that $(v_0,p_0,\R{R}_0)(\cdot,t)$ is compactly supported, it is easy to relax this condition and show that it suffices that $\R{R}_0(\cdot,t)$ is compactly supported.

Indeed, let us choose a bounded domain with smooth boundary $U\Supset \Omega$. As we mentioned in \cref{pegado subsoluciones} and Remark~\ref{condiciones triviales}, any subsolution defined in all $\RR^3$ automatically satisfies the conditions in \cref{extender subs con soporte compacto}. Hence, there exists a subsolution $(\T{v}_0,\T{p}_0,\R{\T{R}}_0)\in C^\infty(\RR^3\times[0,T])$ that extends $(v_0,p_0,\R{R}_0)$ outside $\overline{U}$ and that is compactly supported in each time slice. Since it is an extension, we see that the Reynolds stress vanishes in $\overline{U}\backslash \Omega$. 
	
We can now apply \cref{teorema integración convexa} in the case that the subsolution is compactly supported for each time slice, thus obtaining a weak solution $(v,p)$ with the appropriate regularity and such that $(v,p)=(v_0,p_0)$ in $(\overline{U}\backslash \Omega)\times[0,T]$. Since $\Omega\Subset U$, we may glue this region back into $(v_0,p_0)$, obtaining a weak solution that equals $(v_0,p_0)$ outside $\Omega\times[0,T]$. 
	
	Finally, since $\Omega$ and $\RR^3\backslash \overline{U}$ are disjoint, when we apply \cref{teorema integración convexa} we may fix the energy profile in each region independently. Therefore, we can prescribe the energy profile in $\Omega$ of the final weak solution to be any smooth function satisfying condition (\ref{condición energía 1}).

\section{Proof of Proposition~\ref{prop iteracion integracion convexa}}\label{proof of prop}

The different steps in the proof of Proposition~\ref{prop iteracion integracion convexa} are presented in Subsection~\ref{SS.stages}. These steps are elaborated in the next subsections, and we complete the proof of the proposition in Subsection~\ref{SS.prupro}. Roughly speaking, our proof adapts the arguments of~\cite{Onsager_final} to the nonperiodic setting. To do so, we introduce an additional step in the iteration: a gluing in space that ensures that the error does not spread out to the whole space when we glue in time. 

\subsection{Stages of the proof}\label{SS.stages}

\begin{enumerate}
	\setlength\itemsep{4pt}
	\item Preparing the subsolution. We mollify our subsolution $(v_q,p_q,\mathring{R}_q)$ to avoid the \textit{loss of derivatives problem}, obtaining a new subsolution $(v_\ell,p_\ell,\mathring{R}_\ell)$. It is convenient to glue it \textit{in space} to the original subsolution far from the turbulent zone.
	\item Gluing in space. We pick a collection of times $\{t_i\}\subset [0,T]$ and we consider the solutions $(v_i,p_i)$ of the Euler equations with initial data $\widetilde{v}_\ell(t_i)$. We glue them in space to $(\widetilde{v}_\ell,\widetilde{p}_\ell,\mathring{\widetilde{R}}_\ell)$, obtaining new subsolutions $(\widetilde{v}_i,\widetilde{p}_i,\mathring{\widetilde{R}}_i)$. The error $\mathring{\widetilde{R}}_i$ is small and these subsolutions equal $(v_0,p_0,\mathring{R}_0)$ near $\partial\Omega$.
	\item Gluing in time. We glue together the subsolutions $(\widetilde{v}_i,\widetilde{p}_i,\mathring{\widetilde{R}}_i)$, obtaining a new subsolution $(\overline{v}_q,\overline{p}_q,\mathring{\overline{R}}_q)$ in which most of the error is concentrated in temporally disjoint regions. The error remains localized within $\Omega$ owing to the fact that the differences between the subsolutions $(\widetilde{v}_i,\widetilde{p}_i,\mathring{\widetilde{R}}_i)$ vanish near $\partial\Omega$.
	\item Perturbation. We add a highly oscillatory perturbation to reduce the error. In fact, we add many corrections, each of them reducing the error in one of the temporally disjoint regions. These perturbations do not interact with each other, which yields the optimal estimates for the new subsolution $(v_{q+1},p_{q+1},\R{R}_{q+1})$. 
\end{enumerate}

Throughout the iterative process we will use the notation $x\lesssim y$ to denote $x\leq C y$ for a sufficiently large constant $C>0$ that is independent of $a$, $b$ and $q$. However, the constant is allowed to depend on $\alpha$, $\beta$, $\Omega$ and $\{y_N\}_{N=0}^\infty$ and it may change from line to line. 

\subsection{Preparing the subsolution}
The first step consists in mollifying the field in order to avoid the \textit{loss of derivatives} problem, which is typical of convex integration. The problem is the following: to control a Hölder norm of $(v_{q+1},p_{q+1},\mathring{R}_{q+1})$ we need estimates of higher-order Hölder norms of $(v_q,p_q,\mathring{R}_q)$. As the iterative process goes on, we need to estimate higher and higher Hölder norms of the initial terms of the sequence to control just the first few Hölder norms of the subsolution. However, if we mollify the subsolution $(v_q,p_q,\mathring{R}_q)$, we can control all the derivatives in terms of the first few Hölder norms and the mollification parameter.

Note that this process changes the subsolution in the whole space, yet mollification is only strictly necessary in $A_q\times[0,T]$, as $(v_q,p_q,\mathring{R}_q)$ equals $(v_0,p_0,\mathring{R}_0)$ outside of this set. Furthermore, it will be convenient for later estimates that the resulting subsolution equals the initial subsolution far from the turbulent zone, as this is a property that we want to impose unto $(v_{q+1},p_{q+1},\R{R}_{q+1})$. Hence, our approach consists in gluing the mollified subsolution to the initial subsolution, which is not quite demanding because both subsolutions are very close far from the turbulent zone. 

We begin by fixing a mollification kernel \textit{in space}  $\psi\in C^\infty_c(\RR^3)$ and we introduce the mollification parameter 
\begin{equation}
    \label{def ell}
    \ell\coloneqq \frac{\delta^{1/2}_{q+1}}{\delta_q^{1/2}\lambda_q^{1+3\alpha}}.
\end{equation}
Since $(\delta_{q+1}/\delta_q)^{1/2}=\lambda_q^{-\beta(b-1)}$ and our assumption (\ref{condicion b}) implies that $\beta(b-1)<1/2$, we see that we may choose $\alpha$ sufficiently small and $a$ sufficiently large so as to have 
\begin{equation}
    \label{tamaño ell}
    \lambda_q^{-3/2}\leq \ell\leq \lambda_q^{-1}.
\end{equation}
Note that our definition of $\ell$ differs from the one in \cite{Onsager_final} by a factor of $\lambda_q^{-3\alpha/2}$. The coefficients of $\alpha$ used in \cite{Onsager_final} are fined-tuned to their proof. Since we do some things different, it is not surprising that we have to change some of these coefficients. In general the factor $\lambda_q^{-3\alpha/2}$ is quite harmless, as only simple relationships like (\ref{tamaño ell}) are used throughout most of the paper, and these are the same here and in \cite{Onsager_final}. The actual value of $\ell$ is only used at the very end, when fine relationships between the parameters are needed to estimate $\R{R}_{q+1}$. We will study these situations when they arise, but in any case we will see that the extra factor is essentially irrelevant. Indeed, $\alpha$ is assumed to be so small that in those inequalities the term containing $\alpha$ is negligible. Our definition of $\ell$ leads to a different coefficient multiplying $\alpha$, but this only changes how small $\alpha$ has to be chosen, so it is not important.

After this brief digression, we define
\begin{align*}
    v_\ell&\coloneqq v_q\ast \psi_\ell,\\
    p_\ell&\coloneqq p_q\ast \psi_\ell+\abs{v_q}^2\ast \psi_\ell-\abs{v_\ell}^2,\\
    \mathring{R}_\ell&\coloneqq \mathring{R}_q\ast \psi_\ell-(v_q\mathring{\otimes}v_q)\ast \psi_\ell+v_\ell\mathring{\otimes}v_\ell,
\end{align*}
where the convolution with $\psi_\ell$ is in space only and $f\mathring{\otimes}g$ denotes the traceless part of the tensor $f\otimes g$. It is easy to check that the triplet $(v_\ell,p_\ell,\mathring{R}_\ell)$ is a subsolution and by~\cite[Proposition~2.2]{Onsager_final} we have the following estimates:
\begin{align*}
    \norm{v_\ell-v_q}_0&\lesssim \delta_{q+1}^{1/2}\lambda_q^{-\alpha}, \\
    \norm{v_\ell}_{N+1}&\lesssim \delta_q^{1/2}\lambda_q\,\ell^{-N}\hspace{11pt} \forall N\geq 0, \\
    \|\mathring{R}_\ell\|_{N+\alpha}&\lesssim \delta_{q+1}\ell^{-N+3\alpha} \quad \forall N\geq 0,\\
    \abs{\int_\Omega \abs{v_q}^2-\abs{v_\ell}^2\,dx}&\lesssim \delta_{q+1}\ell^\alpha.
\end{align*}
Note that we have an extra factor $\ell^{2\alpha}$ in the estimate for the Reynolds stress in comparison to \cite{Onsager_final}. This comes from the extra factor $\lambda_q^{-3\alpha}$ in $\|\R{R}_q\|_0$ and from our definition of $\ell$. They cause an extra factor $\lambda_q^{-3\alpha}$ to appear in the estimate, which yields an extra factor $\ell^{2\alpha}$ by (\ref{tamaño ell}). 

Obtaining an extra factor $\ell^{2\alpha}$ (actually $\ell^\alpha$ would suffice) is our reason for modifying the inductive estimate of the Reynolds stress and the definition of $\ell$. We will use this extra factor to compensate a suboptimal estimate that we will be forced to use in \cref{section gluing in space}.

Once we have mollified the subsolution, we will glue it to the initial subsolution far from the turbulent zone. It will be convenient to divide $A_{q+1}\backslash A_q$ into several pieces because we will have to do several constructions in this region. We define
\begin{equation}
    \sigma = \frac{1}{5}(d_q-d_{q+1}),
\end{equation}
where $d_q$ was defined in (\ref{def d_q}). Using the elementary inequalities
\begin{equation}
	\label{elementary inequalities}
	2\pi\leq \frac{\lambda_q}{a^{b^q}}\leq 4\pi
\end{equation}
we deduce $\delta_{q+2}\gtrsim \lambda_q^{-2\beta b^2}$. By (\ref{condicion b}) we have $b^2<5/4$, so for sufficiently small $\alpha$ we have
\[d_q\gtrsim\lambda_q^{-1/11}.\]
Therefore,
\begin{equation}
	\label{tamaño sigma}
	\sigma^{-1}\lesssim \lambda_q^{1/11}.
\end{equation}
In particular, $\sigma\gg \ell$. For $1\leq j\leq 5$ we define
\[B_j\coloneqq \RR^3\backslash[A_q+B(0,j)]\]
By hypothesis, $(v_q,p_q,\mathring{R}_q)$ equals $(v_0,p_0,\mathring{R}_0)$ outside $A_q\times[0,T]$. Hence, it follows from (\ref{est v no cambia con el rescalado})$-$(\ref{est R no cambia con el rescalado}) that
\begin{align}
    \norm{v_\ell-v_0}_{N;B_1}+\norm{\partial_t v_\ell-\partial_t v_0}_{N;B_1}&\lesssim \ell^2, \label{mollifying good zone 1}\\
    \norm{p_\ell-p_0}_{N;B_1}&\lesssim \ell^2, \\
    \|\mathring{R}_\ell-\mathring{R}_0\|_{N;B_1}&\lesssim \ell^2. \label{mollifying good zone 3}
\end{align}
Thus, both subsolutions are very close in this region, which makes gluing them much easier. Taking into account (\ref{mollifying good zone 1}), it follows from \cref{invertir matrices} that there exists a potential $A\in C^\infty(B_1\times[0,T],\mathcal{A}^3)$ such that $\Div A=v_\ell-v_0$ in $B_1\times[0,T]$ and
\[\norm{A}_{N+1+\alpha;B_1}+\norm{\partial_t A}_{N+1+\alpha;B_1}\lesssim \ell^2\qquad \forall N\geq 0.\]
Therefore, the matrices $S_1$ and $M$ that appear in \cref{pegado subsoluciones} satisfy the estimates $\norm{S_1}_{N+\alpha}+\norm{M}_{N+\alpha}\lesssim \ell^2$. We introduce this estimates into \cref{pegado subsoluciones} to perform a gluing in the region $B_1\backslash B_2$. Since $\sigma\gg \ell$ we may essentially ignore any factor coming from the derivatives of the cutoff by absorbing it into the $\ell$ factor. Carrying out the rest of the construction of \cref{pegado subsoluciones}, we conclude that there exists a smooth subsolution such that
\begin{equation}
	\label{tvell igual inicial}
	(\T{v}_\ell,\T{p}_\ell,\R{\T{R}}_\ell)(x,t)=\begin{cases} (v_\ell,p_\ell,\R{R}_\ell)(x,t) \qquad \text{if } x\in A_q+B(0,\sigma), \\
		(v_0,p_0,\R{R}_0)(x,t) \hspace{18pt} \text{if } x\in B_2\end{cases}
\end{equation}
and satisfies the estimates
\begin{align}
    \norm{\T{v}_\ell-v_q}_0&\lesssim \delta_{q+1}^{1/2}\lambda_q^{-\alpha}, \label{est moll 1}\\
    \norm{\T{v}_\ell}_{N+1}&\lesssim \delta_q^{1/2}\lambda_q\,\ell^{-N}\hspace{15pt} \forall N\geq 0, \label{est moll 2}\\
    \|\R{\T{R}}_\ell\|_{N+\alpha}&\lesssim \delta_{q+1}\ell^{-N+3\alpha} \quad \forall N\geq 0,\label{est moll 3}\\
    \abs{\int_\Omega \abs{\T{v}_q}^2-\abs{v_\ell}^2\,dx}&\lesssim \delta_{q+1}\ell^\alpha. \label{est moll 4}
\end{align}
That is, we have obtained a new subsolution satisfying the same estimates as $(v_\ell,p_\ell,\R{R}_\ell)$ but with the additional property of being equal to the initial subsolution far from the turbulent zone.

Although many constructions in~\cite{Onsager_final} work in $\RR^3$ with very little or no modification, it is more convenient to work with periodic fields so that we may use results from~\cite{Onsager_final} directly. Note that the subsolution that we have just obtained is supported in $(0,1)^3\times[0,T]$ because it equals the initial subsolution far form the turbulent zone. This allows us to consider its periodic extension to $\RR^3\slash \ZZ^3$, which we denote the same. From now on, we consider that we are working in this setting.

\subsection{Overview of the gluing in space} \label{overview gluing space}
The key idea introduced by Isett in~\cite{Isett2018} is to glue \textit{in time} exact solutions of the Euler equations to obtain a new subsolution $(\overline{v}_q,\overline{p}_q,\mathring{\overline{R}}_q)$ such that $\overline{v}_q$ is close to $v_q$ but $\mathring{\overline{R}}_q$ is supported in a series of disjoint temporal regions of the appropriate length. This allows the use of Mikado flows, leading to the optimal regularity $C^\beta$ for  any $\beta<1/3$ (as in Onsager's conjecture).

More specifically, we define the length
\begin{equation}
\label{def tau}
    \tau_q=\frac{\ell^{2\alpha}}{\delta_q^{1/2}\lambda_q}
\end{equation}
and we consider the smooth solutions of the Euler equations
\begin{equation}
	\label{ec vi}
	\begin{cases}
		\partial_t v_i+\Div(v_i\otimes v_i)+\nabla p_i=0,\\
		\Div v_i=0,\\
		v_i(\cdot,t_i)= \T{v}_\ell(\cdot,t_i),
	\end{cases}
\end{equation}
where $t_i\coloneqq i\tau_q$. We will see that they are defined in the time interval $[t_i-\tau_q,t_i+\tau_q]$. The pressure is recovered as the unique solution to the equation $-\Delta p_i=\tr(\nabla v_i\nabla v_i)$ with the normalization
\begin{equation}
    \int_{\TT^3}p_i(x,t)\,dx=\int_{\TT^3}\T{p}_\ell(x,t)\,dx.
\end{equation}
We will see that these solutions remain sufficiently close to $\T{v}_\ell$ in their respective intervals. Following Isett's ideas, we would like to glue them in time to obtain a subsolution in the whole interval $[0,T]$. The velocity field would remain sufficiently close to $\T{v}_\ell$ while the Reynolds stress would be localized to the intersection of the consecutive time intervals.

However, even if $\mathring{R}_\ell$ is well localized, the solutions $v_i$ will immediately differ from $\widetilde{v}_\ell$ in the whole space. Furthermore, different solutions $v_i$, $v_{i+1}$ will also differ in the whole space during the intersection of their temporal domains. If we tried to apply Isett's procedure to them we would obtain a Reynolds stress that spreads throughout the whole space. This is not suitable for our purposes, so we must modify Isett's approach.

What we will do is to glue \textit{in space} the exact solutions $v_i$ to $\widetilde{v}_\ell$ in the region where $\R{\T{R}}_\ell$ is small, obtaining subsolutions $(\widetilde{v}_i,\widetilde{p}_i,\mathring{\widetilde{R}}_i)$. The Reynolds stress will no longer be 0, but it will so small that we may ignore it in the current iteration. Since these subsolutions will coincide far from the turbulent zone, we will be able to glue them in time while keeping the Reynolds stress localized.

The actual process is not so simple because the difference between the exact solutions $(v_i,p_i)$ and the subsolution $(\T{v}_\ell,\T{p}_\ell,\R{\T{R}}_\ell)$ is too big, leading to an unacceptably large error if we try to glue them. Our approach consist in producing a series of intermediate subsolutions that act as a sort of interpolation between them. Instead of a single gluing we perform a big number of them, going from $v_i$ to $\T{v}_\ell$ far from the turbulent zone. The difference between two consecutive intermediate subsolutions will be very small so that the error introduced in each of these middle gluings is small.

At the end of this process we will obtain subsolutions $(\T{v}_i,\T{p}_i,\R{\T{R}}_i)$ such that
\begin{align}
	(\T{v}_i,\T{p}_i,\R{\T{R}}_i)&=(v_0,p_0,\R{\T{R}}_0) \text{ in } B_3\times[t_i-\tau_q,t_i+\tau_q], \label{vi igual inicial}\\
    \T{v}_i(t_i,\cdot)&=\T{v}_\ell(t_i,\cdot). \label{vi en ti}
\end{align}
In addition, for $\abs{t-t_i}\leq \tau_q$ and any $N\geq 0$ we will have the following estimates:
    \begin{align}
    	\|\R{\T{R}}_i\|_0&\leq \frac{1}{2}\delta_{q+2}\lambda_{q+1}^{-6\alpha}, \label{estimate vi 1}\\
        \norm{\T{v}_i-\T{v}_\ell}_{N+\alpha}&\lesssim \tau_q\delta_{q+1}\ell^{-N-1+\alpha}, \label{estimate vi 2}\\
        \norm{D_{t,\ell}(\T{v}_i-\T{v}_\ell)}_{N+\alpha}&\lesssim \delta_{q+1}\ell^{-N-1+\alpha}, \label{estimate vi 3}
    \end{align}
where we write
\[D_{t,\ell}\coloneqq \partial_t+\T{v}_\ell\cdot \nabla\]
for the transport derivative. Furthermore, there exist smooth vector potentials $\T{z}_i$ defined in $[t_i-\tau_q,t_i+\tau_q]$ such that $\T{v}_i=\curl \T{z}_i$ and at the intersection of two intervals we have:
    \begin{align}
        \norm{\T{z}_i-\T{z}_{i+1}}_{N+\alpha}&\lesssim \tau_q\delta_{q+1}\ell^{-N+\alpha}, \label{estimate zi 1} \\
        \|D_{t,\ell}(\T{z}_i-\T{z}_{i+1})\|_{N+\alpha}&\lesssim\delta_{q+1}\ell^{-N+\alpha}. \label{estimate zi 2}
    \end{align}
These estimates are completely analogous to the ones in \cite[Proposition 3.3, Proposition 3.4]{Onsager_final} but we have the additional benefit of the fields being equal to the initial subsolution far from the turbulent zone. We do have to pay a price because now we have subsolutions instead of solutions of the Euler equations. Nevertheless, the errors $\R{\T{R}}_i$ are so small that we may ignore them until the $(q+1)$-th iteration.

All of the steps summarized here will be discussed in full detail in \cref{section gluing in space}, where we derive the claimed estimates.

\subsection{Overview of the gluing in time}\label{overview gluing time}
Once we have our subsolutions $(\T{v}_i,\T{p}_i,\R{\T{R}}_i)$ we will glue them in time. The goal is to obtain a subsolution defined in all $[0,T]$ that remains close to $\widetilde{v}_\ell$ but in which most of the Reynolds stress is localized to temporally disjoint regions of the appropriate length. This will allow us to correct the error in each region separately using Mikado flows. Let 
\[t_i\coloneqq i\tau_q, \qquad I_i\coloneqq\left[t_i+\frac{1}{3}\tau_q, t_i+\frac{2}{3}\tau_q\right]\cap[0,T],\]
\[J_i\coloneqq\left(t_i-\frac{1}{3}\tau_q,t_i+\frac{1}{3}\tau_q\right)\cap[0,T].\]
Note that $\{I_j,J_i\}$ is a pairwise disjoint decomposition of $[0,T]$. We choose a smooth partition of unity $\{\chi_i\}$ such that:
\begin{itemize}
	\item $\sum_i \chi_i=1$.
	\item $\supp \chi_i\cap \supp \chi_{i+2}=\varnothing$. Furthermore, 
	\begin{equation}
		\begin{aligned}
			\supp \chi_i\subset \left(t_i-\frac{2}{3}\tau_q,t_i+\frac{2}{3}\tau_q\right),\\
			\chi_i(t)=1 \qquad \forall t\in J_i.
		\end{aligned}
	\end{equation}
	\item For any $i$ and $N\geq 0$ we have
	\begin{equation}
		\norm{\partial_t^N\chi_i}_0\lesssim \tau_q^{-N}.
	\end{equation}
\end{itemize} 
We define
\[\overline{v}_q\coloneqq \sum_i \chi_i\T{v}_i, \qquad \overline{p}^{(1)}_q\coloneqq \sum_i \chi_i\T{p}_i, \qquad \overline{\mathcal R}^{(1)}_q\coloneqq \sum_i\chi_i\R{\T{R}}_i.\]
It is clear that $\overline{v}_q$ is divergence-free and it equals $v_0$ in $B_3\times[0,T]$ because of (\ref{vi igual inicial}). In addition, it inherits the estimates of $\widetilde{v}_i$:
\begin{proposition}
	The velocity field $\overline{v}_q$ satisfies the following estimates:
	\begin{align}
		\norm{\overline{v}_q-\widetilde{v}_\ell}_\alpha &\lesssim \delta_{q+1}^{1/2}\ell^{\alpha}, \label{est vbarra 1}\\
		\norm{\overline{v}_q-\widetilde{v}_\ell}_{N+\alpha} &\lesssim \tau_q\delta_{q+1}\ell^{-1-N+\alpha}, \label{est vbarra 2}\\
		\norm{\overline{v}_q}_{1+N}&\lesssim \delta_q^{1/2}\lambda_q\ell^{-N}, \label{est vbarra 3}\\
		\abs{\int_\Omega \abs{\overline{v}_q}^2-\abs{\T{v}_\ell}^2}&\lesssim \delta_{q+1}\ell^{\alpha},
	\end{align}
	for all $N\geq 0$.
\end{proposition}
The proof is exactly the same as the proof of \cite[Proposition 4.3, Proposition 4.5]{Onsager_final} because our fields $\widetilde{v}_i$ satisfy completely analogous estimates to the solutions $v_i$ in~\cite{Onsager_final}. 

We conclude that the new velocity $\overline{v}_q$ equals $v_0$ far from the turbulent zone, it is close to $\widetilde{v}_\ell$ and satisfies suitable bounds. Nevertheless, we must check if it leads to a subsolution. If $t\in J_i$, then in a neighborhood of $t$ we have $\chi_i=1$ while the rest of the cutoffs vanish. Thus, for all $t\in J_i$ we have
\[\overline{v}_q= \T{v}_i, \qquad \overline{p}^{(1)}_q=  \T{p}_i, \qquad \overline{\mathcal R}^{(1)}_q=\R{\T{R}}_i.\]
Since $(\T{v}_i,\T{p}_i,\R{\T{R}}_i)$ is a subsolution, we have
\[\partial_t\overline{v}_q+\Div(\overline{v}_q\otimes\overline{v}_q)+\nabla \overline{p}^{(1)}_q=\Div\overline{\mathcal R}^{(1)}_q\]
On the other hand, if $t\in I_i$, then $\chi_j=0$ for $j\neq i,i+1$ and $\chi_i+\chi_{i+1}=1$. Hence, on $I_i$ we have
\[\overline{v}_q=\chi_i\T{v}_i+\chi_{i+1}\T{v}_{i+1}, \qquad \overline{p}^{(1)}_q=\chi_i\T{p}_i+\chi_{i+1}\T{p}_{i+1}, \qquad \overline{\mathcal R}^{(1)}_q=\chi_i\R{\T{R}}_i+\chi_{i+1}\R{\T{R}}_{i+1}.\]
After a tedious computation we obtain
\begin{align}
	\partial_t\overline{v}_q&+\Div(\overline{v}_q\otimes\overline{v}_q)+\nabla \overline{p}^{(1)}_q-\Div\overline{\mathcal R}^{(1)}_q= \label{ec velocidad pegado} \\
	&=\partial_t\chi_i(\T{v}_i-\T{v}_{i+1})-\chi_i(1-\chi_i)\Div((\T{v}_i-\T{v}_{i+1})\otimes (\T{v}_i-\T{v}_{i+1})), \nonumber
\end{align}
where we have used the fact that $(\T{v}_i,\T{p}_i,\R{\T{R}}_i)$ are subsolutions. 

In conclusion, for $t\in J_i$ the triplet $(\overline{v}_q,\overline{p}_q^{(1)},\overline{\mathcal R}^{(1)}_q)$ is trivially a subsolution, whereas for $t\in I_i$ it suffices to express the right-hand side of \cref{ec velocidad pegado} as the divergence of a symmetric matrix. However, we must do this carefully because it is very important to keep under control the spatial support of the Reynolds stress. Indeed, since $\overline{\mathcal R}^{(1)}_q$ equals $\R{R}_0$ outside $A_{q+1}\times[0,T]$, any perturbation in the Reynolds stress must be contained within $A_{q+1}\times[0,T]$ in order to satisfy the inductive property (\ref{inductive 1}). Fortunately, we will be able to achieve this because the right-hand side of \cref{ec velocidad pegado} is supported in $A_q+B(0,3\sigma)$ due to (\ref{vi igual inicial}). We see that performing the gluing in space in the previous stage is crucial.

The logical approach would be to apply \cref{invertir divergencia matrices} to the right-hand side of (\ref{ec velocidad pegado}) and define the new Reynolds stress as the sum of the obtained matrix and 
$\overline{\mathcal R}^{(1)}_q$. Unfortunately, we cannot simply do that because we cannot obtain the necessary estimates for the transport derivative from \cref{invertir divergencia matrices}. 

Nevertheless, this difficulty can be solved. In \cref{section gluing time} we will find smooth symmetric matrices $\overline{\mathcal R}^{(2)}_q,\overline{\mathcal R}^{(3)}_q$ solving the equation
\[\Div\left(\overline{\mathcal R}^{(2)}_q+\overline{\mathcal R}^{(3)}_q\right)=\partial_t\chi_i(\T{v}_i-\T{v}_{i+1})\]
for $t\in I_i$. We set $\overline{\mathcal R}^{(2)}_q,\overline{\mathcal R}^{(3)}_q=0$ for $t\notin \bigcup_{i}I_i$. Since the source term is supported in $A_q+B(0,3\sigma)$ due to (\ref{vi igual inicial}), we will be able to choose $\overline{\mathcal R}^{(2)}_q,\overline{\mathcal R}^{(3)}_q$ supported in $A_q+B(0,4\sigma)$. The motivation for each matrix is the following:
\begin{itemize}
	\item $\overline{\mathcal R}^{(2)}_q$ solves the equation except for a small error. We have good bounds for the derivative of the material derivative.
	\item $\overline{\mathcal R}^{(3)}_q$ corrects the errors introduced when fixing the support of $\overline{\mathcal R}^{(2)}_q$. We do not have good bounds for its material derivative, but its $C^0$-norm is very small.
\end{itemize}
Using these auxiliary matrices, we define
\begin{align}
	\Rbar{1}&\coloneqq \overline{\mathcal R}^{(1)}_q+\overline{\mathcal R}^{(3)}_q-\frac{1}{3}\tr{\overline{\mathcal R}^{(3)}_q}\Id, \label{def Rbar1}\\
	\Rbar{2}&\coloneqq \overline{\mathcal R}^{(2)}_q-\chi_i(1-\chi_i)(\T{v}_i-\T{v}_{i+1})\mathring{\otimes} (\T{v}_i-\T{v}_{i+1})-\frac{1}{3}\tr{\overline{\mathcal R}^{(2)}_q}\Id, \label{def Rbar2}\\
	\R{\olsi{R}}_q&\coloneqq \Rbar{1}+\Rbar{2}, \\
	\overline{p}_q&\coloneqq \overline{p}_q^{(1)}-\chi_i(1-\chi_i)\abs{\T{v}_i-\T{v}_{i+1}}^2-\frac{1}{3}\tr\left(\overline{\mathcal R}^{(2)}_q+\overline{\mathcal R}^{(3)}_q\right).
\end{align}
By construction of $\overline{\mathcal R}^{(2)}_q$ and $\overline{\mathcal R}^{(3)}_q$, we see that $(\overline{v}_q,\overline{p}_q,\R{\olsi{R}}_q)$ is a smooth subsolution. Furthermore, we have
\[(\overline{v}_q,\overline{p}_q,\R{\olsi{R}}_q)=(v_0,p_0,\R{R}_0) \text{ in } B_4\times[0,T]\]
because of (\ref{vi igual inicial}) and the fact that $\overline{\mathcal R}^{(2)}_q, \overline{\mathcal R}^{(3)}_q$ are supported in $A_q+B(0,4\sigma)$. We emphasize again the importance of performing a gluing in space to use $\T{v}_i$ instead of the solutions $v_i$. Otherwise, we would have no control on the Reynolds stress because $v_i-v_{i+1}$  will in general be spread throughout the whole space.

We summarize the facts about the new Reynolds stress that we will prove in \cref{section gluing time}:
\begin{proposition}
	The smooth symmetric matrices $\Rbar{1},\Rbar{2}$ satisfy
	\begin{align}
		\|\Rbar{1}\|_0&\leq \frac{3}{4}\delta_{q+2}\lambda_{q+1}^{-3\alpha}, \label{estimate Rbar1 final}\\
		\supp \Rbar{2}&\subset [A_q+B(0,4\sigma)]\times\bigcup_i I_i, \\
		\|\Rbar{2}\|_N&\lesssim \delta_{q+1}\ell^{-N+\alpha}, \label{estimate Rbar2 final}\\
		\|(\partial_t+\overline{v}_q\cdot\nabla)\Rbar{2}\|_N&\lesssim\delta_{q+1}\delta_q^{1/2}\lambda_q\ell^{-N-\alpha}. \label{estimate Rbar2 mat der final}
	\end{align}
\end{proposition}

Comparing (\ref{estimate Rbar2 final}) with the analogous estimate in \cite{Onsager_final}, we see that we estimate the $C^N$-norm instead of the $C^{N+\alpha}$-norm. This difference is immaterial because they merely use the $C^{N+\alpha}$-norm to estimate the $C^N$-norm, which leads to the same bound as (\ref{estimate Rbar2 final}).

In conclusion, $\Rbar{1}$ is so small that we may ignore it for the present iteration, whereas $\Rbar{2}$ is big but it is supported in temporally disjoint regions and it satisfies good estimates. The next stage of the process is aimed at correcting $\Rbar{2}$ by means of highly oscillatory perturbations.

All of the steps summarized here will be discussed in full detail in \cref{section gluing time}, where we derive the claimed estimates.

\subsection{Overview of the perturbation step.}
We have localized most of the Reynolds stress, that is, $\Rbar{2}$, to small disjoint temporal regions but to reduce it we must resort to convex integration. 

First of all, we fix a cutoff $\phi_q\in C^\infty(\RR^3,[0,1])$ that equals $1$ in $A_q+\overline{B}(0,4\sigma)$ and whose support is contained in $A_q+B(0,5\sigma)$. In particular, $\phi_q\equiv 1$ on the support of $\Rbar{2}(\cdot,t)$ for all $t\in [0,T]$.

We follow the construction of \cite{Onsager_final} with Mikado flows, but there are some differences:
\begin{itemize}
	\item We control the support of the perturbation by introducing the cutoff $\phi_q$.
	\item To obtain the desired energy, we must use a slightly different normalization coefficient for the perturbation to account for the presence of the cutoff $\phi_q$ when integrating.
	\item We construct the new Reynolds stress using \cref{invertir divergencia matrices} so that we have control on its support. To apply this lemma we must introduce a minor correction $w_L$ to ensure that the perturbation has vanishing angular momentum.
\end{itemize}
Since $\norm{\phi_q}_N$ are much smaller than the $C^N$-norms of the other maps involved, the presence of the cutoff does not affect the estimates. In addition, $w_L$ will be negligible because its size is determined by an integral quantity, which is very small for a highly oscillating perturbation.

The form of the perturbation $\T{w}_{q+1}= v_{q+1}-v_q$ is $\T{w}_{q+1}=w_0+w_c+w_L$, where $w_0$ is the main perturbation term and it is used to cancel the Reynolds stress, $w_c$ is a small correction to ensure that the perturbation is divergence-free and $w_L$ is a tiny correction that ensures that the perturbation has vanishing angular momentum. This term is not present in \cite{Onsager_final} because they do not control the support of the Reynolds stress, so it suffices to use $w_{q+1}\coloneqq w_0+w_c$.

At the end of the process we will obtain a new subsolution $(v_{q+1},p_{q+1},\R{R}_{q+1})$ that equals $(v_0,p_0,\R{R}_0)$ outside $A_{q+1}\times[0,T]$ and satisfying the following estimates:
\begin{align}
	\norm{v_{q+1}\overline{v}_q}_0+\frac{1}{\lambda_{q+1}}\norm{v_{q+1}\overline{v}_q}_1&\leq \frac{3}{4}M\delta_{q+1}^{1/2}, \label{final iteración 1}\\
	\|\R{R}_{q+1}\|_0&\leq \delta_{q+2}\lambda_{q+1}^{-6\alpha},\label{final iteración 2} \\
	\delta_{q+2}\lambda_{q+1}^{-7\alpha}\leq e(t)-\int_\Omega\abs{v_{q+1}(x,t)}^2dx&\leq \delta_{q+1}. \label{final iteración 3}
\end{align}

\subsection{Proof of \cref{prop iteracion integracion convexa}}\label{SS.prupro}
Finally, we are ready to complete the proof of the proposition. The estimate (\ref{estimación cambio iteración}) is a consequence of (\ref{est moll 1}), (\ref{est moll 2}), (\ref{est vbarra 1}), (\ref{est vbarra 3})  and (\ref{final iteración 1}):
\[\norm{v_{q+1}-v_q}_0+\lambda_{q+1}^{-1}\norm{v_{q+1}-v_q}_1\leq \frac{3}{4}M\delta_{q+1}^{1/2}+C\delta_{q+1}^{1/2}\ell^\alpha+C\delta_q^{1/2}\lambda_q\lambda_{q+1}^{-1},\]
where the constant $C$ depends on $\alpha$, $\beta$, $\Omega$ and $(v_0,p_0,\R{R}_0)$, but not on $a$, $b$ or $q$. Thus, for any fixed $b$ (\ref{estimación cambio iteración}) holds for sufficiently large $a$. Regarding (\ref{inductive 3}), we use the inequality at level $q$ to get
\[\norm{v_{q+1}}_1\leq M\delta_q^{1/2}\lambda_q+\frac{3}{4}M\delta_{q+1}^{1/2}+C\delta_{q+1}^{1/2}\ell^\alpha\lambda_{q+1}+C\delta_q^{1/2}\lambda_q.\]
Hence, if we choose $a$ large enough we obtain (\ref{inductive 3}). Finally, (\ref{inductive v0}) follows from
\[\norm{v_{q+1}}_0\leq \norm{v_q}_0+\norm{v_{q+1}-v_q}_0\leq 1-\delta_q^{1/2}+M\delta_{q+1}^{1/2}.\]
The inductive hypotheses (\ref{inductive 1}), (\ref{inductive 2}) and (\ref{inductive 5}) were obtained in the perturbation step, and we are done. \hfill\qedsymbol

\section{Gluing in space}\label{section gluing in space}

In this section we develop the second step of the proof of Proposition~\ref{prop iteracion integracion convexa}, which was summarized in \cref{overview gluing space}. We do it in three subsections. 

\subsection{Interpolating sequence}
We begin by describing the intermediate subsolutions that we will use. First of all, we recall the following local existence result. It is standard, but we provide a proof for the sake of completeness.
\begin{proposition}
	\label{classical existence}
    For any $\alpha>0$ there exists a constant $c(\alpha)>0$ with the following property. Given any $C^\infty$ initial data $u_0\in H^3(\RR^3)$, any $C^\infty$~force $f\in L^1_{\mathrm{loc}}(\RR,H^3(\RR^3))$, let us fix a constant $T>0$ such that
    \begin{equation}
    T\norm{u_0}_{1+\alpha}+T^2\norm{f}_{1+\alpha}\leq c(\alpha).
    \label{condición T Euler}
    \end{equation} 
    Then there exists a unique solution $u\in C^\infty(\RR^3\times[-T,T])\cap C([-T,T],H^3(\RR^3))$ to the Euler equation
    \[\begin{cases}
    \partial_t u+\Div(u\otimes u)+\nabla p=f,\qquad
    \Div u=0,\\
    u(\cdot,0)=u_0.
    \end{cases}\]
    Moreover, $u$ obeys the bounds
    \[\norm{u}_{N+\alpha}\lesssim \norm{u_0}_{N+\alpha}+T\norm{f}_{N+\alpha}\]
    for all $N\geq 1$, where the implicit constant depends on $N$ and $\alpha>0$.
\end{proposition}
\begin{proof}
It is classical~\cite{Swann,Kato,Temam} that the 3d Euler equation is locally wellposed on $H^3(\RR^3)$, and that the solution, which is defined a priori for some time $T=T(\|u_0\|_{H^3(\RR^3)})>0$, can be continued (and stay smooth, provided that $u_0\in C^\infty$) as long as the norm $\|u(t)\|_{H^3(\RR^3)}$ remains bounded. Furthermore, the weak Beale--Kato--Majda criterion shows~\cite{BKM} that this norm is controlled by $\|\nabla u\|_{L^1L^\infty}$. More precisely, one has
\begin{multline*}
    \|u(t)\|_{H^3(\RR^3)}\leq \|u_0\|_{H^3(\RR^3)}+C \int_{-|t|}^{|t|}\|f(\tau)\|_{H^3(\RR^3)}\, d\tau\\
    + C \int_{-|t|}^{|t|}\|\nabla u(\tau)\|_{L^\infty(\RR^3)}\,\|f(\tau)\|_{H^3(\RR^3)}\,e^{C \int_{-|\tau|}^{|\tau|}\|f(\tau')\|_{H^3(\RR^3)}\, d\tau'}\, d\tau\,.
\end{multline*}

Therefore, we only need to provide a uniform a priori estimate for $\|\nabla u\|_{L^\infty}$ for all times $t< T$, where the maximum value $T>0$ is yet to be determined. To this end, note that for any multi-index $\theta$ with $\abs{\theta}=N$ we have
	\[\partial_t\partial^\theta u+u\cdot\nabla \partial^\theta u+[\partial^\theta, u\cdot \nabla]u+\nabla\partial^\theta p=\partial^\theta f.\]
	Since the pressure satisfies the equation $-\Delta p=\tr(\nabla u\nabla u)-\Div f$, it follows that
	\[\norm{\nabla \partial^\theta p}_\alpha\lesssim \norm{\tr(\nabla u\nabla u)}_{N-1+\alpha}+\norm{f}_{N+\alpha}\lesssim \norm{u}_{1+\alpha}\norm{u}_{N+\alpha}+\norm{f}_{N+\alpha}.\]
	Hence, we have
	\[\norm{(\partial_t\partial^\theta+u\cdot\nabla \partial^\theta )u}_\alpha\lesssim \norm{u}_{1+\alpha}\norm{u}_{N+\alpha}+\norm{f}_{N+\alpha}.\]
	Thus, by \cref{transporte}:
	\begin{equation}
		\norm{u(\cdot,t)}_{N+\alpha}\lesssim \norm{u_0}_{N+\alpha}+T\norm{f}_{N+\alpha}+\int_0^t\norm{u(\cdot,s)}_{1+\alpha}\norm{u(\cdot,s)}_{N+\alpha}ds.
		\label{normas Hölder solución usando transporte}
	\end{equation}
	
Specializing to the case $N=1$, we note that for all $|t|<T$ one has
\begin{equation}
\norm{u(\cdot,t)}_{1+\alpha}\leq \tilde{c}(\alpha)\left(\norm{u_0}_{1+\alpha}+T\norm{f}_{1+\alpha}\right)+\tilde{c}(\alpha)\int_0^t\norm{u(\cdot,s)}_{1+\alpha}^2ds,
\label{aux existencia Euler}
\end{equation}
for some constant $\tilde{c}(\alpha)>0$. We define the constant $c(\alpha)$ that appears in the statement of the proposition as $c(\alpha)\coloneqq [2\tilde{c}(\alpha)]^{-2}$ and then assume that $T>0$ satisfies (\ref{condición T Euler}). Let $y_0$ be the first term on the right-hand side of (\ref{aux existencia Euler}) and let $y(t)$ be the solution to the ODE $y'= \tilde{c}(\alpha)\hspace{1pt}y^2$ with $y(0)=y_0$. One finds
\[
\norm{u}_{1+\alpha}\leq y(t)= \frac{y_0}{1-\tilde{c}(\alpha)y_0t}.
\]
Since $y(t)\lesssim1$ for all $|t|<T$, by our choice of $T$, we conclude that the solution is well defined for $|t|< T$. Furthermore, inserting this estimate in (\ref{normas Hölder solución usando transporte}) and applying Grönwall's inequality, we obtain the desired H\"older bounds for $N>1$.
\end{proof}

To construct the desired sequence of subsolutions, we define
\begin{align}
	m&\coloneqq \left\lceil\lambda_q^{1/2}\right\rceil, \\
	\R{R}_i^k&\coloneqq \frac{k}{m}\R{\T{R}}_\ell.
\end{align}
For $i\geq 0$ and $0\leq k\leq m$ we consider the smooth solutions of the forced Euler equations
\begin{equation}
	\label{ec vik}
	\begin{cases}
		\partial_t v^k_i+\Div(v^k_i\otimes v^k_i)+\nabla p^k_i=\Div\R{R}_i^k,\\
		\Div v^k_i=0,\\
		v^k_i(\cdot,t_i)= \T{v}_\ell(\cdot,t_i).
	\end{cases}
\end{equation}
Thus, the triplet $(v_i^k,p_i^k,\R{R}_i^k)$ is a subsolution. The pressure is recovered as the unique solution to the equation \[-\Delta p^k_i=\tr(\nabla v^k_i\nabla v^k_i)-\Div\Div\R{R}_i^k\] 
with the normalization
\begin{equation*}
	\int_{\TT^3}p^k_i(x,t)\,dx=\int_{\TT^3}\T{p}_\ell(x,t)\,dx.
\end{equation*}
Note that (in their common interval of existence) the pair $(v_i^0,p^0_i)$ equals the solutions $(v_i,p_i)$ considered in (\ref{ec vi}), while the pair $(v_i^m,p^m_i)$ equals $(\T{v}_\ell,\T{p}_\ell)$. The rest of the subsolutions form a sort of interpolating sequence between them.

We claim that $(v_i^k,p^k_i)$ are defined in the time interval $[t_i-\tau_q,t_i+\tau_q]$, where $\tau_q$ was defined in (\ref{def tau}). Indeed, it follows from (\ref{est moll 2}) that
\[\tau_q\|\T{v}_\ell\|_{1+\alpha}\lesssim \tau_q\delta_{q}^{1/2}\lambda_q\ell^{-\alpha}\lesssim \ell^{\alpha}.\]
On the other hand, by (\ref{est moll 3}) we have \[\tau_q^2\|\R{\T{R}}_\ell\|_{2+\alpha}\lesssim \tau_q^2\delta_{q+1}\ell^{-2+3\alpha}\lesssim \left(\frac{\ell^{4\alpha}}{\delta_q\lambda_q^2}\right)\left(\delta_{q+1}\ell^{3\alpha}\frac{\delta_q\lambda_q^{2+6\alpha}}{\delta_{q+1}}\right)=\ell^{7\alpha}\lambda_q^{6\alpha}\leq \ell^\alpha.\]
Hence, for sufficiently large $a$ we have 
\[\tau_q\|\T{v}_\ell\|_{1+\alpha}+\tau_q^2\|\R{\T{R}}_\ell\|_{2+\alpha}<\frac{1}{2}.\]
By \cref{classical existence}, we conclude that solutions $(v^k_i,p^k_i)$ of the corresponding forced Euler equations exist in the claimed interval and they are unique. Furthermore, we have the following bounds:
\begin{corollary}
	\label{corol estimate vik 1}
	If $a$ is sufficiently large, for $\abs{t-t_i}\leq \tau_q$ and any $N\geq 1$ we have
	\begin{equation}
		\label{estimate vik 1}
		\norm{v_i^k}_{N+\alpha}\lesssim \delta_q^{1/2}\lambda_q\ell^{1-N-\alpha}\lesssim \tau_q^{-1}\ell^{1-N+\alpha}.
	\end{equation}
\end{corollary}
\begin{proof}
	It follows from \cref{classical existence} and estimates (\ref{est moll 2}) and (\ref{est moll 3}) that
	\[\norm{v_i^k}_{N+\alpha}\lesssim\norm{\T{v}_\ell(t_i)}_{N+\alpha}+\tau_q\|\R{\T{R}_\ell}\|_{N+1+\alpha}\lesssim\delta_q^{1/2}\lambda_q\ell^{1-N-\alpha}+\tau_q\delta_{q+1}\ell^{-N-1+\alpha}\]
	Using definitions (\ref{def ell}) and (\ref{def tau}) and the comparison (\ref{tamaño ell}), we see that
	\[\frac{\tau_q\delta_{q+1}\ell^{-2+\alpha}}{\delta_{q}^{1/2}\lambda_q}=\frac{\delta_{q+1}\ell^{-2+3\alpha}}{\delta_q\lambda_q^2}=\lambda_q^{3\alpha}\ell^{3\alpha}\leq 1,\]
	from which the first inequality follows. For the second one, we use (\ref{def tau}) again.
\end{proof}

\subsection{Estimates for the interpolating sequence}
Now that we have defined our subsolutions, we must ensure that they remain close to $\T{v}_\ell$ and to each other. Taking into account that $(\T{v}_\ell,\T{p}_\ell,\R{\T{R}}_\ell)$ and $(v_i^k,p_i^k,\R{R}_i^k)$ are subsolutions, we see that the difference satisfies the following transport equation:
\[\partial_t(\T{v}_\ell-v_i^k)+\T{v}_\ell\cdot\nabla(\T{v}_\ell-v_i^k)=(v_i^k-\T{v}_\ell)\cdot\nabla v_i^k-\nabla(\T{p}_\ell-p_i^k)+\left(1-\frac{k}{m}\right)\Div\R{\T{R}}_\ell.\]
Hence, the difference $\T{v}_\ell-v_i^k$ satisfies the same equation as $v_\ell-v_i$ in \cite{Onsager_final} except for a factor multiplying $\Div\R{\T{R}}_\ell$, but it is less or equal than 1. Furthermore, by (\ref{est moll 2}) and (\ref{estimate vik 1}) the fields $\T{v}_\ell$ and $v_i^k$ satisfy the same estimates as $v_\ell$ and $v_i$. Therefore, we may argue as in \cite[Proposition 3.3]{Onsager_final}, obtaining:
\begin{proposition}
	\label{prop estimates vik}
    If $a$ is sufficiently large, for $\abs{t-t_i}\leq \tau_q$ and any $N\geq 0$ we have
    \begin{align}
        \norm{v_i^k-\T{v}_\ell}_{N+\alpha}&\lesssim \tau_q\delta_{q+1}\ell^{-N-1+3\alpha}, \label{estimate vik 2}\\
        \norm{\nabla(\T{p}_\ell-p_i^k)}_{N+\alpha}&\lesssim \delta_{q+1}\ell^{-N-1+3\alpha},\\
        \norm{D_{t,\ell}(v_i^k-\T{v}_\ell)}_{N+\alpha}&\lesssim \delta_{q+1}\ell^{-N-1+3\alpha}, \label{estimate vik 4}
    \end{align}
    where we write
    \[D_{t,\ell}\coloneqq \partial_t+\T{v}_\ell\cdot \nabla\]
    for the transport derivative.
\end{proposition}
Note that our extra factor $\ell^{2\alpha}$ in (\ref{est moll 3}) is inherited by these estimates. Next, let us consider the vector potential associated to the field $v_i^k$:
\[
    z_i^k=\mathcal{B}\hspace{0.5pt}v_i^k\coloneqq (-\Delta)^{-1}\text{curl\,}v_i^k,
\]
where $\mathcal{B}$ is the Biot-Savart operator. We have
\[\Div z_i^k=0 \qquad \text{and} \qquad \curl ,z_i^k=v_i^k-\int_{\TT^3}v_i^k.\]
Since $(v_i^k,p_i^k,\frac{k}{m}\R{\T{R}}_\ell)$ is a subsolution, we have
\[\frac{d}{dt}\int_{\TT^3}v_i^k=-\int_{\TT^3}\left(\Div(v_i^k\otimes v_i^k)+\nabla p_i^k-\frac{k}{m}\Div\R{\T{R}}_\ell\right)=0.\]
On the other hand, the average of $v_i^k$ at time $t=t_i$ is the average of $\T{v}_\ell$ at that time, which vanishes because $\T{v}_\ell(\cdot,t_i)$ is divergence-free and its support is contained in $(0,1)^3$. Therefore, $v_i^k$ has zero mean and we have $\curl z_i^k=v_i^k$.

Since our fields satisfy the same estimates as the fields in \cite{Onsager_final}, we can again argue as in \cite[Proposition 3.4]{Onsager_final}, obtaining:
\begin{proposition}
	\label{prop estimates zik}
    For $\abs{t_i-\tau_q}\leq \tau_q$ and any $N\geq 0$ we have
    \begin{align}
        \norm{z_i^k-z_i^m}_{N+\alpha}&\lesssim \tau_q\delta_{q+1}\ell^{-N+3\alpha}, \label{estimate zik 1} \\
        \|D_{t,\ell}(z_i^k-z_i^m)\|_{N+\alpha}&\lesssim\delta_{q+1}\ell^{-N+3\alpha}, \label{estimate zik 2}
    \end{align}
    where $D_{t,\ell}=\partial_t+\T{v}_\ell\cdot\nabla$.
\end{proposition}
In summary, the difference $v_i^k-\T{v}_\ell$ satisfies the same stability estimates as the ones in \cite{Onsager_final} plus an additional factor $\ell^{2\alpha}$ due to the difference in the definition of $\ell$ and in the inductive estimate (\ref{inductive 2}). Furthermore, we will see that the difference between consecutive fields is much smaller, which will allow us to obtain suitable bounds for the gluing. We argue as in \cite[Proposition 3.3]{Onsager_final}. Subtracting the equation for each field and rearranging we obtain 

\begin{equation}
	\label{ec difference vik}
	\partial_t(v_i^{k+1}-v_i^k)+v_i^{k+1}\cdot\nabla(v_i^{k+1}-v_i^k)=(v_i^k-v_i^{k+1})\cdot \nabla v_i^k-\nabla(p_i^{k+1}-p_i^k)+\frac{1}{m}\Div\R{\T{R}}_\ell.
\end{equation}
Taking the divergence, we have
\[\Delta(p_i^{k+1}-p_i^k)=\Div[\nabla v_i^{k+1}(v_i^k-v_i^{k+1})]+\Div[\nabla v_{i}^k(v_i^k-v_i^{k+1})]+\frac{1}{m}\Div\Div\R{\T{R}}_\ell.\]
Since $(-\Delta)^{-1}\Div\Div$ is a Calderón-Zygmund operator, we obtain
\[\norm{\nabla(p_i^{k+1}-p_i^k)}_\alpha\lesssim \tau_q^{-1}\norm{v_i^k-v_i^{k+1}}_\alpha+\frac{1}{m}\delta_{q+1}\ell^{-1+\alpha}\]
where we have used (\ref{est moll 3}) and (\ref{estimate vik 1}). The additional factor $\ell^{2\alpha}$ is not needed here, so we just omit it. Inserting this estimate into \cref{ec difference vik} and using (\ref{est moll 3}) and (\ref{estimate vik 1}) again, we obtain
\[\norm{\partial_t(v_i^{k+1}-v_i^k)+v_i^{k+1}\cdot\nabla(v_i^{k+1}-v_i^k)}_\alpha\lesssim \frac{1}{m}\delta_{q+1}\ell^{-1+\alpha}+\tau_q^{-1}\norm{v_i^k-v_i^{k+1}}_\alpha.\]
Applying \cref{transporte} yields
\[\norm{(v_i^{k+1}-v_i^k)(\cdot,t)}_\alpha\lesssim \frac{1}{m}\abs{t-t_i}\delta_{q+1}\ell^{-1+\alpha}+\int_{t_i}^t\tau_q^{-1}\norm{(v_i^{k+1}-v_i^k)(\cdot,s)}_\alpha \,ds.\]
Using Grönwall's inequality and the assumption $\abs{t-t_i}\leq \tau_q$ we conclude
\[\norm{v_i^{k+1}-v_i^k}_\alpha\lesssim \frac{1}{m}\tau_q\delta_{q+1}\ell^{-1+\alpha}.\]
If we carry on arguing as in \cite[Proposition 3.3]{Onsager_final} we obtain the following higher-order estimates:
\begin{equation}
	\label{obtenido diferencia vik}
	\norm{v_i^{k+1}-v_i^k}_{N+\alpha}\lesssim \frac{1}{m}\tau_q\delta_{q+1}\ell^{-N-1+\alpha}.
\end{equation}
Let us use this bound to estimate the other fields. We may rewrite the equation for the pressure as
\begin{align}
	\label{ec pik difference}
	\Delta(p_i^{k+1}-p_i^k)=\Div\Div\bigg(&\frac{1}{2}(v_{i}^k+v_i^{k+1})\otimes(v_{i}^k-v_i^{k+1})\\&+\frac{1}{2}(v_{i}^k-v_i^{k+1})\otimes (v_{i}^k+v_i^{k+1})+\frac{1}{m}\R{\T{R}}_\ell\bigg) \nonumber
\end{align}
because
\[v_i^k\otimes v_i^k-v_i^{k+1}\otimes v_i^{k+1}=\frac{1}{2}(v_{i}^k+v_i^{k+1})\otimes(v_{i}^k-v_i^{k+1})+\frac{1}{2}(v_{i}^k-v_i^{k+1})\otimes (v_{i}^k+v_i^{k+1}).\]
Interpolating between (\ref{inductive 2}) and (\ref{inductive 3}) and between (\ref{est moll 1}) and (\ref{est moll 2}) we have
\begin{align*}
	\norm{v_{i}^k+v_i^{k+1}}_\alpha &\leq  2\norm{v_q}_\alpha+2\norm{v_q-\T{v}_\ell}_\alpha+\norm{v_{i}^k-\T{v}_\ell}_\alpha+\norm{v_i^{k+1}-\T{v}_\ell}_\alpha \\
	&\lesssim (\delta_q^{1/2}\lambda_q)^\alpha+ (\delta_{q+1}\lambda_q^{-\alpha})^{1-\alpha}(\delta_q^{1/2}\lambda_q)^\alpha+\frac{1}{m}\tau_q\delta_{q+1}\ell^{-1+\alpha} \lesssim \lambda_q^\alpha.
\end{align*}
For the higher-order bounds we use \cref{corol estimate vik 1}. Therefore, from (\ref{obtenido diferencia vik}) we conclude
\begin{align}
	\nonumber \norm{v_i^k\otimes v_i^k-v_i^{k+1}\otimes v_i^{k+1}}_\alpha&\lesssim \lambda_q^\alpha(m^{-1}\tau_q\delta_{q+1}\ell^{-N-1+\alpha})
	\\ \label{quadratic difference} &\hspace{10pt}+\sum_{j=1}^N(\tau_q^{-1}\ell^{1-j+\alpha})(m^{-1}\tau_q\delta_{q+1}\ell^{-(N-j)-1+\alpha}) \\ &\nonumber \lesssim\frac{1}{m}\tau_q\delta_{q+1}\ell^{-N-1}\lesssim \frac{1}{m}\delta_{q+1}^{1/2}\ell^{-N}.
\end{align}
Introducing this estimate and (\ref{est moll 3}) in \cref{ec pik difference}, we finally obtain
\begin{equation}
	\label{difference pik}
	\norm{p_i^{k+1}-p_i^k}_{N+\alpha}\lesssim\frac{1}{m}\tau_q\delta_{q+1}\ell^{-N-1}\lesssim \frac{1}{m}\delta_{q+1}^{1/2}\ell^{-N}.
\end{equation}
Let us now estimate the difference in the vector potentials. We recall the identity $\curl\curl=-\Delta+\nabla\Div$ and that $\Div z_i^k=0$. Hence, taking the $\curl$ in \cref{ec difference vik} and rearranging we arrive at
\begin{align}
	\label{ec zik difference}
	-\Delta\left[\partial_t(z_i^{k+1}-z_i^k)\right]=\curl\Div\bigg(&\frac{1}{2}(v_{i}^k+v_i^{k+1})\otimes(v_{i}^k-v_i^{k+1})\\&+\frac{1}{2}(v_{i}^k-v_i^{k+1})\otimes (v_{i}^k+v_i^{k+1})+\frac{1}{m}\R{\T{R}}_\ell\bigg). \nonumber
\end{align}
Reasoning as in the case of the pressure, we conclude
\begin{equation}
	\label{difference partial t zik}
	\norm{\partial_t(z_i^{k+1}-z_i^k)}_{N+\alpha}\lesssim\frac{1}{m}\tau_q\delta_{q+1}\ell^{-N-1}\lesssim \frac{1}{m}\delta_{q+1}^{1/2}\ell^{-N}.
\end{equation}
Since $z_i^k(\cdot,t_i)=z_i^{k+1}(\cdot,t_i)=\mathcal{B}\,\T{v}_\ell(\cdot,t_i)$, the difference vanishes at $t=t_i$. Using the assumption $\abs{t-t_i}\leq \tau_q$ we deduce
\begin{equation}
	\label{difference zik}
	\norm{z_i^{k+1}-z_i^k}_{N+\alpha}\lesssim\frac{1}{m}\tau_q^2\delta_{q+1}\ell^{-N-1}\lesssim \frac{1}{m}\tau_q\delta_{q+1}^{1/2}\ell^{-N}.
\end{equation}

\subsection{Gluing the interpolating sequence}
Now that we have the appropriate estimates, we will start gluing the subsolutions $(v_i^k,p_i^k,\R{R}_i^k)$ to one another to construct a subsolution that equals $(v_i,p_i,0)$ in $A_q+B(0,2\sigma)$ and $(v_0,p_0,\R{R}_0)$ in $B_3$. For the sake of clarity, we will do it inductively. Let
\[r\coloneqq \lambda_q^{-3/5}\]
and for $k\geq 0$ consider the sets
\begin{align*}
	\Omega^k&\coloneqq A_q+B(0,2\sigma+kr), \\ U^k&\coloneqq \{x\in \RR^3:2\sigma+kr<\dist(x,A_q)<2\sigma+(k+1)r\}
\end{align*}
and we fix smooth cutoff functions $\varphi^k\in C^\infty(\Omega^{k+1},[0,1])$ that equal 1 in a neighborhood of $\Omega^k$. By \cref{cutoff}, we may assume the bounds $\|\varphi^k\|_N\lesssim r^{-N}$.

We will construct a sequence of subsolutions $(\vtildei{k},\ptildei{k},\Rtildei{k})$ and potentials $\ztildei{k}$ with $\vtildei{k}=\curl \ztildei{k}$ such that
\begin{align}
	(\vtildei{k},\ptildei{k},\Rtildei{k},\ztildei{k})(x,t)&=(v_i^k,p_i^k,\R{R}_i^k,z_i^k)(x,t) \quad  \forall x\notin \Omega^k, \; \abs{t-t_i}\leq \tau_q, \label{vtildeik igual inicial}\\
	\vtildei{k}(\cdot,t_i)&=\T{v}_\ell(\cdot,t_i), \label{vtildeik en ti} 
\end{align}
Furthermore, for $\abs{t-t_i}\leq \tau_q$ and any $N\geq 0$ they satisfy
\begin{align}
	\|\Rtildei{k}\|_0&\leq \frac{1}{2}\delta_{q+2}\lambda_{q+1}^{-6\alpha}, \label{estimate Tvik 1}\\
	\norm{\vtildei{k}-\T{v}_\ell}_N&\lesssim \tau_q\delta_{q+1}\ell^{-N-1+3\alpha}, \label{estimate Tvik 2}\\
	\norm{D_{t,\ell}(\vtildei{k}-\T{v}_\ell)}_N&\lesssim \delta_{q+1}\ell^{-N-1+3\alpha}, \label{estimate Tvik 3} \\
	\norm{\ztildei{k}-z_i^m}_N&\lesssim \tau_q\delta_{q+1}\ell^{-N+3\alpha}, \label{estimate ztildeik 1} \\
	\|D_{t,\ell}(\ztildei{k}-z_i^m)\|_N&\lesssim\delta_{q+1}\ell^{-N+3\alpha}, \label{estimate ztildeik 2}
\end{align}
where we write
\[D_{t,\ell}\coloneqq \partial_t+\T{v}_\ell\cdot \nabla\]
for the transport derivative. 

If we could construct such a sequence, setting $(\T{v}_i,\T{p}_i,\R{\T{R}}_i,\T{z}_i)\coloneqq (\vtildei{m},\T{p}_i^{\hspace{1pt}m},\Rtildei{m},\ztildei{m})$ would yield the subsolution and potential claimed in \cref{overview gluing space}. Indeed, (\ref{vi en ti}) and (\ref{estimate vi 1}) are just (\ref{vtildeik en ti}) and (\ref{estimate Tvik 1}). The estimates (\ref{estimate vi 2}) and (\ref{estimate vi 3}) follow from (\ref{estimate Tvik 2}) and (\ref{estimate Tvik 3}) by interpolation, that is, we use (\ref{interpolation inequality}) to estimate the $C^{N+\alpha}$ seminorm using the $C^0$ and $C^{N+1}$ norms:
\begin{align*}
    \big[\vtildei{k}-\T{v}_\ell\big]_{N+\alpha}&\lesssim\norm{\vtildei{k}-\T{v}_\ell}_0^{1-\frac{N+\alpha}{N+1}}\norm{\vtildei{k}-\T{v}_\ell}_{N+1}^{\frac{N+\alpha}{N+1}}\\&\lesssim \left(\tau_q\delta_{q+1}\ell^{-1+3\alpha}\right)^{1-\frac{N+\alpha}{N+1}}\left(\tau_q\delta_{q+1}\ell^{-1-(N+1)+3\alpha}\right)^{\frac{N+\alpha}{N+1}} \\&\lesssim\tau_q\delta_{q+1}\ell^{-1-N+2\alpha}.
\end{align*}
This, combined with the estimate for the $C^N$ norm, yields (\ref{estimate vi 2}). Note that here we have obtained an extra factor $\ell^\alpha$, but we discard it because it will not be necessary. Estimate (\ref{estimate vi 1}) follows from (\ref{estimate Tvik 3}) by an analogous argument. Meanwhile, to obtain (\ref{estimate zi 1}) and (\ref{estimate zi 2}) from (\ref{estimate ztildeik 1}) and (\ref{estimate ztildeik 2}), we also need to apply the triangle inequality and use the fact that $z_{i+1}^m-z_i^m=0$. Recall that both vector potentials are just the restriction of $\mathcal{B}\hspace{0.5pt} \T{v}_\ell$ to their respective intervals. Finally, (\ref{vi igual inicial}) follows from (\ref{vtildeik igual inicial}) because $mr\lesssim \lambda_q^{-1/10}$, so by (\ref{tamaño sigma}) it must be smaller than $\sigma$ for sufficiently large $a$.

Let us then construct this sequence. We define the initial term as 
\[(\vtildei{0},\ptildei{0},\Rtildei{0},\ztildei{0})\coloneqq(v_i^0,p_i^0,\R{R}_i^0)=(v_i,p_i,0,z_i^0).\]
It follows from \cref{corol estimate vik 1}, \cref{prop estimates vik} and \cref{prop estimates zik} that this term satisfies \Crefrange{vtildeik igual inicial}{estimate ztildeik 2}. Next, let us suppose that we have defined the $k$-th term $(\vtildei{k},\ptildei{k},\Rtildei{k},\ztildei{k})$ satisfying these inductive hypotheses. We will construct the $(k+1)$-th term satisfying them, too. To do so, we will apply \cref{pegado subsoluciones} to glue $(\vtildei{k},\ptildei{k},\Rtildei{k})$ and $(v_i^{k+1},p_i^{k+1},\R{R}_i^{k+1})$ in the region $U^k$. Since these subsolutions are defined in the whole space, by \cref{condiciones triviales} we do not need to check the compatibility conditions (\ref{comp1}) and (\ref{comp2}).

Note that in \cref{pegado subsoluciones} we use skew-symmetric matrices instead of potential vectors because it is stated in any dimension $\dim\geq 2$. However, it is completely equivalent: given a potential vector $z$, we simply define $A_{ij}=\varepsilon_{ijk}z_k$, where $\varepsilon_{ijk}$ is the usual Levi-Civita symbol. It is easy to check that $A$ is skew-symmetric and $\curl z=\Div A$.

Hence, applying \cref{pegado subsoluciones} we obtain a subsolution satisfying:
\[(\vtildei{k+1},\ptildei{k+1},\Rtildei{k+1})(\cdot,t)=\begin{cases}(\vtildei{k},\ptildei{k},\Rtildei{k}) \hspace{40pt} \text{in } \Omega^k, \\[5pt] (v_i^{k+1},p_i^{k+1},\R{R}_i^{k+1}) \quad \text{outside } \Omega^{k+1}\end{cases}\]
for $\abs{t-t_i}\leq \tau_q$. In addition, there exists a smooth vector potential $\ztildei{k+1}$ with $\vtildei{k+1}=\curl \ztildei{k+1}$ and such that $\ztildei{k+1}(\cdot,t)\equiv z_i^{k+1}(\cdot,t)$ outside $\Omega^{k+1}$. Thus, (\ref{vtildeik igual inicial}) is satisfied. 

Furthermore, by definition of $v_i^{k+1}$ and the inductive hypothesis (\ref{vtildeik en ti}) we have
\[v_i^{k+1}(\cdot,t_i)=\T{v}_\ell(\cdot,t_i)=\vtildei{k}(\cdot,t_i).\]
In addition, by (\ref{vtildeik igual inicial}) we know that $\ztildei{k}$ equals $z_i^{k}$ outside $\Omega^k$. Since $z_i^{k}$ and $z_i^{k+1}$ both equal $\mathcal{B}\,\T{v}_\ell$ at time $t=t_i$, we see that the difference $z_i^{k+1}-\ztildei{k}$ vanishes outside $\Omega^k$ at $t=t_i$. Inspecting \cref{pegado subsoluciones} and replacing skew-symmetric matrices by the equivalent potential vectors, we see that
\[(\vtildei{k+1}-w_L)(\cdot,t_i)=\left[\varphi^k\,\vtildei{k}+(1-\varphi^k)v_i^{k+1}+\nabla\varphi^k\times(z_i^{k+1}-\ztildei{k})\right](\cdot,t_i)=\T{v}_\ell(\cdot,t_i).\]
Since only the time derivative of $w_L$ matters, we may assume that it vanishes at $t=t_i$, that is, we start the integration in (\ref{def Lij}) at $t=t_i$. Hence, (\ref{vtildeik en ti}) holds.

Concerning the estimates, we see in \cref{pegado subsoluciones} that we only need to consider the bounds in $U^k$, so by (\ref{vtildeik igual inicial}) we only need to study the difference between $(v_i^k,p_i^k,\R{R}_i^k,z_i^k)$ and $(v_i^{k+1},p_i^{k+1},\R{R}_i^{k+1},z_i^{k+1})$. Since
\[\|\R{R}_i^k-\R{R}_i^{k+1}\|_0=\frac{1}{m}\|\R{\T{R}}_\ell\|_0\lesssim\frac{1}{m}\delta_{q+1},\]
it follows from \cref{quadratic difference,difference pik,difference partial t zik} that the matrix $M$ that appears in \cref{pegado subsoluciones} satisfies
\[\norm{M}_{0;\hspace{0.3pt}U^k}\leq \frac{1}{m}\delta_{q+1}^{1/2}.\]
With this bound we may use the estimates from \cref{pegado subsoluciones}. Let us focus on the vector potential:
\begin{align}
	\norm{\ztildei{k+1}-(\varphi^k\ztildei{k}+(1-\varphi^k)z_i^{k+1})}_{N}&\lesssim \frac{1}{m}\tau_q\delta_{q+1}^{1/2}r^{-N},\\
	\norm{\partial_t(\ztildei{k+1}-\varphi^k\ztildei{k}-(1-\varphi^k)z_i^{k+1})}_{N}&\lesssim \frac{1}{m}\delta_{q+1}^{1/2}r^{-N}, \label{partialt correccion z}
\end{align}
To use this estimates, the following inequality will be useful:
\begin{equation}
	\label{desigualdad m}
	m\delta_{q+2}\delta_{q+1}^{-1/2}=\lambda_q^{1/2-2b^2\beta+b\beta}\geq \lambda_{q}^{1/20},
\end{equation}
where we have used that the exponent in the middle term is greater than $1/20$ for any $0<\beta<1/3$ and $1<b<11/10$. We compute
\begin{align*}
	\norm{\ztildei{k+1}-z_i^m}_N&\leq \norm{\ztildei{k+1}-(\varphi^k\ztildei{k}+(1-\varphi^k)z_i^{k+1})}_{N}\\&\hspace{12pt}+\norm{\varphi^k(\ztildei{k}-z_i^m)}_N+\norm{(1-\varphi^k)(z_i^{k+1}-z_i^m)}_N\\
	&\lesssim \frac{1}{m}\tau_q\delta_{q+1}^{1/2}r^{-N}+\sum_{j=0}^Nr^{-j}\left(\norm{\ztildei{k}-z_i^m}_{N-j}+\norm{z_i^{k+1}-z_i^m}_{N-j}\right) \\
	&\lesssim \frac{1}{m}\tau_q\delta_{q+1}^{1/2}r^{-N}+\tau_q\delta_{q+1}\ell^{-N+3\alpha}\lesssim \tau_q\delta_{q+1}\ell^{-N+3\alpha},
\end{align*}
where we have used (\ref{desigualdad m}) and we have assumed $\alpha$ to be sufficiently small. Hence, (\ref{estimate ztildeik 1}) holds for the $(k+1)$-th term. In addition, (\ref{estimate Tvik 2}) clearly follows from it.

Note that we have used (\ref{estimate zik 1}) to estimate $\norm{\cdot}_N$, that is, we lose an $\alpha$. This is clearly not optimal, but it cannot be avoided for $N=0$. Thus, we pay the prize of losing a factor $\ell^\alpha$ in the estimates for $\norm{\cdot}_{N+\alpha}$. To compensate this, we gain an extra factor $\ell^{2\alpha}$ from (\ref{inductive 2}) and (\ref{def ell}), which are different from their counterparts in \cite{Onsager_final}.

Let us now estimate the material derivative. By the triangle inequality:
\begin{align}
	\label{estimate material derivative correction z}
	\norm{D_{t,\ell}(\ztildei{k+1}-z_i^m)}_N&\leq \norm{\partial_t(\ztildei{k+1}-\varphi^k\ztildei{k}-(1-\varphi^k)z_i^{k+1})}_{N}\\ \nonumber&\hspace{12pt}+\norm{\T{v}_\ell\cdot\nabla(\ztildei{k+1}-\varphi^k\ztildei{k}-(1-\varphi^k)z_i^{k+1})}_{N}\\ \nonumber&\hspace{12pt}+\norm{D_{t,\ell}[\varphi^k(\ztildei{k}-z_i^m)]}_N\\ \nonumber&\hspace{12pt}+\norm{D_{t,\ell}[(1-\varphi^k)(v_i^{k+1}-z_i^m)]}_N.
\end{align}
The last two terms can be estimated in the same manner, so we just study one:
\begin{align*}
	\norm{D_{t,\ell}[\varphi^k(\ztildei{k}-z_i^m)]}_N&\lesssim \sum_{j=0}^N\big(\norm{\varphi^k}_{j+1}\norm{\ztildei{k}-z_i^m}_{N-j}\\&\hspace{35pt}+\norm{\varphi^k}_{j}\norm{D_{t,\ell}(\ztildei{k}-z_i^m)}_{N-j}\big) \\ 
	&\lesssim r^{-1}\tau_q\delta_{q+1}\ell^{-N+3\alpha}+\delta_{q+1}\ell^{-N+3\alpha}\lesssim \delta_{q+1}\ell^{-N+3\alpha},
\end{align*}
where we have used $\norm{\varphi^k}_N\lesssim r^{-N}\lesssim \ell^{-N}$ and $\tau_qr^{-1}\lesssim 1$. Next, taking into account that $\norm{\T{v}_\ell}_N\lesssim \ell^{-N}$, we have
\begin{align*}
	\|\T{v}_\ell\cdot\nabla(\ztildei{k+1}-&\varphi^k\ztildei{k}-(1-\varphi^k)z_i^{k+1})\|_{N}\lesssim \\
	&\lesssim\sum_{j=0}^N\ell^{-j}\norm{\ztildei{k+1}-\varphi^k\ztildei{k}-(1-\varphi^k)z_i^{k+1}}_{N+1-j}\\
	&\lesssim \frac{1}{mr}\tau_q\delta_{q+1}^{1/2}\ell^{-N}\lesssim \delta_{q+1}\ell^{-N+3\alpha},
\end{align*}
where we have used $\tau_qr^{-1}\lesssim 1$, the inequality (\ref{desigualdad m}) and we have assumed $\alpha$ to be sufficiently small. Using (\ref{partialt correccion z}) along with the same tricks, we obtain the same estimates for the first term in \cref{estimate material derivative correction z}, and we conclude that (\ref{estimate ztildeik 2}) holds for the $(k+1)$-th term.

To obtain (\ref{estimate Tvik 3}), we estimate the commutator $[D_{t,\ell},\curl]$. We fix an arbitrary vector field $u$ and we compute
\[\left(\curl(D_{t,\ell}u)-D_{t,\ell}\curl u\right)_i=\varepsilon_{ijk}\partial_j(\T{v}_\ell)_l\partial_l u_k.\]
Therefore, we have
\begin{align*}
	\norm{[D_{t,\ell},\curl](\ztildei{k+1}-z_i^m)}_N&\lesssim \sum_{j=0}^N\norm{\T{v}_\ell}_{j+1}\norm{\ztildei{k+1}-z_i^m}_{N-j+1}, \\
	&\lesssim \sum_{j=0}^N (\delta_q^{1/2}\lambda_q\ell^{-j})(\tau_q\delta_{q+1}\ell^{-N-1+j+3\alpha})\\
	&\lesssim \delta_q^{1/2}\lambda_q\tau_q\delta_{q+1}\ell^{-N-1+3\alpha}\lesssim \delta_{q+1}\ell^{-N-1+3\alpha}
\end{align*}
where we have used (\ref{est moll 2}) and that $\delta_q^{1/2}\lambda_q\tau_q=\ell^{2\alpha}\lesssim 1$ by definition of $\tau_q$. Hence, 
\begin{align*}
	\norm{D_{t,\ell}(\vtildei{k+1}-\T{v}_\ell)}_N&\leq \norm{\curl[D_{t,\ell}(\ztildei{k+1}-z_i^m)]}_N+\norm{[D_{t,\ell},\curl](\ztildei{k+1}-z_i^m)}_N \\
	&\lesssim \delta_{q+1}\ell^{-N-1+3\alpha}.
\end{align*}

Finally, let us consider the size of the new Reynolds stress. Taking into account that $\ztildei{k}$ equals $z_i^k$ in $U^k$, it follows from (\ref{difference zik}) and \cref{pegado subsoluciones} that
\[\norm{w}_0=\norm{\vtildei{k+1}-\varphi^k\vtildei{k}-(1-\varphi^k)v_i^{k+1}}_0\lesssim \tau_qr^{-1}\frac{1}{m}\delta_{q+1}^{1/2}\lesssim \frac{1}{m}\delta_{q+1}^{1/2},\]
where we have used $\tau_qr^{-1}\lesssim 1$. Taking into account that $\vtildei{k}$ equals $v_i^k$ in $U^k$, it follows from (\ref{obtenido diferencia vik}) and \cref{pegado subsoluciones} that
\[\|\Rtildei{k+1}-(\varphi^k\Rtildei{k}+(1-\varphi^k)\R{R}_i^{k+1})\|_0\lesssim \frac{1}{m}\delta_{q+1}^{1/2}.\]
Taking into account (\ref{desigualdad m}), we see that for sufficiently small $\alpha$ and sufficiently large $a$ we have 
\[\|\Rtildei{k+1}-(\varphi^k\Rtildei{k}+(1-\varphi^k)\R{R}_i^{k+1})\|_0\leq\frac{1}{4}\delta_{q+2}\lambda_{q+1}^{-3\alpha}.\]
If $x\in \Omega^k$, then $\varphi^k(x)=1$ and by (\ref{estimate Tvik 1}) we have
\[|\Rtildei{k+1}(x,t)|=|\Rtildei{k}(x,t)|\leq \frac{1}{2}\delta_{q+2}\lambda_{q+1}^{-3\alpha}.\]
On the other hand, if $x\notin \Omega^k$, then $\Rtildei{k}(x,t)=\frac{k}{m}\R{\T{R}}_\ell(x,t)$ by (\ref{vi igual inicial}). Furthermore, we also have $\R{\T{R}}_\ell(x,t)=\R{R}_0(x,t)$ because $\dist(x,A_q)\geq 2\sigma$. Since $x\notin A_q$, it follows from (\ref{tamaño R0 cerca frontera}) that $|\R{R}_0(x,t)|\leq \frac{1}{4}\delta_{q+2}\lambda_{q+1}^{-6\alpha}$, so
\begin{align*}
	|\Rtildei{k+1}(x,t)|&\leq \varphi^k\|\Rtildei{k}\|_0+(1-\varphi^k)\|\R{R}_i^{k+1}\|_0+\|\Rtildei{k+1}-(\varphi^k\Rtildei{k}+(1-\varphi^k)\R{R}_i^{k+1})\|_0\\
	&\leq \varphi^k\frac{1}{4}\delta_{q+2}\lambda_{q+1}^{-6\alpha}+(1-\varphi^k)\frac{1}{4}\delta_{q+2}\lambda_{q+1}^{-6\alpha}+\frac{1}{4}\delta_{q+2}\lambda_{q+1}^{-6\alpha} \\
	&\leq \frac{1}{2}\delta_{q+2}\lambda_{q+1}^{-6\alpha}
\end{align*}
for $x\notin\Omega^k$. We conclude that (\ref{estimate Tvik 1}) holds for the $(k+1)$-th term.


\section{Gluing in time}\label{section gluing time}

In this short section we develop the third step of the proof of Proposition~\ref{prop iteracion integracion convexa}, which was summarized in \cref{overview gluing time}.

By \cite[Proposition 4.4]{Onsager_final} the matrix $S\coloneqq \mathscr{R}[\partial\chi_i(\T{v}_i-\T{v}_{i+1})]$ satisfies the following bounds for any $N\geq 0$:
\begin{align*}
	\|S\|_{N+\alpha}&\lesssim \delta_{q+1}\ell^{-N+\alpha},\\
	\|(\partial_t+\overline{v}_q\cdot\nabla)S\|_{N+\alpha}&\lesssim \delta_{q+1}\delta_q^{1/2}\lambda_q\ell^{-N-\alpha}.
\end{align*}

We define
\[m\coloneqq \lceil\lambda_q^{1/2}\rceil, \qquad r=\lambda_q^{-3/5}\]
and we fix smooth cutoff functions $\theta_j$ such that
\begin{itemize}
	\item $\theta_j\equiv 1$ in a neighborhood of $A_q+B(0,3\sigma+(j-1)r)$,
	\item the support of $\theta_j$ is contained in $A_q+B(0,3\sigma+jr)$
\end{itemize}
for $1\leq j\leq m$. We define
\[\overline{\mathcal R}^{(2)}_q\coloneqq \sum_{i=1}^m\frac{1}{m}\theta_j\hspace{0.5pt} S.\]
Note that $\supp \overline{\mathcal R}^{(2)}_q(\cdot,t)$ is contained in $A_q+B(0,4\sigma)$ because $mr<\sigma$ for $a$ sufficiently large. Since $r^{-1}\lesssim \ell^{-1}$, we have
\begin{equation}
	\label{estimate Rbar2}
	\|\overline{\mathcal R}^{(2)}_q\|_N\lesssim \delta_{q+1}\ell^{-N+\alpha}.
\end{equation}

Let us estimate the material derivative. We compute
\[(\partial_t+\overline{v}_q\cdot\nabla)\overline{\mathcal R}^{(2)}_q=\sum_{i=1}^m\frac{1}{m}\theta_j(\partial_t+\overline{v}_q\cdot\nabla)\hspace{0.5pt}S+\sum_{i=1}^m\frac{1}{m}\overline{v}_q\cdot\nabla\theta_j\hspace{0.5pt}S.\]
Regarding the second term, it follows from (\ref{est vbarra 2}) and $\norm{\T{v}_\ell}_N\lesssim\ell^{-N}$ that the new field also satisfies $\norm{\overline{v}_q}_N\lesssim \ell^{-N}$. Thus,
\[\norm{\overline{v}_q\cdot\nabla\theta_j}_N\lesssim r^{-1}\ell^{-N}.\]
Since the support of the $\nabla\theta_j$ are pairwise disjoint, we have
\begin{align}
	\|(\partial_t+\overline{v}_q\cdot\nabla)\overline{\mathcal R}^{(2)}_q\|_N&\lesssim \delta_{q+1}\delta_q^{1/2}\lambda_q\ell^{-N-\alpha}+\frac{1}{m}r^{-1}\delta_{q+1}\ell^{-N+\alpha}\label{estimate Rbar2 material derivative}\\&\lesssim\delta_{q+1}\delta_q^{1/2}\lambda_q\ell^{-N-\alpha}.\nonumber
\end{align}

We conclude that the matrix $\overline{\mathcal R}^{(2)}_q$ satisfies the right estimates but we have changed the equation:
\[\Div\overline{\mathcal R}^{(2)}_q=\Div S+\sum_{j=1}^m\frac{1}{m}\nabla\theta_j\cdot S,\]
where we have used that the support of $\Div S$ is contained in $A_q+B(0,3\sigma)$ by (\ref{vtildeik igual inicial}). Let us correct this. We define $\rho_j\coloneqq \nabla\theta_j\cdot S$, whose support is contained in
\[\{x\in \RR^3:3\sigma+(i-1)r<\dist(x,A_q)<3\sigma+ir\}.\]
Therefore, by \cref{frecuencias bajas de dominio pequeño} we have
\[\norm{\rho_j}_{B_{\infty,\infty}^{-1+\alpha}}\lesssim r^{1-\alpha}\norm{\nabla\theta_j\cdot S}_0\lesssim \delta_{q+1}.\]
We wish to apply \cref{invertir divergencia matrices}, so we have to check the compatibility conditions. We fix a Killing field $\xi$ and we compute
\[\int \xi\cdot\rho_j=\int\xi\cdot\Div(\theta_j S)-\int \xi\cdot\theta_j\Div S=-\int\xi\cdot\Div S=-\partial_t\chi_i\int\xi\cdot(\T{v}_i-\T{v}_{i+1}),\]
where we have used (\ref{green id}) and the fact that $\theta_j=1$ on the support of $\Div S$. To show that this integral vanishes, we first compute
\[\frac{d}{dr}\int \xi\cdot(\T{v}_i-\T{v}_\ell)=\int\xi\cdot\Div\left(\T{v}_\ell\otimes\T{v}_\ell-\T{v}_i\otimes\T{v}_i+\T{p}_\ell\Id-\T{p}_i\Id+\R{\T{R}}_i-\R{\T{R}}_\ell\right)=0\]
because of (\ref{green id}) and the fact that the matrix in parentheses is compactly supported due to (\ref{tvell igual inicial}) and (\ref{vi igual inicial}). Since $\T{v}_i(\cdot,t_i)=\T{v}_\ell(\cdot,t_i)$, we see that $\int \xi\cdot(\T{v}_i-\T{v}_\ell)=0$. Repeating this for $\T{v}_{i+1}$ and subtracting, we conclude
\[\int \xi\cdot \rho_j=-\partial_t\chi_j\int\xi\cdot(\T{v}_i-\T{v}_{i+1})=0\]
for any Killing field $\xi$. Therefore, by \cref{invertir divergencia matrices} there exists a smooth symmetric matrix $M_j$ such that $\Div M_j=\rho_j$ and whose support is contained in
\[\{x\in \RR^3:3\sigma+(j-1)r<\dist(x,A_q)<3\sigma+jr\}.\]
Furthermore, we have the estimate
\[\norm{M_j}_0\lesssim r^{-\alpha}\delta_{q+1}.\] 
We define
\[\Div\overline{\mathcal R}^{(3)}_q\coloneqq -\frac{1}{m}\sum_{j=1}^mM_j.\]
By construction of the $M_j$ we have
\[\Div\left(\overline{\mathcal R}^{(2)}_q+\overline{\mathcal R}^{(3)}_q\right)=\partial_t\chi_j(\T{v}_i-\T{v}_{i+1}),\]
as we wanted. Since the supports of the $M_j$ are pairwise disjoint, we see that
\[\|\overline{\mathcal R}^{(3)}_q\|_0=\frac{1}{m}\max_j\norm{M_j}_0\lesssim \frac{1}{m}r^{-\alpha}\delta_{q+1}.\]
It follows from our assumption $b-1<1/10$ that
\[\frac{\delta_{q+1}}{\delta_{q+2}}=\lambda_q^{\beta b(b-1)}<\lambda_q^{1/10}.\]
Therefore, for $a$ is sufficiently large and $\alpha$ sufficiently small the matrix $\Rbar{1}$ defined in (\ref{def Rbar1}) satisfies
\[\|\Rbar{1}\|_0\leq \frac{3}{4}\delta_{q+2}\lambda_{q+1}^{-6\alpha},\]
where we have used that
\[\|\overline{\mathcal R}^{(3)}_q\|_0\leq \frac{1}{2}\delta_{q+2}\lambda_{q+1}^{-6\alpha}\]
due to (\ref{estimate Tvik 1}). Thus, $\Rbar{1}$ is so small that it may be ignored until the next iteration.

Finally, let us conclude the estimates for $\Rbar{2}$, defined in (\ref{def Rbar2}). It follows from \cite[Proposition 4.4]{Onsager_final} that
\begin{align*}
	\|\chi_i(1-\chi_i)(\T{v}_i-\T{v}_{i+1})\R{\otimes}(\T{v}_i-\T{v}_{i+1})\|_{N+\alpha}&\lesssim\delta_{q+1}\ell^{-N+\alpha},\\
	\|(\partial_t+\overline{v}_q\cdot\nabla)[\chi_i(1-\chi_i)(\T{v}_i-\T{v}_{i+1})\R{\otimes}(\T{v}_i-\T{v}_{i+1})]\|_{N+\alpha}&\lesssim\delta_{q+1}\delta_q^{1/2}\lambda_q\ell^{-N+\alpha}.
\end{align*}
Combining this with (\ref{estimate Rbar2}) and (\ref{estimate Rbar2 material derivative}), we infer (\ref{estimate Rbar2 final}) and (\ref{estimate Rbar2 mat der final}).

\section{The perturbation step}\label{section perturbation}

In this section we complete the final step in the proof of Proposition~\ref{prop iteracion integracion convexa}, which is done in Subsections~\ref{SS.squi} to~\ref{SS.ener}. The last part, Subsection~\ref{SS.prolem}, contains the proof of Lemma~\ref{lema integración convexa una iteración}.

For simplicity, we will assume that $\Omega$ is connected. If it had several connected components $\Omega^j$ and we wanted to fix an energy profile $e^j$ in each of them, we would simply carry out the construction of this section in each connected component, taking into account \cref{varias componentes}. Note that $\phi_q$ can be split into cutoffs $\phi_q^j$ associated to each $\Omega^j$.

\subsection{Squiggling stripes and the stress $\T{R}_{q,i}$}\label{SS.squi}
Before we can define the perturbation, we need to introduce several objects. By \cite[Subsection 5.2]{Onsager_final} there exist nonnegative cutoff functions $\eta_i$ with the following properties:
\begin{enumerate}[(i)]
	\item $\eta_i\in C^\infty(\TT^3\times[0,T], [0,1])$.
	\item $\supp\eta_i\cap\supp\eta_j=\varnothing$.
	\item $\eta_i(x,t)=1$ for any $x$ and $t\in I_i$.
	\item $\supp\eta_i\subset \TT^3\times(t_i-\frac{1}{3}\tau_q,t_{i+1}+\frac{1}{3}\tau_q)\cap[0,T]$.
	\item There exists a positive geometric constant $c_0>0$ such that for any $t\in[0,T]$
	\[\sum_i \int_{\TT^3}\eta_i(x,t)^2\,dx\geq c_0.\]
	\item For any $k,N\geq 0$ there exists constants depending on $k,N$ such that
	\[\norm{\partial_t^k\eta_i}_N\lesssim \tau_q^{-N}.\]
\end{enumerate}

We replace $\eta_i(x,t)$ by $\eta_i(mx,t)$ for sufficiently large $m\in \NN$, so that we may assume that the cutoffs $\eta_i$ satisfy
\begin{equation}
	\label{cotas integral cutoffs}
	\widetilde{c}_0\leq \sum_i \int_{\TT^3}\phi_q(x)^2\eta_i(x,t)^2\,dx\leq 2\abs{\Omega}
\end{equation}
for some constant $\widetilde{c}_0$ depending on $\Omega$. This can be done because $A_0$ will contain one of the cubes of a grid of sidelength $m^{-1}$ for sufficiently large $m$. This settles the first inequality. Regarding the second inequality, note that at most 2 of the cutoffs are nonzero at any given time. In addition, the new cutoffs will still satisfy $(\text{i})-(\text{vi})$ but the constants that appear in $(\text{vi})$ will now depend on $\Omega$, too. 

We proceed analogously to \cite{Onsager_final}, defining the amplitudes
\begin{align}
	\rho_q(t)&\coloneqq \frac{1}{3}\left(e(t)-\frac{1}{2}\delta_{q+2}-\int_\Omega\abs{\overline{v}_q}^2dx\right),\\
	\rho_{q,i}(x,t)&\coloneqq \frac{\eta_i(x,t)^2}{\sum_i \int_{\TT^3}\phi_q(x)^2\hspace{0.5pt}\eta_i(x,t)^2\hspace{0.5pt}dx}\hspace{1pt}\rho_q(t). \label{def rho_qi}
\end{align} Note that our definition of $\rho_{q,i}$ differs from the one in \cite{Onsager_final} in the normalization. Next, we define the backwards flows $\Phi_i$ for the velocity field $\overline{v}_q$ as the solution of the transport equation
\[\begin{cases}
	(\partial_t +\overline{v}_q\cdot\nabla)\Phi_i=0,\\
	\Phi_i(x,t_i)=x.
\end{cases}\]
Finally, we define
\begin{align}
	R_{q,i}&\coloneqq \rho_{q,i}\Id-\eta_i^2\Rbar{2},\\
	\T{R}_{q,i}&\coloneqq \frac{\nabla\Phi_i R_{q,i}(\nabla\Phi_i)^t}{\rho_{q,i}}.
\end{align}
It follows from properties $(\text{i})-(\text{iv})$ of $\eta_i$ that
\begin{itemize}
	\item $\supp R_{q,i}\subset \supp \eta_i$,
	\item we have $\sum_i\eta_i^2=1$ on $\supp \Rbar{2}$,
	\item $\supp\T{R}_{q,i}\subset \TT^3\times(t_i-\frac{1}{3}\tau_q,t_{i+1}+\frac{1}{3}\tau_q)$,
	\item $\supp \T{R}_{q,i}\cap\supp \T{R}_{q,j}=\varnothing$ for all $i\neq j$.
\end{itemize}
We collect some other useful estimates:
\begin{lemma}
	\label{propiedades varias}
	For $a\gg 1$ sufficiently large we have
	\begin{equation}
		\label{misc 1}
		\norm{\nabla \Phi_i-\Id}_0\leq \frac{1}{2}\quad \text{for } t\in\supp(\eta_i).
	\end{equation}
	Furthermore, for any $N\geq 0$
	\begin{align}
		\label{misc 2}
		\frac{\delta_{q+1}}{8\lambda_q^\alpha}\leq \abs{\rho_q(t)}&\leq \delta_{q+1} \quad \forall t, \\
		\label{misc 3}
		\norm{\rho_{q,i}}_0&\leq\frac{\delta_{q+1}}{\T{c}_0},\\
		\label{misc 4}
		\norm{\rho_{q,i}}_N&\lesssim \delta_{q+1},\\
		\label{misc 5}
		\norm{\partial_t\rho_{q}}_0&\lesssim \delta_{q+1}\delta_q^{1/2}\lambda_q,\\
		\label{misc 6}
		\norm{\partial_t\rho_{q,i}}_N&\lesssim \tau_q^{-1}\delta_{q+1},\\
		\label{misc 7}
		\|\T{R}_{q,i}\|_N&\lesssim \ell^{-N},\\
		\label{misc 8}
		\|D_{t,q}\T{R}_{q,i}\|_N&\lesssim\tau_q^{-1}\ell^{-N},
	\end{align}
	where $D_{t,q}\coloneqq \partial_t+\overline{v}_q\cdot\nabla$. Moreover, for all $(x,t)$ we have $\T{R}_{q,i}(x,t)\in B\left(\Id,\frac{1}{2}\right)\subset \mathcal{S}^3$.
\end{lemma}
\begin{proof}
	Since our subsolution $(\overline{v}_q,\overline{p}_q,\R{\olsi{R}}_q)$ satisfies analogous estimates to the ones in \cite{Onsager_final}, we may use the same argument to infer (\ref{misc 1}), (\ref{misc 2}) and (\ref{misc 5}). Estimate (\ref{misc 3}) then follows from (\ref{misc 2}) and (\ref{cotas integral cutoffs}). In fact, since our $\rho_{q,i}$ only differs from the one in \cite{Onsager_final} in a time-dependent normalization coefficient that is bounded above and below, the bounds for $\norm{\rho_{q,i}}_N$ are the same except for the implicit constant. Next, it follows from the property $(\text{vi})$ of $\eta_i$ that
	\[\abs{\frac{d}{dt}\sum_i\int_\Omega\phi_q(y,t)^2\eta_i(y,t)^2dy}\lesssim \tau_q^{-1}.\]
	Using this estimate and arguing as in \cite{Onsager_final} we obtain (\ref{misc 6}). Finally, taking into account these estimates for $\rho_{q,i}$ and the bounds (\ref{estimate Rbar2 final}) and (\ref{estimate Rbar2 mat der final}), the facts regarding $\T{R}_{q,i}$ follow as in \cite{Onsager_final}.
\end{proof}

\subsection{Definition of the perturbation} 
The building blocks of the perturbation are Mikado flows, introduced in \cite[Lemma 2.3]{Daneri}:
\begin{lemma}
	\label{mikado}
	For any compact subset of positive-definite matrices $\mathcal{N}\subset \mathcal{S}^3$ there exists a smooth vector field
	\[W:\mathcal{N}\times\TT^3\to\RR^3\]
	such that, for every $R\in \mathcal{N}$
	\begin{equation}
		\label{eq mikado}
		\begin{cases}
			\Div_\xi[W(R,\xi)\otimes W(R,\xi)]=0, \\
			\Div_\xi W(R,\xi)=0
		\end{cases}
	\end{equation}
	and 
	\begin{align*}
		\fint_{\TT^3}W(R,\xi)\,d\xi&=0,\\
		\fint_{\TT^3}W(R,\xi)\otimes W(R,\xi)\,d\xi&=R.
	\end{align*}
\end{lemma}

Since $W(R,\cdot)$ is $\TT^3$-periodic and has zero mean, we may write
\begin{equation}
	\label{fourier W}
	W(R,\xi)=\sum_{k\in \ZZ^3\backslash \{0\}}a_k(R)e^{ik\cdot\xi}
\end{equation}
for some $a_k\in C^\infty(\mathcal{N},\RR^3)$. Similarly, for some $C_k\in C^\infty(\mathcal{N},\mathcal{S}^3)$ we have
\begin{equation}
    \label{fourier WxW}
    W(R,\xi)\otimes W(R,\xi)=R+\sum_{k\neq 0}C_k(R)e^{ik\cdot\xi}.
\end{equation}
It follows from (\ref{eq mikado}) that
\begin{equation}
    \label{ortogonalidad fourier}
    a_k(R)\cdot k=0, \qquad C_k(R)k=0.
\end{equation}
In addition, because of the smoothness we have
\begin{equation}
	\label{est coeff fourier mikado}
	\sup_{R\in\mathcal{N}}\abs{D^N_Ra_k(R)}+\sup_{R\in\mathcal{N}}\abs{D^N_RC_k(R)}\leq \frac{C(\mathcal{N},N,m)}{\abs{k}^m}.
\end{equation}
In the construction this estimates are used with a particular choice of $\mathcal{N}$ (namely $B(\Id,1/2)\subset\mathcal{S}^3$) and $m$. This choice determines the constant $M$ appearing in \cref{prop iteracion integracion convexa}.

With these building blocks, we define the main perturbation term $w_0$ as
\[w_0\coloneqq \sum_i\phi_q(x,t)(\rho_{q,i}(x,t))^{1/2}(\nabla\Phi_i)^{-1} W(\T{R}_{q,i},\lambda_{q+1}\Phi_i)=\sum_i w_{0,i}.\]
Recall that in~\cref{propiedades varias} we saw that $\T{R}_{q,i}$ takes values in the compact subset of positive-definite matrices $\mathcal{N}\coloneqq B(\Id,1/2)\subset \mathcal{S}^3$. Therefore, the previous expression is well-defined. To shorten the notation, we set
\[b_{i,k}(x,t)\coloneqq \phi_q(x,t)(\rho_{q,i}(x,t))^{1/2}a_k(\T{R}_{q,i}(x,t)).\]
Thus, using (\ref{fourier W}) we may write
\[w_0=\sum_{i,k\neq 0}(\nabla\Phi_i)^{-1}b_{i,k}e^{i\lambda_{q+1} k\cdot\Phi_i}.\]

Although Mikado flows are divergence-free, the perturbation $w_0$ will not be solenoidal, in general, due to the other factor. We must add a small correction term $w_c$ so that $w_{q+1}\coloneqq w_0+w_c$ is divergence-free. We set
\[w_{q+1}=\frac{-1}{\lambda_{q+1}}\curl\left(\sum_{i,k\neq 0}(\nabla\Phi)^t\left(\frac{ik\times b_{k,i}}{\abs{k}^2}\right)e^{i\lambda_{q+1}k\cdot\Phi_i}\right).\]
It is clear that $w_{q+1}$ is divergence-free and the correction $w_{q+1}-w_0$ can be seen to be a lower-order term in $\lambda_{q+1}$.

Unlike in \cite{Onsager_final}, $w_{q+1}$ is not the final correction. We must add another small perturbation $w_L$ to ensure the final correction $\T{w}_{q+1}=w_{q+1}+w_L$ has no angular momentum. We define $L\in C^\infty([0,T],\RR^3)$ as
\[L(t)\coloneqq\frac{1}{\lambda_{q+1}}\sum_{i,k\neq 0}\int(\nabla\Phi_i)^t\left(\frac{ik\times b_{k,i}}{\abs{k}^2}\right)e^{i\lambda_{q+1}k\cdot\Phi_i}\,dx,\]
where the integration is undestood to be component-wise. We fix a ball $B\subset A_0$ and we choose $\psi\in C^\infty_c(B,\RR)$ such that $\int\psi=1$. We define the correction $w_L$ as
\[w_L\coloneqq \curl(\psi L).\]
By construction, the correction $\T{w}_{q+1}=w_{q+1}+w_L$ is given by $\T{w}_{q+1}=\curl z$ for a potential vector $z$ such that $\int z\,dx=0$ component-wise. Using the vector identity $\Div(\xi\times z)=z\cdot\curl \xi-\xi\cdot\curl z$, we have
\[\int \xi \cdot \T{w}_{q+1}=\int \xi\cdot \curl z=-\int z\cdot \curl \xi.\]
Since the curl of any Killing field is constant, we conclude
\begin{equation}
	\label{momento angular Twq+1}
	\int \xi\cdot \T{w}_{q+1}\,dx=0\qquad \forall t\in[0,T], \;\forall \xi\in \ker\gradsim.
\end{equation}

\subsection{Estimates on the perturbation} \label{subsect estimates perturbation}
Aside from $w_L$, our perturbation differs from the one in \cite{Onsager_final} in the presence of the factor $\phi_q$. Nevertheless, the factor $\phi_q$ appears multiplying $a_k(\T{R}_{q,i})$. Thus, if we show that $\phi_q\hspace{0.5pt} a_k(\T{R}_{q,i})$ satisfies the same bounds as $a_k(\T{R}_{q,i})$, the same estimates that are derived in~\cite{Onsager_final} will apply here. It follows from (\ref{est coeff fourier mikado}), (\ref{misc 7}) and (\ref{misc 8}) that
\[\|a_k(\T{R}_{q,i})\|_N\lesssim \frac{\ell^{-N}}{\abs{k}^6}, \qquad \|D_{t,q}a_k(\T{R}_{q,i})\|_N\lesssim \frac{\tau_q^{-1}\ell^{-N}}{\abs{k}^6}.\]
Since $\norm{\phi_q}_N\lesssim \lambda_q^{-N/10} \lesssim \tau_q^{-N}\lesssim \ell^{-N}$ we also have
\[\|\phi_q\hspace{0.5pt}a_k(\T{R}_{q,i})\|_N\lesssim \frac{\ell^{-N}}{\abs{k}^6}, \qquad \|D_{t,q}[\phi_q\hspace{0.5pt}a_k(\T{R}_{q,i})]\|_N\lesssim \frac{\tau_q^{-1}\ell^{-N}}{\abs{k}^6}.\]
We conclude that all of the estimates in~\cite{Onsager_final} are also valid here. In particular, we have
\begin{lemma}
	Assuming $a$ is sufficiently large, the perturbations $w_0$, $w_c$ and $w_q$ satisfy the following estimates:
	\begin{align*}
		\norm{w_0}_0+\frac{1}{\lambda_{q+1}}\norm{w_0}_1&\leq \frac{M}{4}\delta_{q+1}^{1/2}, \\
		\norm{w_c}_0+\frac{1}{\lambda_{q+1}}\norm{w_c}_1&\lesssim \delta_{q+1}^{1/2}\ell^{-1}\lambda_{q+1}^{-1},\\
		\norm{w_{q+1}}_0+\frac{1}{\lambda_{q+1}}\norm{w_{q+1}}_1&\leq \frac{M}{2}\delta_{q+1}^{1/2},
	\end{align*}
	where the constant $M$ depends solely on the constant $\T{c}_0$ in (\ref{cotas integral cutoffs}).
\end{lemma}

We carry out the estimates for $w_L$ with more detail. First, due to (\ref{est coeff fourier mikado}) we have $\norm{b_{i,k}}_0\lesssim \norm{\rho_{i,q}}_0^{1/2}\abs{k}^{-6}\lesssim \delta_{q+1}^{1/2}\abs{k}^{-6}$. Next, it follows from (\ref{misc 1}) that $\norm{\nabla\Phi}_0\lesssim 1$. Introducing these estimates in the definition of $L$ we obtain
\[\abs{L(t)}\lesssim \sum_{k\neq 0} \frac{1}{\lambda_{q+1}}\norm{\nabla\Phi_i}_0\norm{b_{i,k}}_0\lesssim \sum_{k\neq 0}\frac{\delta_{q+1}^{1/2}}{\abs{k}^6\lambda_{q+1}}\lesssim \delta_{q+1}^{1/2}\lambda_{q+1}^{-1},\]
where we have used that at most two of the $b_{i,q}$ are nonzero at any given time. Since $\psi$ is fixed throughout the iterative process, we conclude
\begin{equation}
	\label{estimate wL}
	\norm{w_L}_N\lesssim \delta_{q+1}^{1/2}\lambda_{q+1}^{-1}.
\end{equation}
Therefore, the correction $w_L$ is really small. We see that for sufficiently large $a$ the perturbation $\T{w}_{q+1}$ then satisfies
\[
	\norm{w_{q+1}}_0+\frac{1}{\lambda_{q+1}}\norm{w_{q+1}}_1\leq \frac{3}{4}M\delta_{q+1}^{1/2}.
\]

Regarding $\partial_t w_L$, we must first estimate
\begin{align*}
	\|D_{t,q}b_{i,k}\|_0&\lesssim\|\partial \rho_{q,i}+\overline{v}_q\cdot\nabla\rho_{q,i}\|_0\,\|\phi_q\hspace{0.5pt}a_k(\T{R}_{q,i})\|_0+\norm{\rho_{q,i}}_0\,\|D_{t,q}[\phi_q\hspace{0.5pt}a_k(\T{R}_{q,i})]\|_0\\
	&\lesssim \tau_q^{-1}\delta_{q+1}^{1/2}\abs{k}^{-6}.
\end{align*}
Next, we compute
\begin{align*}
	L'(t)&=\frac{1}{\lambda_{q+1}}\sum_{i,k\neq 0}\int D_{t,q}\left[(\nabla\Phi)^t\left(\frac{ik\times b_{k,i}}{\abs{k}^2}\right)e^{i\lambda_{q+1}k\cdot\Phi_i}\right]dx\\
	&=\frac{1}{\lambda_{q+1}}\sum_{i,k\neq 0}\int\Biggl[-\nabla \overline{v}_q(\nabla\Phi)^t\left(\frac{ik\times b_{k,i}}{\abs{k}^2}\right)\\ &\hspace{100pt}+(\nabla\Phi)^t\left(\frac{ik\times D_{t,q}b_{k,i}}{\abs{k}^2}\right)\Biggr]e^{i\lambda_{q+1}k\cdot\Phi_i}\,dx,
\end{align*}
where we have computed the material derivative of $\nabla\Phi_i$ taking into account that $\Phi_i$ solves the transport equation. Since $\norm{\nabla\overline{v}_q}_0\lesssim \delta_q^{1/2}\lambda_q$ by (\ref{est vbarra 3}), we have
\[\abs{L'(t)}\lesssim \lambda_{q+1}^{-1}\delta_{q+1}^{1/2}(\delta_q^{1/2}\lambda_q+\tau_q^{-1})\lesssim \frac{\delta_{q+1}^{1/2}\delta_q^{1/2}\lambda_q}{\lambda_{q+1}^{1-3\alpha}}\]
because $\tau_q^{-1}=\delta_q^{1/2}\lambda_q\ell^{-2\alpha}\lesssim \delta_q^{1/2}\lambda_q\lambda_{q+1}^{-3\alpha}$ by (\ref{tamaño ell}). Since $\psi$ is fixed, we conclude
\begin{equation}
	\label{partialt wL}
	\norm{\partial_t w_L}_N\lesssim \frac{\delta_{q+1}^{1/2}\delta_q^{1/2}\lambda_q}{\lambda_{q+1}^{1-3\alpha}}.
\end{equation}

\subsection{The final Reynolds stress}
Taking into account that $\phi_q$ and $\sum_i\eta_i^2$ equal $1$ on $\supp \Rbar{2}$, it follows from the definition of $R_{q,i}$ that
\[\sum_i\phi_q^2R_{q,i}=-\Rbar{2}+\sum_i\phi_q^2\rho_{q,i}\Id.\]
Using this along with the fact that $\left(\overline{v}_q,\overline{p}_q, \Rbar{1}+\Rbar{2}\right)$ is a subsolution, we obtain:
\begin{align*}
	\partial_t v_{q+1}+\Div(v_{q+1}\otimes v_{q+1})+\nabla\overline{p}_q=&\Div\left(\Rbar{1}+w_L\otimes \overline{v}_q+\overline{v}_q\otimes w_L+w_L\otimes w_L\right)\\
	&+\nabla\left(\sum_i\phi_q^2\rho_{q,i}\right)+(\partial_t w_{q+1}+\overline{v}_q\cdot \nabla w_{q+1})
	\\&+w_{q+1}\cdot \nabla \overline{v}_q+\Div\left(w_{q+1}\otimes w_{q+1}-\sum_i \phi_q^2R_{q,i}\right).
\end{align*} 
Hence, we conclude that we may construct a new subsolution $(v_{q+1},p_{q+1},\R{R}_{q+1})$ by setting
\begin{align*}
	\R{R}_{q+1}&\coloneqq \Rbar{1}+w_L\mathring{\otimes} \overline{v}_q+\overline{v}_q\mathring{\otimes} w_L+w_L\mathring{\otimes} w_L+S-\frac{1}{3}(\tr S)\Id,\\
	p_{q+1}&\coloneqq \overline{p}_q-\sum_i\phi_q^2\rho_{q,i}-\abs{w_L}^2-2\,\overline{v}_q\cdot w_L-\frac{1}{3}\tr S,
\end{align*}
where the smooth symmetric matrix $S$ satisfies
\begin{align*}
	\Div S = \partial_t w_L+\underbrace{(\partial_t w_{q+1}+\overline{v}_q\cdot \nabla w_{q+1})}_\text{transport error}+\underbrace{w_{q+1}\cdot \nabla \overline{v}_q}_\text{Nash error}+\underbrace{\Div\left(w_{q+1}\otimes w_{q+1}-\sum_i \phi_q^2R_{q,i}\right)}_\text{oscillation error} 
\end{align*}
and the support of $S(\cdot,t)$ is contained in $A_q+B(0,5\sigma)$ for all $t\in[0,T]$. If such a matrix existed, the new subsolution $(v_{q+1},p_{q+1},\R{R}_{q+1})$ would equal $(\overline{v}_q,\overline{p}_q,\R{\olsi{R}}_q)$ in $A_{q+1}\times[0,T]$ because the support of $\phi_q$ and $\T{w}_{q+1}(\cdot,t)$ is contained in $A_q+B(0,5\sigma)$. Therefore, it equals $(v_0,p_0,\R{R}_0)$ in $A_{q+1}\times[0,T]$. 

We will show that it is possible to construct $S$ and we will derive the necessary estimates. Let $f\coloneqq \Div S$. Note that the support of $f(\cdot,t)$ is contained in $A_q+B(0,5\sigma)$ for all $t\in[0,T]$ because so is the support of $\phi_q$ and the perturbation. Next, we see that for any Killing field $\xi$ we have
\begin{align*}
	\int \xi\cdot f\,dx&=\frac{d}{dt}\int \xi\cdot \T{w}_{q+1}+\int \xi\cdot\Div\left(\overline{v}_q\otimes w_{q+1}+w_{q+1}\otimes \overline{v}_q+w_{q+1}\otimes w_{q+1}-\sum_i \phi_q^2R_{q,i}\right)\\&=0
\end{align*}
because of (\ref{momento angular Twq+1}) and (\ref{green id}). Therefore, by~\cref{invertir divergencia matrices} there exists a symmetric matrix $S\in C^\infty(\RR^3\times[0,T],\mathcal{S}^3)$ such that $\Div S=f$ and with the stated support. We will now estimate the $C^{\alpha}$-norm of the potential theoretic solution of the equation. Arguing as in~\cref{invertir divergencia matrices}, this yields a bound for the $C^0$-norm of $S$.

We begin by studying the first term in $f$. It follows from (\ref{partialt wL}) and the fact that $\mathscr{R}$ is bounded on Hölder spaces that
\begin{equation}
    \label{inverting partialt wL}
    \norm{\mathscr{R}(\partial_t w_L)}_\alpha\leq \norm{\mathscr{R}(\partial_t w_L)}_{1+\alpha}\lesssim \norm{\partial_t w_L}_\alpha\lesssim  \frac{\delta_{q+1}^{1/2}\delta_q^{1/2}\lambda_q}{\lambda_{q+1}^{1-3\alpha}}.
\end{equation}

The remaining three error terms are analogous to the ones in~\cite{Onsager_final}. Since our fields satisfy the same estimates as in~\cite{Onsager_final}, the estimates for the potential-theoretic solution are completely analogous. We do have to take into account that the oscillation error has a slightly different expression than the one in \cite{Onsager_final}:
\begin{align*}
	&\Div\left(w_{q+1}\otimes w_{q+1}-\sum_i \phi_q^2R_{q,i}\right)=\\ &\hspace{40pt}\Div\left(w_{0}\otimes w_{0}-\sum_i \phi_q^2R_{q,i}\right)+\Div(w_0\otimes w_c+w_c\otimes w_0+w_c\otimes w_c)\equiv \mathcal{O}_1+\mathcal{O}_2.
\end{align*}
The second term $\mathcal{O}_2$ is the same as in \cite{Onsager_final}. Regarding $\mathcal{O}_1$, it follows from the fact that the cutoffs $\eta_i$ have pairwise disjoint support that
\[\mathcal{O}_1=\sum_i\Div(w_{0,i}\otimes w_{0,i}-\phi_q^2R_{q,i}).\]
We use the definition of $w_{0,i}$ and (\ref{fourier WxW}) to write
\begin{align}
    w_{0,i}\otimes w_{0,i} &= \phi_q^2\rho_{q,i}\nabla\Phi_i^{-1}(W\otimes W)(\T{R}_{q,i}, \lambda_{q+1}\Phi_i)\nabla\Phi_i^{-t} \nonumber \\
    &=\phi_q^2\nabla\Phi_i^{-1}\T{R}_{q,i}\nabla\Phi_i^{-t}+\sum_{k\neq 0}\phi_q^2\rho_{q,i}\nabla\Phi_i^{-1}C_k(\T{R}_{q,i})\nabla\Phi_i^{-t}e^{i\lambda_{q+1}k\cdot\Phi_i} \nonumber \\
    &=\phi_q^2 R_{q,i}+\sum_{k\neq 0}\phi_q^2\rho_{q,i}\nabla\Phi_i^{-1}C_k(\T{R}_{q,i})\nabla\Phi_i^{-t}e^{i\lambda_{q+1}k\cdot\Phi_i}. \label{expresion w0xw0}
\end{align}
On the other hand, it follows from (\ref{ortogonalidad fourier}) that
\[\nabla\Phi_i^{-1}C_k\nabla\Phi_i^{-t}\nabla\Phi_i^tk=0,\]
so
\[\mathcal{O}_1=\sum_{i,k\neq 0}\Div(\phi_q^2\rho_{q,i}\nabla\Phi_i^{-1}C_k(\T{R}_{q,i})\nabla\Phi_i^{-t})e^{i\lambda_{q+1}k\cdot\Phi_i},\]
which is the same as in~\cite{Onsager_final} except for the presence of $\phi_q^2$. Nevertheless, it is easy to check, as in~\cref{subsect estimates perturbation}, that $\phi_q^2C_k(\T{R}_{q,i})$ satisfies the same estimates as $C_k(\T{R}_{q,i})$. Hence, we obtain the same bounds for $\mathcal{O}_1$.

Aside from the presence of the cutoff, there is another subtle difference that we have to take into account. The proof in~\cite{Onsager_final} uses that the following inequality holds for a suitable choice of the parameters:
\[\frac{1}{\lambda_{q+1}^{N-\alpha}\ell^{N+\alpha}}\leq \frac{1}{\lambda_{q+1}^{1-\alpha}}.\]
Since our definition of $\ell$ is slightly different, we must check that this inequality holds. Remember that $\lambda_{q+1}\gtrsim\lambda_q^b$ by (\ref{elementary inequalities}). Hence, we have
\[\frac{\lambda_{q+1}^{N-\alpha}\ell^{N+\alpha}}{\lambda_{q+1}^{1-\alpha}}=\lambda_{q+1}^{N-1-\beta(N+\alpha)}\lambda_q^{-(N+\alpha)(1-\beta+3\alpha)}\gtrsim \lambda_q^{[(b-1)(1-\beta)-3\alpha]N-b(1+\beta\alpha)-\alpha(1-\beta+3\alpha)}.\]
Note that $(b-1)(1-\beta)>0$ so by choosing $\alpha$ sufficiently small we can ensure that the coefficient multiplying $N$ is positive. Thus, for sufficiently large $N$ the exponent is positive and choosing $a$ sufficiently large beats any geometrical constant. We conclude that with this choice of parameters the claimed inequality holds.

A similar argument is used several times, for instance, to obtain the inequality $\lambda_{q+1}\ell\geq 1$. Nevertheless, in all of them the difference in the definition of $\ell$ is quite harmless and it only leads to choosing a smaller $\alpha$ and a slightly larger $N$.

In conclusion, the estimates from~\cite{Onsager_final} apply to our case. Combining them with (\ref{inverting partialt wL}) we obtain
\[\norm{\mathscr{R}\,f}_{\alpha}\lesssim \frac{\delta_{q+1}^{1/2}\delta_q^{1/2}\lambda_q}{\lambda_{q+1}^{1-4\alpha}}.\]
Since $\int \xi\cdot f\,dx=0$ for all $t\in[0,T]$ and any Killing field $\xi$, we may use the construction of~\cref{invertir divergencia matrices} to modify $\mathscr{R}\,f$ into a smooth symmetric matrix $S$ such that
\begin{itemize}
    \item $\Div S=f$,
    \item $\supp S(\cdot,t)\subset A_q+B(0,5\sigma)$ for all $t\in[0,T]$,
    \item $\norm{S}_0\lesssim \delta_{q+1}^{1/2}\delta_q^{1/2}\lambda_q\lambda_{q+1}^{-(1-4\alpha)}$.
\end{itemize}
Combining this with (\ref{estimate wL}) we have
\begin{equation}
	\label{resto R}
	\norm{w_L\mathring{\otimes} \overline{v}_q+\overline{v}_q\mathring{\otimes} w_L+w_L\mathring{\otimes} w_L+S-\frac{1}{3}(\tr S)\Id}_0\lesssim  \frac{\delta_{q+1}^{1/2}\delta_q^{1/2}\lambda_q}{\lambda_{q+1}^{1-4\alpha}}.
\end{equation}
We claim that with a suitable choice of the parameters
\begin{equation}
	\label{parameters final}
	\frac{\delta_{q+1}^{1/2}\delta_q^{1/2}\lambda_q}{\lambda_{q+1}}\leq \delta_{q+2}\lambda_{q+1}^{-11\alpha}.
\end{equation}
In that case, (\ref{estimate Rbar1 final}) and (\ref{resto R}) would yield (\ref{final iteración 2}) for sufficiently large $a$, as we wanted. To prove the claimed inequality, we compute
\[\frac{\lambda_{q+1}^{1-11\alpha}\delta_{q+2}}{\delta_{q+1}^{1/2}\delta_q^{1/2}\lambda_q}=\lambda_{q+1}^{1+\beta-11\alpha}\lambda_{q+2}^{-2\beta}\lambda_q^{-1+\beta}\gtrsim\lambda_q^{b(1+\beta-11\alpha)-2b^2\beta-1+\beta}.\]
The condition (\ref{condicion b}) ensures that
\[b-1+\beta+b\beta-2b^2\beta>0,\]
so the exponent is positive for sufficiently small $\alpha>0$. Thus, choosing $a$ sufficiently large beats any numerical constant and (\ref{parameters final}) follows.

\subsection{The new energy profile}\label{SS.ener}
By definition
\[\int_{\Omega}\abs{v_{q+1}}^2dx=\int_{\Omega}\abs{\overline{v}_q}^2dx+2\int_{\Omega}\overline{v}_q\cdot\T{w}_{q+1}dx+\int_{\Omega}\abs{\T{w}_{q+1}}^2dx.\]
Note that in the last two terms we can integrate on the whole $\TT^3$ because the perturbation is supported in $\Omega$. Arguing as in~\cite{Onsager_final} yields the estimate
\[\abs{\int_{\TT^3}\left(2\,\overline{v}_q\cdot w_{q+1}+2w_0\cdot w_c+\abs{w_c}^2\right)dx}\lesssim \frac{\delta_q^{1/2}\delta_{q+1}^{1/2}\lambda_q^{1+2\alpha}}{\lambda_{q+1}}.\]
On the other hand, any term containing $w_L$ will be smaller than this bound by (\ref{estimate wL}). The remaining term is
\[\int_{\TT^3}\abs{w_0}^2dx=\sum_i\int_{\TT^3}\phi_q^2 \tr R_{q,i}dx+\int_{\TT^3}\sum_{i,k\neq 0}\phi_q^2\rho_{q,i}\nabla\Phi_i^{-1}\tr C_k(\T{R}_{q,i})\nabla\Phi_i^{-t}e^{i\lambda_{q+1}k\cdot\Phi_i},\]
where we have used (\ref{expresion w0xw0}). The second term can be estimated as in~\cite{Onsager_final} because $\phi_q^2 C_k(\T{R}_{q,i})$ satisfies the same estimates as $C_k(\T{R}_{q,i})$, as argued several times. Regarding the first term:
\[\sum_i\int_{\TT^3}\phi_q^2(x) \tr R_{q,i}(x,t)\,dx=3\sum_i\int_{\TT^3}\phi_q^2(x)\rho_{q,i}(x,t)\,dx=3\rho_q(t)=e(t)-\frac{1}{2}\delta_{q+2}-\int_\Omega\abs{\overline{v}_q}^2.\]
We conclude that
\[\abs{e(t)-\int_\Omega \abs{v_{q+1}}^2dx-\frac{\delta_{q+2}}{2}}\lesssim \frac{\delta_q^{1/2}\delta_{q+1}^{1/2}\lambda_q^{1+2\alpha}}{\lambda_{q+1}} \stackrel{(\ref{parameters final})}{\lesssim}\delta_{q+2}\lambda_{q+1}^{-9\alpha},\]
which yields (\ref{final iteración 3}) for sufficiently large $a$.

\subsection{Proof of \cref{lema integración convexa una iteración}}\label{SS.prolem}
This lemma is just a simplified version of the construction presented in the previous subsections, so we will just sketch its proof. Since the initial Reynolds stress $\mathring{R}_0$ and its derivatives vanish at $\partial \Omega\times[0,T]$, for any $k\in\NN$ there exist a constant $C_k$ such that for any $x\in \Omega$ we have
\[\vert\mathring{R}_0(x,t)\vert\leq C_k\dist(x,\partial \Omega)^k.\]
Therefore, if $\dist(x,\partial\Omega)<3\lambda^{-1/12}$ we have
\[\vert\mathring{R}_0(x,t)\vert\leq C_{12}(3\lambda^{-1/12})^{12}\lesssim \frac{1}{\lambda^{1/2}}.\]
We fix a smooth cutoff function $\phi\in C^\infty_c(\Omega,[0,1])$ such that
\[\phi(x)=\begin{cases}
	1 \qquad \text{if }\dist(x,\partial \Omega)\geq 3\lambda^{-1/12}, \\ 0 \qquad \text{if }\dist(x,\partial \Omega)\leq 2\lambda^{-1/12}.
\end{cases}\]
It can be chosen so that $\norm{\phi}_N\lesssim \lambda^{N/12}$. This function will control the support of the perturbation.

Next, we define the backwards flows $\Phi$ for the velocity field $v_0$ as the solution of the transport equation
\[\begin{cases}
	(\partial_t +v_0\cdot\nabla)\Phi_i=0,\\
	\Phi_i(x,t_i)=x
\end{cases}\]
and we define
\[\T{R}\coloneqq \nabla\Phi\left(\Id-(2\|\R{R}_0\|_0)^{-1}\R{R}_0\right)\nabla\Phi^{t}.\]
Let $\mathcal{N}\subset \mathcal{S}^3$ be the compact subset of positive definite matrices whose eigenvalues take values between $1/2$ and $3/2$. We see that $\T{R}$ takes values in $\mathcal{N}$. Thus, we may apply \cref{mikado} and define
\[w_0\coloneqq (2\|\R{R}_0\|_0)^{1/2}\phi\,(\nabla \Phi)^{-1}W(\T{R},\lambda \Phi)=\sum_{k\neq 0}(\nabla \Phi)^{-1}b_k e^{i\lambda k\cdot \Phi},\]
with $b_k\coloneqq  (2\|\R{R}_0\|_0)^{1/2}\phi\, a_k(\T{R})$. We also have
\begin{equation}
	\label{una iteracion w0xw0}
	w_0\otimes w_0=2\|\R{R}_0\|_0\,\phi^2\,\T{R}+\sum_{k\neq 0}2\|\R{R}_0\|_0\,\phi^2\nabla\Phi^{-1}C_k(\T{R})\nabla\Phi^{-t}e^{i\lambda k\cdot\Phi}.
\end{equation}
Next, the correction $w_c$ is then defined so that $w\coloneqq w_0+w_c$ is divergence-free:
\begin{equation}
	w_0+w_c=\frac{-1}{\lambda}\curl\left(\sum_{k\neq0}(\nabla\Phi)^t\left(\frac{ik\times b_k}{\abs{k}^2}\right)e^{i\lambda k\cdot\Phi}\right).
	\label{w0+wc una iteración}
\end{equation}
Regarding the angular momentum, we define
\[L(t)\coloneqq\frac{1}{\lambda}\sum_{k\neq 0}\int(\nabla\Phi)^t\left(\frac{ik\times b_{k}}{\abs{k}^2}\right)e^{i\lambda_{q+1}k\cdot\Phi}\,dx.\]
We fix a ball $B\Subset\Omega$ and we choose $\psi\in C^\infty_c(B)$ such that $\int \psi=1$. We add the correction $w_L$ so that the perturbation has vanishing angular momentum:
\[w_L\coloneqq \curl(\psi L).\]
If $\Omega$ has several connected components $\Omega^j$, we will have to consider the partial angular momentum $L^j$. We will need to add one such vortex to each $\Omega^j$ to cancel the angular momentum in each connected component of $\Omega$. We still denote the total correction as $w_L$.

The new velocity field is $v\coloneqq v_0+w_0+w_c+w_L$.  Note that by taking a larger $\lambda$ we may force $B\subset \supp\phi$, so the perturbation vanishes for $\dist(x,\partial\Omega)\leq 2\lambda^{-1/12}$. 

Since $\norm{\phi}_N\lesssim \lambda^{N/12}$, the dominant term is the exponential. Hence, from (\ref{w0+wc una iteración}) and the definition of $L(t)$ we conclude 
\[\norm{v}_N\lesssim \lambda^{N}.\]
Let us denote $D_t\coloneqq\partial_t+v_0\cdot\nabla$. 
Since $D_t\Phi=0$, we see that
$\abs{L'(t)}\lesssim \lambda^{-(1-1/12)}$. We conclude
\[\norm{\partial_t w_L}_N\lesssim \lambda^{-(1-1/12)}.\]

Finally, we define
\begin{align*}
	\R{R}&\coloneqq (1-\phi^2)\R{R}_0+v_0\otimes w_L+w_L\otimes v_0+S-\frac{1}{3}\tr(S)Id, \\
	p&\coloneqq p_0-2\|\R{R}_0\|_0\phi^2 -\abs{w_L}^2-2v_0\cdot w_L-\frac{1}{3}\tr(S),
\end{align*}
where the smooth symmetric matrix $S$ satisfies
\begin{equation}
	\label{ec S una iteración}
	\Div S=\partial_t w_L+\left[D_t w+w\cdot\nabla v_0+\Div\left(w\otimes w-2\|\R{R}_0\|_0\phi^2\T{R}\right)\right]\equiv f.
\end{equation}
It is easy to check that $(v,p,\R{R})$ is a subsolution. Furthermore, the fact that the perturbation has vanishing angular momentum ensures that we may choose $S$ with support contained in $A_*$ by using \cref{invertir divergencia matrices}. Therefore, the $(v,p,\R{R})$ equals the initial subsolution outside $A_*$.

Regarding the estimates, by (\ref{una iteracion w0xw0}) we may write the term in brackets in (\ref{ec S una iteración}) as
\[\sum_{k\neq 0} c_ke^{i\lambda k\cdot\Phi}\]
for certain vectors $c_k$ such that $\norm{c_k}_N\lesssim \abs{k}^{-6}\lambda^{(N+1)/12}$. The standard stationary phase lemma (see \cite{Onsager_final}) yields
\[\norm{\mathscr{R}f}_{1/4}\lesssim \frac{\lambda^{2/12}}{\lambda^{1-1/4}}\leq \lambda^{-1/2}.\]
Continuing the construction of \cref{invertir divergencia matrices}, the claimed bound for $\R{R}$ follows.

Regarding the energy, (\ref{una iteracion w0xw0}) and a standard stationary phase lemma (see \cite{Onsager_final}) yield:
\begin{align*}
	\int_\Omega\abs{v}^2dx&=\int_\Omega \abs{v_0}^2dx+\int_\Omega \abs{w_0}^2dx+O\left(\frac{1}{\lambda^{1-1/12}}\right)\\
	&=\int_\Omega \abs{v_0}^2dx+\int_\Omega 2\|\R{R}_0\|\phi^2\tr(\T{R})\,dx.
\end{align*}
Since $\R{R}_0$ is trace-free, $\tr{\T{R}}=3$. We conclude that (\ref{energía una iteración}) holds for sufficiently large~$\lambda$.

\section{Proof of Theorem~\ref{teorema gordo}}\label{section small}

We are ready to prove our main theorem using Theorems~\ref{extender subs con soporte compacto} and~\ref{teorema integración convexa}. First, we show that the conditions are necessary. Suppose that such a weak solution $(v,p)$ exists. We fix a connected component $\Sigma$ of $\partial \Omega$ and $a\in\RR^3$. In order to give us some room to mollify the subsolution, we also fix a smooth surface $\Sigma'\subset \Omega$ that will be used to approximate $\Sigma$ from the inside of $\Omega$. Given $\varepsilon>0$, if follows from the smoothness of $(v_0,p_0)$ that $\Sigma'$ can be chosen sufficiently close to $\Sigma$ so that 
\begin{align*}
	&\abs{\int_{\Sigma}v_0\cdot \normal-\int_{\Sigma'}v_0\cdot \normal}<\varepsilon,\\[5pt]
	& \abs{\int_{\Sigma}\left[(a\cdot x)\partial_tv_0+(a\cdot v)v_0+p_0 a\right]\cdot \normal-\int_{\Sigma'}\left[(a\cdot x)\partial_tv_0+(a\cdot v)v_0+p_0 a\right]\cdot \normal}<\varepsilon
\end{align*}
for any $t\in[0,T]$.

Next, we fix a mollification kernel $\psi\in C^\infty_c(\RR^3\times\RR)$ whose support is contained in the unit ball and for $0<\ell<\varepsilon$ we define
	\begin{align*}
		v_\ell&\coloneqq v\ast \psi_\ell,\\
		p_\ell&\coloneqq p\ast \psi_\ell+\abs{v_q}^2\ast \psi_\ell-\abs{v_\ell}^2,\\
		\mathring{R}_\ell&\coloneqq v_\ell\mathring{\otimes}v_\ell-(v\mathring{\otimes}v)\ast \psi_{\ell},
	\end{align*}
where $f\mathring{\otimes}g$ denotes the traceless part of the tensor $f\otimes g$. Since $(v,p)$ is a weak solution, it is easy to see that $(v_\ell,p_\ell,\R{R}_\ell)$ is a smooth subsolution in $\RR^3\times(\varepsilon,T-\varepsilon)$. On the other hand, the values of $(v_\ell,p_\ell,\R{R}_\ell)$ on $\Sigma'$ depend only on $(v_0,p_0)$ for $\ell<\dist(\Sigma',\partial\Omega)$ because $(v,p)$ equals $(v_0,p_0)$ on $\overline{\Omega}\times[0,T]$. In particular, we have,
\begin{align*}
	&\lim_{\ell \to 0}(v_\ell,p_\ell,\R{R}_\ell)(x,t)= (v_0,p_0,0)(x,t) \hspace{11pt} \text{uniformly in }(x,t)\in \Sigma'\times[\varepsilon,T-\varepsilon], \\[5pt]
	& \lim_{\ell \to 0}\partial_t v_\ell(x,t)= \partial_t v_0(x,t)\hspace{60pt} \text{uniformly in }(x,t)\in \Sigma'\times[\varepsilon,T-\varepsilon].
\end{align*}
Hence, for sufficiently small $\ell$, for any $t\in[\varepsilon,T-\varepsilon]$ we have
\begin{align*}
	&\abs{\int_{\Sigma'}v_0\cdot \normal-\int_{\Sigma'}v_\ell\cdot \normal}<\varepsilon,\\[5pt]
	& \abs{\int_{\Sigma'}\left[(a\cdot x)\partial_tv_0+(a\cdot v)v_0+p_0 a\right]\cdot \normal-\int_{\Sigma'}\left[(a\cdot x)\partial_tv_\ell+(a\cdot v)v_\ell+p_\ell a -a^t\R{R}_\ell\right]\cdot \normal}<\varepsilon.
\end{align*}
However, these integrals vanish:
\[\int_{\Sigma'}v_\ell\cdot \normal=\int_{\Sigma'}\left[(a\cdot x)\partial_tv_\ell+(a\cdot v)v_\ell+p_\ell a -a^t\R{R}_\ell\right]\cdot \normal=0\qquad \forall t\in(\varepsilon,T-\varepsilon)\]
because $(v_\ell,p_\ell,\R{R}_\ell)$ is a smooth subsolution (see the proof of \cref{pegado subsoluciones}, equations (\ref{integral subs 1}) and (\ref{integral subs 2})). We conclude that for all $t\in(\varepsilon,T-\varepsilon)$ we have
\[\abs{\int_{\Sigma}v_0\cdot \normal}+\abs{\int_{\Sigma}\left[(a\cdot x)\partial_tv_0+(a\cdot v)v_0+p_0 a\right]\cdot \normal}<4\varepsilon.\]
Since $\varepsilon>0$ is arbitrary and $(v_0,p_0)$ is smooth up to the endpoints of the interval, we deduce
\[\int_{\Sigma}v_0\cdot \normal=\int_{\Sigma}\left[(a\cdot x)\partial_tv_0+(a\cdot v)v_0+p_0 a\right]\cdot \normal=0 \qquad \forall t\in[0,T].\]
Taking into account that $a\in \RR^3$ and the connected component $\Sigma$ of $\partial\Omega$ are arbitrary, we see that the conditions in \cref{teorema gordo} are, indeed, necessary.
	
Next, we prove that the conditions are also sufficient. First of all, we may assume that $\Omega'\supset \overline{\Omega}$ is a bounded set with smooth boundary and with a finite number of connected components. Next, by Theorem~\ref{extender subs con soporte compacto} there exists a subsolution $(\T{v}_0,\T{p}_0,\R{\T{R}}_0)\in C^\infty(\RR^3\times[0,T])$ that extends $(v_0,p_0,0)$ outside $\overline{\Omega}$ and whose spatial support is contained in $\Omega'$. In particular, $\supp \R{\T{R}}_0(\cdot,t)$ will be contained in $\overline{\Omega}'\backslash \Omega$. Applying Theorem~\ref{teorema integración convexa}, we obtain a weak solution of the Euler equations $(v,p)$ that equals $(\T{v}_0,\T{p}_0)$ outside $\overline{\Omega}'\backslash \Omega$. Therefore, its support is contained in $\overline{\Omega}'$ and it extends $(v_0,p_0)$. On the other hand, when applying Theorem~\ref{teorema integración convexa} we may prescribe any energy profile $e\in C^\infty([0,T])$ such that
 \[e(t)>\int_{\RR^3}\abs{\T{v}_0(x,t)}^2dx+6\|\R{\T{R}}_0\|_0\abs{\Omega'\backslash \Omega}.\]
Hence, we can define the constant $e_0$ that appears in the statement of \cref{teorema gordo} as any number greater than the right-hand side of the previous inequality. Finally, $v\in C^{\beta}(\RR^3\times[0,T])$, as we wanted, thus completing the proof of the theorem.

    \begin{proof}[Sketch of the proof of \cref{R.teoremagordo}]
    If we construct the spatial extension $(\T{v}_0,\T{p}_0,\R{\T{R}}_0)\in C^\infty(\RR^3\times[0,T])$ so that $\T{v}_0$ has vanishing total angular momentum, we can easily extend in time to a subsolution $(\hat{v}_0,\hat{p}_0,\R{\hat{R}}_0)\in C^\infty(\RR^3\times[0,+\infty)$ whose temporal support is contained in $[0,T')$. By taking the energy profile $e$ larger, if necessary, we may extend it to $[0,T']$ maintaining an analogue of the condition (\ref{condición energía 1}).

    We then carry out the same construction as in \cref{teorema integración convexa} with some minor modifications:
    \begin{itemize}
        \item in \cref{lema integración convexa una iteración} and in the perturbation step of \cref{prop iteracion integracion convexa} we introduce a temporal cutoff so that we do not modify the subsolution at the times when the Reynolds stress is identically $0$ and 
        \item  in the perturbation step of \cref{prop iteracion integracion convexa} we define $\rho_{q,i}\coloneqq \eta_i^2\delta_{q+1}$ instead of (\ref{def rho_qi}) if the interval $(t_i-\frac{1}{3}\tau_q,t_{i+1}+\frac{1}{3}\tau_q)$ is disjoint from $[0,T]$.
    \end{itemize}
    With such an scheme we only prescribe the energy profile in $[0,T]$. However, the energy profile in $[T,T']$ does not differ much (depending on the initial Reynolds stress) from $\int \abs{\hat{v}_0(x,t)}^2dx$. Hence, if we choose $e(0)$ sufficiently large, the final weak solution will be admissible.
    \end{proof}

\section{Open time interval}\label{section intervalo abierto}

In this section we study what happens when the fields are defined in an open interval $(0,T)$ and there is some singular behavior at the endpoints of the interval. So far we have always considered the supremum in time of the spatial Hölder norms of our fields. This is not an option for the situation that we have in mind, which is the setting for our applications.

\subsection{Main result}

Our approach consists in decomposing $(0,T)$ as a countable union of closed intervals $\{\mathcal{I}_k\}_{k=0}^\infty$ meeting only at their endpoints and trying to work in each of them independently. While some parts of the iterative scheme that we have discussed in the previous sections depend only on what is happening at the current time (solving the symmetric divergence equation, for instance), many others do not (whenever we have dealt with transport). Therefore, if the Reynolds stress is nonzero at the endpoints of the intervals, when we try to correct it the subsolution in one interval will affect its neighbor. This propagates the bad estimates from near $t=0$ and $t=T$ to any $\mathcal{I}_k$ after enough iterations of the scheme.

Hence, if we want to isolate the closed intervals and work in each of them independently, we must ensure that the Reynolds stress vanishes identically at their endpoints: 


\begin{lemma}
	\label{anular R en t0}
	Let $(v,p,R)\in C^\infty(\RR^3\times[0,T])$ be a subsolution of the Euler equations. Let $t_0\in (0,T)$ and let $s>0$ be sufficiently small. Suppose that the support of $R(\cdot,t_0)$ is contained in an open set $\Omega$ and that $\Div\Div R(\cdot,t_0)=0$. Then, there exists a smooth subsolution $(\T{v},\T{p},\T{R})$ such that $\T{R}(\cdot,t_0)\equiv 0$ and such that $(\T{v},\T{p},\T{R})=(v,p,R)$ outside $\Omega\times(t_0-s,t_0+s)$. Furthermore, we have the following estimates:
	\begin{align*}
		\norm{\T{v}-v}_N&\leq C(N)\,s\norm{R(\cdot,t_0)}_{C^{N+1}}, \\
		\|\T{R}-R\|_0&\leq C\left(\norm{R(\cdot,t_0)}_{C^0}+s\norm{R(\cdot,t_0)}_{C^1}\norm{v}_0+s^2\norm{R(\cdot,t_0)}^2_{C^1}\right)
	\end{align*}
	for some constants independent of $s$ and $(v,p,R)$.
\end{lemma}
\begin{proof}
	We fix a smooth cutoff function $\chi\in C^\infty_c((-1,1),\RR)$ that equals $1$ in a neighborhood of the origin. Consider the field
	\[\T{v}(x,t)= v(x,t)+w(x,t)\coloneqq v(x,t)-(t-t_0)\,\chi\left(\frac{t-t_0}{s}\right)\Div R(x,t_0).\]
	The condition $\Div\Div R(\cdot,t_0)=0$ ensures that $\T{v}$ is divergence-free. By definition of $\chi$, we see that $w$ vanishes unless $\abs{t-t_0}<s$ so we deduce $\norm{w}_N\leq C\,s\norm{R(\cdot,t_0)}_{C^{N+1}}$. Next, we define
	\begin{align*}
		\T{R}(x,t)\coloneqq &R(x,t)-\chi\left(\frac{t-t_0}{s}\right)R(x,t_0)-\frac{t-t_0}{s}\,\chi'\left(\frac{t-t_0}{s}\right)R(x,t_0)\\
		&+w\otimes v+v\otimes w+w\otimes w.
	\end{align*}
	It is easy to see that $(\T{v},p,\R{R})$ is a subsolution. Furthermore, it follows from our choice of $\chi$ and the fact that $w$ vanishes at $t=t_0$ that $\R{R}(\cdot,t_0)$ is identically $0$. On the other hand, the bound for $\norm{w}_0$ yields the claimed estimate for $\T{R}$.
	
	On the other hand, since the support of $\chi$ is contained in $(-1,1)$ and the support of $R(\cdot,t_0)$ is contained in $\Omega$, we see that the support of $\T{v}-v$ and $\T{R}-R$ is contained in $\Omega\times(t_0-s,t_0+s)$.
	
	Finally, we absorb the trace of $\R{R}$ into the pressure, which preserves the other properties that we have discussed.
\end{proof}

\begin{remark}
	The condition $\Div\Div R(\cdot,t_0)=0$ is quite restrictive, but it can be removed if one is willing to relinquish spatial control of the velocity field. Indeed, we may decompose $\Div R(\cdot,t_0)$ as the sum of a divergence-free field and a gradient, which we absorb into the pressure. The divergence-free part is canceled using the previous lemma. The issue is that the divergence-free component of $\Div R(\cdot,t_0)$ does not have compact support, in general. Thus, we modify the subsolution outside $\supp R$ and we loose the spatial control, which has been our main concern so far. This approach could yield interesting applications in $\TT^3$, nevertheless. However, in $\RR^3$ we would have to modify the construction to address the fact that we have to add perturbations in the whole space. We do not pursue this path here.
	
\end{remark}

Since our construction relies on keeping the velocity fixed at the endpoints of the intervals $\mathcal{I}_k$, we cannot expect to obtain a nonincreasing energy profile for the final weak solution. Indeed, by weak-strong uniqueness it should equal the smooth solution with that initial data (in its domain of definition) and that is not what will will obtain with the convex integration scheme. 

Thus, instead of trying to fix the energy with this construction, we will focus on keeping the changes small after each iteration. Our goal is to ensure that the energy profile can be extended to a continuous function in $[0,T]$. In that case, we may use \cref{teorema integración convexa} to add a (nonsingular) perturbation elsewhere so that the total energy achieves the desired profile.

Hence, the main result that we will prove in this section, which is a nontrivial variation of Theorem~\ref{teorema integración convexa}, is:
\begin{theorem}
	\label{teorema integración convexa segunda parte}
	Let $0<\beta<1/3$. Let $T>0$ and let $\{\mathcal{I}_k\}_{k=0}^\infty$ be a sequence of closed intervals meeting only at their endpoints and such that $(0,T)=\bigcup_{k\geq0}\mathcal{I}_k$. Let $\Omega_k\Subset (0,1)^3$ be a bounded domain with smooth boundary for $k\geq 0$. Let $(v_0,p_0,\R{R}_0)\in C^\infty(\RR^3\times(0,T))$ be a subsolution of the Euler equations such that $\supp\R{R}_0(\cdot,t)\subset \overline{\Omega}_k$ for $t\in\mathcal{I}_k$. In addition, assume that $\Div\Div \R{R}_0$ vanishes. Then, there exists a weak solution of the Euler equations $v \in C^\beta_\text{loc}(\RR^3\times(0,T))$ that equals $v_0$ in $(\RR^3 \backslash \Omega_k)\times\mathcal{I}_k$ for any $k\geq 0$. In addition, $v=v_0$ at the endpoints of the intervals $\mathcal{I}_k$. Furthermore,
	\[\norm{(v-v_0)(\cdot,t)}_{C^0}\leq C\sup_{t\in \mathcal{I}_k}\|\R{R}(\cdot,t)\|_{C^0}^{1/2}\qquad \forall k\geq 0\]
	for some universal constant $C$.
\end{theorem}

\subsection{The iterative process}
As in \cref{prop iteracion integracion convexa}, we are given an initial subsolution $(v_0,p_0,\mathring{R}_0)\in C^\infty(\RR^3\times[0,T])$ and we will iteratively construct a sequence of subsolutions $\{(v_q,p_q,\mathring{R}_q)\}_{q=0}^\infty$ whose limit will be the desired weak solution. To construct the subsolution at step $q$ from the one in step $q-1$, we will add an oscillatory perturbation with frequency $\lambda_q$. Meanwhile, the size of the Reynolds stress will be measured by an amplitude $\delta_q$. These parameters are given by
\begin{align}
	\lambda_q&=2\pi \lceil a^{b^q}\rceil, \label{def lambda segunda parte}\\
	\delta_q&=\lambda_q^{-2\beta}, \label{def delta segunda parte}
\end{align}
The parameters $a,b>1$ are very large and very close to 1, respectively. They will be chosen depending on the exponent $0<\beta<1/3$ that appears in \cref{teorema integración convexa segunda parte}, on $\Omega$ and on the initial subsolution. We introduce another parameter $\alpha>0$ that will be very small. The necessary size of all the parameters will be discovered in the proof.

We will assume that the support of $(v_0,p_0,\R{R}_0)(\cdot,t)$ is contained in $(0,1)^3$. Meanwhile, the support of $\R{R}$ is contained in $\overline{\Omega}\times[0,T]$ for a potentially smaller domain $\Omega$ with smooth boundary. The main difference with respect to \cref{prop iteracion integracion convexa} is that we also assume that $\R{R}_0(\cdot,t)$ vanishes for $t=0$ and $t=T$.

It will be convenient to do an additional rescaling in our problem. In the rescaled problem the initial subsolution will depend on $a$, but we assume that nevertheless there exists a sequence $\{y_N\}_{N=0}^\infty$ independent of the parameters such that
\begin{align}
	\norm{v_0}_N+\norm{\partial_t v_0}_N&\leq y_N, \label{est v no cambia con el rescalado segunda parte}\\
	\norm{p_0}_N&\leq y_N, \\
	\|\R{R}_0\|_N+\|\partial_t \R{R}_0\|_N&\leq y_N. \label{est R no cambia con el rescalado segunda parte}
\end{align}

Since the initial Reynolds stress $\mathring{R}_0$ and its derivatives vanish at $\partial \Omega\times[0,T]$, for any $k\in\NN$ there exist a constant $C_k$ such that for any $x\in \Omega$ we have
\[\vert\mathring{R}_0(x,t)\vert\leq C_k\dist(x,\partial \Omega)^k.
	\label{R inicial segunda parte}\]
Note that the constants $C_k$ are independent of $a$ by (\ref{est R no cambia con el rescalado segunda parte}). We define
\begin{align}
	d_q&\coloneqq\left(\frac{\delta_{q+2}\lambda_{q+1}^{-6\alpha}}{4C_{10}}\right)^{1/10},  \\
	A_q&\coloneqq \{x\in\Omega:\dist(x,\partial\Omega)\geq d_q\}. \label{def Aq segunda parte}
\end{align}
Hence, we have
\[\vert\mathring{R}_0(x,t)\vert\leq \frac{1}{4}\delta_{q+2}\lambda_{q+1}^{-6\alpha} \qquad \text{if }x\notin A_q.\]

On the other hand, since
$\R{R}_0(\cdot,t)$ vanishes for $t=0$ and $t=T$ there exists a constant $C_t$ depending on $\|\partial_t \R{R}_0\|_0$ such that
\[|\R{R}_0(x,t)|\leq C_t\min\{t,T-t\}.\]
Again, the constant $C_t$ does not depend on $a$ because of (\ref{est R no cambia con el rescalado segunda parte}). It depends only on the initial subsolution (before the rescaling). We define 
\begin{align}
	s_q&\coloneqq \frac{\delta_{q+2}\lambda_{q+1}^{-6\alpha}}{4C_t}, \\
	\mathcal{J}_q&\coloneqq [s_q,T-s_q] \label{def Jq}
\end{align}
so that 
\begin{equation}
	|\R{R}_0(x,t)|\leq \frac{1}{4}\delta_{q+2}\lambda_{q+1}^{-6\alpha} \qquad \text{if }t\in [0,T]\backslash \mathcal{J}_q.
	\label{tamaño R fuera intervalo}
\end{equation}

At step $q$ the perturbation will be localized in $A_q\times \mathcal{J}_q$
so that $(v_q,p_q,\mathring{R}_q)$ equals the initial subsolution in $(\RR^3\backslash A_q)\times\left([0,T]\backslash \mathcal{J}_q\right)$. In this region the Reynolds stress is so small that we will ignore it. We will focus on reducing the error in $A_q\times \mathcal{J}_q$.

The complete list of inductive estimates is the following:
\begin{align}
	(v_q,p_q,\mathring{R}_q)&=(v_0,p_0,\mathring{R}_0) \qquad \text{outside } A_q\times \mathcal{J}_q, \label{inductive 1 segunda parte}\\
	\|\mathring{R}_q\|_0&\leq \delta_{q+1}\lambda_q^{-6\alpha}, \label{inductive 2 segunda parte}\\ 
	\norm{v_q}_1&\leq M\delta_q^{1/2}\lambda_q, \label{inductive 3 segunda parte}\\
	\norm{v_q}_0&\leq 1-\delta_q^{1/2}, \label{inductive 4 segunda parte}
\end{align}
where $M$ is a geometric constant that depends on $\Omega$ and is fixed throughout the iterative process. The following instrumental result is key to the proof of Theorem~\ref{teorema integración convexa segunda parte}, and is analogous to Proposition~\ref{prop iteracion integracion convexa}:

\begin{proposition}
	\label{prop iteracion integracion convexa segunda parte}
	There exists a universal constant $M$ with the following property:
	\newline Let $T\geq 1$ and let $\Omega\subset (0,1)^3\subset \RR^3$ be a bounded domain with smooth boundary. Let $(v_0,p_0,\mathring{R}_0)\in C^\infty(\RR^3\times[0,T])$ be a subsolution whose support is contained in $(0,1)^3\times[0,T]$ and such that $\supp\R{R}_0\subset \overline{\Omega}\times[0,T]$. Furthermore, assume that (\ref{est v no cambia con el rescalado segunda parte})$-$(\ref{est R no cambia con el rescalado segunda parte}) are satisfied for some sequence of positive numbers $\{y_N\}_{N=0}^\infty$. Let $0<\beta<1/3$ and
	\begin{equation}
		1<b^2<\min\left\{\frac{1-\beta}{2\beta},\;\frac{11}{10}\right\}.
		\label{condición b segunda parte}
	\end{equation}
	Then, there exists an $\alpha_0$ depending on $\beta$ and $b$ such that for any $0<\alpha<\alpha_0$ there exists an $a_0$ depending on $\beta$, $b$, $\alpha$, $\Omega$ and $\{y_N\}_{N=0}^\infty$ such that for any $a\geq a_0$ the following holds: 
	\newline Given a subsolution $(v_q,p_q,\mathring{R}_q)$ satisfying (\ref{inductive 1 segunda parte})$-$(\ref{inductive 4 segunda parte}), there exists a subsolution $(v_{q+1},p_{q+1},\mathring{R}_{q+1})$ satisfying (\ref{inductive 1 segunda parte})$-$(\ref{inductive 4 segunda parte}) with $q$ replaced by $q+1$. Furthermore, we have
	\begin{equation}
		\label{estimación cambio iteración segunda parte}
		\norm{v_{q+1}-v_q}_0+\frac{1}{\lambda_{q+1}}\norm{v_{q+1}-v_q}_1\leq M\delta_{q+1}^{1/2}.
	\end{equation}
\end{proposition}

As in \cref{teorema integración convexa}, we need an auxiliary lemma to start the iterative process. This is the analogue of Lemma~\ref{lema integración convexa una iteración}:

\begin{lemma}
	\label{lema integración convexa una iteración segunda parte}
	Let $T>0$ and let $\Omega\subset (0,1)^3\subset\RR^3$ be a bounded domain with smooth boundary. Let $(v_0,p_0,\mathring{R}_0)\in C^\infty(\RR^3\times[0,T])$ be a subsolution whose support is contained in $(0,1)^3\times[0,T]$ and such that $\supp\R{R}_0\subset \overline{\Omega}\times[0,T]$. Let $\lambda>0$ be sufficiently large. There exists a subsolution $(v,p,\mathring{R})\in C^\infty(\RR^3\times[0,T])$ that equals the initial subsolution outside the set
	\[\left\{x\in \Omega: \dist(x,\partial \Omega)>\lambda^{-1/12}\right\}\times(\lambda^{-1/3},T-\lambda^{-1/3})\] 
	and such that
	\begin{align*}
		\norm{v}_N&\lesssim\lambda^N  \qquad \forall N\geq 0,\\
		\|\R{R}\|_0&\leq \lambda^{-1/2},
	\end{align*}
	where the implicit constants are independent of $\lambda$. Furthermore, there exist geometric constants $K_1$, $K_2$ such that if $K_1T\norm{v_0}_1\leq 1$, then 
	\[\norm{v-v_0}_0\leq K_2\|\R{R}_0\|^{1/2}_0.\]
\end{lemma}

\subsection{Proof of \cref{prop iteracion integracion convexa segunda parte}}
	The proof is very similar to the proof of \cref{prop iteracion integracion convexa}, but we will only perturb the subsolution at times $t\in \mathcal{J}_{q+1}$. It will be convenient to define
	\begin{align*}
		\gamma&\coloneqq \frac{1}{2}(s_q-s_{q+1}), \\
		\T{\mathcal{J}}_q&\coloneqq [s_q-\gamma,T-s_q+\gamma].
	\end{align*}
	The parameters $\ell$ and $\tau_q$ are defined as in \cref{prop iteracion integracion convexa}. Let us compare the time parameters:
	\[\frac{s_q}{\tau_q}\gtrsim\delta_{q+2}\lambda_{q+1}^{-6\alpha}\delta_q^{1/2}\lambda_q\ell^{-2\alpha}\gtrsim(\lambda_{q+1}^3\ell)^{-2\alpha}\lambda_q^{1-\beta-2b^2\beta}.\]
	Note that the exponent of $\lambda_q$ is greater than $0$ by our assumption on $b$. Since we may assume that $2d_{q+1}\leq d_q$, we conclude that for sufficiently small $\alpha$ and sufficiently large $a$ we have
	\[\gamma\gg \tau_q.\]
	Hence, the temporal cutoffs that we will need to use will not be too sharp.
	 
	After these definitions, we can prove the result in a similar manner to \cref{prop iteracion integracion convexa}. As before, we divide the proof in four steps: 
	\begin{itemize}
		\item[1.] Preparing the subsolution. The beginning of the iterative process is identical to the one in \cref{prop iteracion integracion convexa}: we use a convolution kernel in space $\psi_\ell$ to mollify $(v_q,p_q,\R{R}_q)$ into $(v_\ell,p_\ell,\R{R}_\ell)$ and then we glue it in space to $(v_0,p_0,\R{R}_0)$, obtaining a subsolution $(\T{v}_\ell,\T{p}_\ell,\R{\T{R}}_\ell)$ that equals $(v_0,p_0,\R{R}_0)$ in $B_2$.
		
		We must, however, introduce a minor modification: we add a correction $w_L$ to ensure that $\T{v}_\ell+w_L$ has the same angular momentum as $v_0$. Note that $v_q$ has the same angular momentum as $v_0$ because they are equal at $t=0$ and subsolutions preserve angular momentum, as argued several times. In addition, it is easy to check that mollifying does not change the total angular momentum, so $v_\ell-v_0$ has $0$ total angular momentum. This may not be true for $\T{v}_\ell-v_0$. Nevertheless, since we are gluing in a (small) region where $v_\ell$ is very similar to $v_0$ (in terms of the parameters), the change in the angular momentum will be very small. Therefore, we may add a small correction $w_L$ to $\T{v}_\ell$ in $A_q$ while keeping the desired estimates. Of course, we modify the pressure and the Reynolds stress accordingly to obtain a subsolution, which we still denote as $(\T{v}_\ell,\T{p}_\ell,\R{\T{R}}_\ell)$, for simplicity.
				
		\item[2.] Gluing in space. In the intervals such that $[t_i-\tau_q,t_i+\tau_q]\cap \mathcal{J}_q \neq \varnothing$ the process is exactly the same as in \cref{prop iteracion integracion convexa}. In the rest of the intervals it suffices to take $(\T{v}_i,\T{p}_i,\R{\T{R}}_i)=(v_0,p_0,\R{R}_0)$ because outside $\mathcal{J}_q$ we have $|\R{R}|_0\leq \frac{1}{4}\delta_{q+2}\lambda_{q+1}^{-6\alpha}$, as required by (\ref{estimate vi 1}).
				
		Regarding the estimates (\ref{estimate vi 2})$-$(\ref{estimate zi 2}), remember that $(v_q,p_q,\R{R}_q)$ equals the initial subsolution for $t\notin \mathcal{J}_q$. Therefore, it follows from (\ref{est v no cambia con el rescalado segunda parte})$-$(\ref{est R no cambia con el rescalado segunda parte}) and standard estimates for mollifiers that for $t\notin \mathcal{J}_q$ we have
		\begin{align*}
			\norm{(v_\ell-v_0)(\cdot,t)}_{C^N}+\norm{(\partial_t v_\ell-\partial_t v_0)(\cdot,t)}_{C^N}&\lesssim \ell^2, \\
			\norm{(p_\ell-p_0)(\cdot,t)}_{C^N}&\lesssim \ell^2, \\
			\|(\mathring{R}_\ell-\mathring{R}_0)(\cdot,t)\|_{C^N}&\lesssim \ell^2.
		\end{align*}
		The same bounds hold for $(\T{v}_\ell,\T{p}_\ell,\R{\T{R}}_\ell)$. Therefore, this choice of $(\T{v}_i,\T{p}_i,\R{\T{R}}_i)$ for the rest of the intervals satisfies (\ref{estimate vi 2})$-$(\ref{estimate zi 2}). Although (\ref{vi en ti}) is not satisfied, its only use in the following step is to ensure that $\T{v}_i$ and $\T{v}_\ell$ have the same angular momentum. This is exactly what we did at the end of the previous step.
		
		\item[3.] Gluing in time. The construction is the same as in \cref{prop iteracion integracion convexa} but due to our choice of $(\T{v}_i,\T{p}_i,\R{\T{R}}_i)$ we are actually gluing only in $\T{\mathcal{J}}_q$ (remember that $\tau_q\ll \gamma$). We see that $(\overline{v}_q,\overline{v}_q,\R{\olsi{R}}_q)$ equals $(v_0,p_0,\R{R}_0)$ for $t\notin \T{\mathcal{J}}_q$ and the support of $\Rbar{2}$ is contained in $[A_q+B(0,4\sigma)]\times\T{\mathcal{J}}_q$.
			
		\item[4.] Perturbation. There are only two changes with respect to \cref{prop iteracion integracion convexa}:
		\begin{itemize}
			\item We define the amplitudes $\rho_{q,i}$ as $\rho_{q,i}(x,t)\coloneqq \eta_i(x,t)^2\hspace{1pt}\delta_{q+1}$.
			\item Instead of the cutoff $\phi_q$ we use $\T{\phi}_q(x,t)\coloneqq \phi_q(x)\theta_q(t)$ for some smooth cutoff $\theta_q\in C^\infty_c(I_{q+1})$ that equals $1$ in $\T{\mathcal{J}}_q$. Hence, $\T{\phi}_q=1$ in the support of $\Rbar{2}$.
		\end{itemize}
		The amplitudes $\rho_{q,i}(x,t)$ clearly satisfy the same estimates as in \cref{prop iteracion integracion convexa}. In particular, $\R{R}_{q,i}$ takes values in $B(\Id,1/2)$. On the other hand, $\T{\phi}_q\hspace{0.5pt}a(\R{R}_{q,i})$ will satisfy the same estimates as $a(\R{R}_{q,i})$ because we may choose $\theta_q$ so that $\abs{\partial_t\theta_q}\lesssim \gamma^{-1}\leq \tau_q^{-1}$.
		
		We conclude that we may carry out the same construction as in the perturbation step of \cref{prop iteracion integracion convexa} (except for  fixing the energy). The new subsolution will satisfy (\ref{inductive 2 segunda parte})$-$(\ref{estimación cambio iteración segunda parte}). In addition, the cutoff $\T{\phi}_q$ ensures that the support of the perturbation is contained in $A_{q+1}\times \mathcal{J}_{q+1}$, as required by the inductive hypothesis (\ref{inductive 1 segunda parte}).
		
		Finally, we emphasize an important difference: the constant $M$ in this case is universal. In \cref{prop iteracion integracion convexa} it depended on $\Omega$ because so did the amplitudes $\rho_{q,i}$. Since here they are independent of $\Omega$, arguing as in \cite{Onsager_final} yields a universal $M$.
	\end{itemize}
\subsection{Proof of \cref{lema integración convexa una iteración segunda parte}}

		We fix a cutoff $\theta\in C^\infty_c((0,T))$ such that
		\[\theta(t)=\begin{cases}1 \quad \text{if }t\in (2\lambda^{-1/3},T-2\lambda^{-1/3}), \\ 0 \quad \text{if }t\notin (\lambda^{-1/3},T-\lambda^{-1/3}). \end{cases}\]
	Since $\R{R}_0$ vanishes at $t=0$ and $t=T$, it is clear that $|(1-\theta)\R{R}_0|\lesssim \lambda^{-1/2}$.
	
	We then carry out the same construction as in \cref{lema integración convexa una iteración} but replacing the cutoff $\phi(x)$ by $\T{\phi}(x,t)\coloneqq \phi(x)\theta(t)$. This ensures that the perturbation vanishes if $t\notin (\lambda^{-1/3},T-\lambda^{-1/3})$. It does worsen the estimates, but it is not catastrophic. The most significant change is that when we write the term in brackets in (\ref{ec S una iteración}) as $\sum_{k\neq 0}c_k e^{i\lambda k\cdot\Phi}$, the vectors $c_k$ now satisfy the estimates
	\[\norm{c_k}_N\lesssim \abs{k}^{-6} \lambda^{1/3+N/12}.\]
	Applying the stationary phase lemma now yields
	\[\norm{\mathscr{R}f}_{1/15}\lesssim \frac{\lambda^{1/3+(1+1/15)/12}}{\lambda^{1-1/15}}<\lambda^{-1/2},\]
	which, continuing the construction in \cref{lema integración convexa una iteración}, yields the desired bound for $\R{R}$. Since the exponent is actually smaller than $-1/2$, we can expend the extra factor in beating any geometric constant.
	
	Finally, let us derive a precise estimate for $v-v_0$. Recall that
	\[w_0= 2\|\R{R}_0\|_0\,\T{\phi}\hspace{0.7pt}{}^2\hspace{0.3pt}(\nabla \Phi)^{-1}W(\T{R},\lambda \Phi).\]
	The last term can be bounded by a geometric constant depending on the compact subset $\mathcal{N}$ used when applying \cref{mikado}. On the other hand, by standard estimates for the transport equation there exists a geometric constant $K_1$ such that
	\[\norm{\nabla\Phi-\Id}_0\leq \frac{1}{2}K_1T\norm{v_0}_1.\]
	By our assumption on $T\norm{v_0}_1$, we have $\norm{\nabla\Phi-\Id}_0\leq 1/2$, so $\|(\nabla\Phi)^{-1}\|_0\leq 2$. We conclude that
	\[\norm{w_0}_0\leq \frac{1}{2}K_2\|\R{R}_0\|_0\]
	for some geometric constant $K_2$. Since $w_c+w_L=O(\lambda^{-1})$, it is clear that the required estimate holds for sufficiently large $\lambda$.

\subsection{Proof of \cref{teorema integración convexa segunda parte}}
To simplify the notation of this proof, given a map $f$ defined in $\RR^3\times I$ for some interval $I$, we will write
\[\norm{f}_{N,I}\coloneqq \supp_{t\in I}\norm{f(\cdot,t)}_{C^N}.\]

By dividing in half the intervals $\mathcal{I}_k$ as many times as necessary, we may assume that
\[2K_1\abs{\mathcal{I}_k}\norm{v_0}_{1,\mathcal{I}_k}\leq 1,\]
where $\abs{\mathcal{I}_k}$ is the length of the interval $\mathcal{I}_k$ and $K_1$ is the constant that appears in \cref{lema integración convexa una iteración segunda parte}. 

We apply \cref{anular R en t0} at the endpoints of each interval $\mathcal{I}_k$, obtaining a new subsolution $(\T{v}_0,\T{p}_0,\R{\T{R}}_0)$ in which the Reynolds stress vanishes at the endpoints of all the intervals $\mathcal{I}_k$. In addition, we have $(\T{v}_0,\T{p}_0,\R{\T{R}}_0)=(v_0,p_0,\R{R}_0)$ outside $\Omega_k$ for $t\in\mathcal{I}_k$. By taking $s$ sufficiently small in each application of \cref{anular R en t0}, we may assume that 
\begin{align*}
	\norm{\T{v}_0-v_0}_{0,\mathcal{I}_k}&\leq \|\R{R}_0\|_{0,\mathcal{I}_k}, \\
	\norm{\T{v}_0}_{1,\mathcal{I}_k}&\leq 2\norm{v_0}_{1,\mathcal{I}_k}, \\
	\|\R{\T{R}}_0\|_{0,\mathcal{I}_k}&\leq 2\|\R{R}_0\|_{0,\mathcal{I}_k}.
\end{align*}

Once we have a subsolution in which the Reynolds stress vanishes at the endpoints of the intervals $\mathcal{I}_k$, we may work in each of them independently. Indeed, our constructions keep the subsolution fixed near the endpoints of the intervals.

We fix $k\geq 0$, and we fix a sequence of positive numbers $\{y_{N,k}\}_{N=0}^\infty$ such that
\begin{align}
	\norm{\T{v}_0}_{N,\mathcal{I}_k}+\norm{\partial_t \T{v}_0}_{N,\mathcal{I}_k}&\leq y_N,\\
	\norm{\T{p}_0}_{N,\mathcal{I}_k}&\leq y_N, \\
	\|\R{\T{R}}_0\|_{N,\mathcal{I}_k}+\|\partial_t \R{\T{R}}_0\|_{N,\mathcal{I}_k}&\leq y_N.
\end{align}
We choose $b$ satisfying (\ref{condición b segunda parte}) and we choose $\alpha$ smaller than the threshold given by \cref{prop iteracion integracion convexa segunda parte}. Let $a_{0,k}$ be the threshold given by \cref{prop iteracion integracion convexa segunda parte} when applied to $\Omega_k$ and our sequence $\{y_{N,k}\}_{N=0}^\infty$. For $a_k\geq a_{0,k}$ we consider the parameters $\lambda_{q,k}$ and $\delta_{q,k}$ defined as in (\ref{def lambda segunda parte}) and (\ref{def delta segunda parte}).

Taking $a_k$ larger if necessary, we apply \cref{lema integración convexa una iteración segunda parte} with the parameter $\lambda_{1,k}^{12\alpha}$, obtaining a subsolution $(v_1,p_1,\R{R}_1)$ that equals $(\T{v}_0,\T{p}_0,\R{\T{R}}_0)$ outside 
\[\left\{(x,t)\in \Omega_k\times\mathcal{I}_k: \dist((x,t),\partial (\Omega_k\times\mathcal{I}_k))>\lambda_{1,k}^{-4\alpha}\right\}\]
and satisfying the following estimates:
\[\norm{v_1}_{N,\mathcal{I}_k}\leq C_{N,k} \lambda_{1,k}^{12N\alpha}, \qquad \|\R{R}_1\|_{0,\mathcal{I}_k}\leq \lambda_{1,k}^{-6\alpha},\]
where the constants $C_{N,k}$ are independent of $\lambda_{1,k}$ but they will depend on $\Omega_k$, $\mathcal{I}_k$ and the initial subsolution. Furthermore, we have
\[\norm{v_1-\T{v}_0}_{0,\mathcal{I}_k}\leq K_2\|\R{\T{R}}_0\|_{0,\mathcal{I}_k}\leq 2K_2\|\R{R}_0\|_{0,\mathcal{I}_k},\]
where we have used that $K_1\abs{\mathcal{I}_k}\norm{\T{v}_0}_{1,\mathcal{I}_k}\leq 1$.

Next, we consider the scale invariance of the Euler equations and subsolutions:
\[v(x,t)\mapsto \Gamma \hspace{0.5pt}  v(x,\Gamma t), \quad p(x,t)\mapsto \Gamma^2 p( x,\Gamma t), \quad \R{R}(x,t)\mapsto \Gamma^2\R{R}( x,\Gamma t).\]
We choose $\Gamma = \delta_{2,k}^{1/2}$ and we begin to work in this rescaled setting, which we will indicate with a superscript $r$. Note that $(\T{v}_0^r,\T{p}_0^r,\R{\T{R}}{}^r_0)$ still satisfies (\ref{est v no cambia con el rescalado segunda parte})$-$(\ref{est R no cambia con el rescalado segunda parte}) with the same sequence $\{y_{N,k}\}_{N=0}^\infty$. Regarding $(v_1^r,p_1^r,\R{R}_1^r)$, it follows from the definition of the rescaling that
\[\|\R{R}_1^r\|_{0,\mathcal{I}_k^r}\leq \delta_{2,k}\lambda_{1,k}^{-6\alpha}.\]
On the other hand, since the constants $C_{N,k}$ are independent of $\lambda_{1,k}$, for sufficiently large $a_k$ we have
\begin{align*}
	\norm{v_1^r}_{0,\mathcal{I}_k^r}&=\delta_{2,k}^{1/2}\norm{v_1}_{0,\mathcal{I}_k}\leq \delta_{2,k}^{1/2}C_{0,k}\leq 1-\delta_{1,k}^{1/2}, \\
	\norm{v_1^r}_{1,\mathcal{I}_k^r}&=\delta_{2,k}^{1/2}\norm{v_1}_{1,\mathcal{I}_k}\leq \delta_{2,k}^{1/2}C_{1,k}\lambda_{1,k}^{12\alpha}\leq  M\delta_{1,k}^{1/2}\lambda_{1,k}.
\end{align*}
Finally, $(v_1^r,p_1^r,\R{R}_1^r)=(\T{v}_0^r,\T{p}_0^r,\R{\T{R}}{}^r_0)$ outside
\[\left\{(x,t)\in \Omega_k\times\mathcal{I}_k^r: \dist((x,t),\partial (\Omega_k\times\mathcal{I}_k^r))>\lambda_1^{-4\alpha}\right\}.\]
Let us consider the sets $A_{q,k}$ and $\mathcal{J}_{q,k}$ that we obtain when we apply the definition of (\ref{def Aq segunda parte}) and (\ref{def Jq}) to $\Omega_k$ and $\mathcal{I}_k^r$. We see that $(v_1^r,p_1^r,\R{R}_1^r)=(\T{v}_0^r,\T{p}_0^r,\R{\T{R}}{}^r_0)$ outside $A_{1,k}\times \mathcal{J}_{1,k}$ for sufficiently small $\alpha$ and sufficiently large $a_k$. 

We conclude that $(v_1^r,p_1^r,\R{R}_1^r)$ satisfies the inductive hypotheses (\ref{inductive 1 segunda parte})$-$(\ref{inductive 4 segunda parte}) in the interval $\mathcal{I}_k^r$ with initial subsolution $(\T{v}_0^r,\T{p}_0^r,\R{\T{R}}{}^r_0)$. In addition, for $t\in\mathcal{I}_k^r$ the initial subsolution satisfies (\ref{est v no cambia con el rescalado segunda parte})$-$(\ref{est R no cambia con el rescalado segunda parte}) with the sequence $\{y_{N,k}\}_{N=0}^\infty$ and the support of $\R{\T{R}}{}^r_0(\cdot,t)$ is contained in $\overline{\Omega}_k$. Finally, by taking $a_k$ even larger, we may assume that $\abs{\mathcal{I}_k^r}\geq 1$.

Applying \cref{prop iteracion integracion convexa segunda parte} in each interval $\mathcal{I}_k^r$ and undoing the scaling, we obtain a sequence of subsolutions $\{(v_q,p_q,\R{R}_q)\}_{q=1}^\infty\in C^\infty(\RR^3\times (0,T))$ that equal $(v_0,p_0,0)$ in $(\RR^3\backslash\Omega_k)\times \mathcal{I}_k$ for any $k\geq 0$. In addition, for $t\in \mathcal{I}_k$ we have: 
\begin{align}
	&\|\R{R}_q(\cdot,t)\|_{C^0}\leq \delta_{q+1,k}, \\
	& \norm{(v_{q+1}-v_q)(\cdot,t)}_{C^0}+\frac{1}{\lambda_{q+1,k}}\norm{(v_{q+1}-v_q)(\cdot,t)}_{C^1}\leq M\delta_{1,k}^{-1/2}\delta_{q+1,k}^{1/2}. \label{cambio iteración segunda parte}
\end{align}
We see that $v_q$ converges uniformly in compact subsets of $\RR^3\times(0,T)$ to some continuous map $v$. On the other hand, note that the pressure is the only compactly supported solution of
\[\Delta p_q=\Div\Div(-v_q\otimes v_q+\R{R}_q).\]
Therefore, $p_q$ also converges to some pressure $p\in L^{s}(\RR^3)$ for any $1\leq s<\infty$. Since $\R{R}_q$ converges to $0$ uniformly in compact subsets of $\RR^3\times(0,T)$, we conclude that $(v,p)$ is a weak solution of the Euler equations in $\RR^3\times(0,T)$.

Furthermore, using (\ref{cambio iteración segunda parte}) we obtain
\begin{align*}
	\sum_{q=1}^\infty \norm{v_{q+1}-v_q}_{\beta',\mathcal{I}_k}&\leq \sum_{q=1}^\infty C(\beta',\beta) \norm{v_{q+1}-v_q}_{0,\mathcal{I}_k}^{1-\beta'}\norm{v_{q+1}-v_q}_{1,\mathcal{I}_k}^{\beta'}\\
	&\leq C(\beta',\beta)\sum_{q=1}^\infty (M\delta_{1,k}^{-1/2}\delta_{q,k}^{1/2})^{1-\beta'}(M\delta_{1,k}^{-1/2}\delta_{q,k}^{1/2})^{\beta'} \\
	&\leq M\,C(\beta',\beta)\delta_{1,k}^{-1/2}\sum_{q=1}^\infty\lambda_{q,k}^{\beta'-\beta},
\end{align*}
so $\{v_q\}_{q=1}^\infty$ is uniformly bounded in $C^0_tC^{\beta'}_x$ in any compact subset $I\subset (0,T)$ for all $\beta'<\beta$. Arguing as in \cite{Onsager_final} we obtain (local) time regularity. We conclude that $v\in C^{\beta''}_\text{loc}(\RR^3\times(0,T))$, with $\beta''<\beta'<1/3$ arbitrary.

Finally, we compute the difference between $v_0$ and $v$. We write $b=1+\gamma$, so that $b^q-1\geq \gamma q$. Taking $a_k$ sufficiently large, we have
\begin{align*}
	\norm{v-v_1}_{0,\mathcal{I}_k}&\leq\sum_{q=1}^\infty \norm{v_{q+1}-v_q}_{0,\mathcal{I}_k}\leq \sum_{q=1}^\infty M\delta_{1,k}^{-1/2}\delta_{q+1,k}^{1/2}\lesssim \sum_{q=1}^\infty a_k^{\beta b(1-b^q)}\lesssim \sum_{q=1}^\infty a_k^{-\beta b\gamma q} \\
	&\lesssim a_k^{-\beta b \gamma}\sum_{q=0}(a_k^{\beta b \gamma})^{-q}\lesssim a_k^{-\beta b \gamma}.
\end{align*}
Thus, by further increasing $a_k$, we have $\norm{v-v_1}_{0,\mathcal{I}_k}\leq \|\R{R}_0\|_{0,\mathcal{I}_k}$. We conclude
\[\norm{v-v_0}_{0,\mathcal{I}_k}\leq \norm{v-v_1}_{0,\mathcal{I}_k}+\norm{v_1-\T{v}_0}_{0,\mathcal{I}_k}+\norm{\T{v}_0-v_0}_{0,\mathcal{I}_k}\leq (2+2K_2)\|\R{R}_0\|_{0,\mathcal{I}_k}.\]

\section{Vortex sheet}\label{section vortex sheet}

This section is devoted to the proof of \cref{teorema vortex sheet}. The key ingredient is \cref{teorema integración convexa segunda parte}, which is applied to an initial subsolution with the appropriate behavior. We fix a parameter $0<\lambda<(4T)^{-1}$. We choose an even function $f\in C^\infty_c((-1,1))$ such that $\int f=2$ and we define
\[F(x)\coloneqq \int_{-1}^xf(s)\,ds, \qquad G(x)\coloneqq \int_{-1}^xsf(s)\,ds.\]
Let $v_0$ and $R_0$ be the periodic extension of
\begin{align*}
	v_0(x,t)&\coloneqq \left[1-F\left(\frac{x_3-1/4}{\lambda t}\right)+F\left(\frac{x_3-3/4}{\lambda t}\right)\right]e_1, \\[3pt]
	\R{R}_0(x,t)&\coloneqq \lambda\left[G\left(\frac{x_3-1/4}{\lambda t}\right)-G\left(\frac{x_3-3/4}{\lambda t}\right)\right](e_1\otimes e_3+e_3\otimes e_1).
\end{align*}
Direct computation shows that $\partial_t v_0=\Div \R{R}_0$. It is also clear that $v_0\cdot\nabla v_0=0$. Therefore, the triplet $(v_0,0,\R{R}_0)$ is a subsolution.

Since $f$ is even, the support of $G$ is contained in $(-1,1)$. For $k\geq 0$ let us define $\mathcal{I}_k\coloneqq [2^{-(k+1)}T,2^{-k}T]$ and
\[\Omega_k\coloneqq \TT^2\times\left[\left(\frac{1}{4}- 2^{-k}\lambda T,\;\frac{1}{4}+2^{-k}\lambda T\right)\cup\left(\frac{3}{4}- 2^{-k}\lambda T,\;\frac{3}{4}+2^{-k}\lambda T\right)+\ZZ\right].\]
We see that $(v_0,0,\R{R}_0)$ equals $(u_0,0,0)$ outside $\Omega_k$ for $t\in\mathcal{I}_k$. Since $\Div\Div\R{R}_0=0$, we may apply \cref{teorema integración convexa segunda parte}, obtaining a weak solution of the Euler equations $v\in C^\beta_\text{loc}(\TT^3\times(0,T))$ that equals $(u_0,0,0)$ outside $\Omega_k$ for $t\in\mathcal I_k$ (in fact, for $t\in(0,2^{-k}T)$ because $\Omega_j\subset\Omega_k$ for all $j\leq k$). In particular, the initial datum is $u_0$. Furthermore, there exists a universal constant $C$ such that for any $k\geq0$ we have
\begin{equation}
	\norm{v-v_0}_{0,\mathcal{I}_k}\leq C\|\R{R}\|_{0,\mathcal{I}_k}\leq C'\lambda.
	\label{diferencia vortex}
\end{equation}

Let us estimate the energy. First, we define
\[\delta\coloneqq 2-\int_{-1}^1[1-F(s)]^2ds.\]
Note that $\delta>0$ because the continuous function $F$ takes values in $[0,2]$ and it satisfies $F(0)=1$, since $f$ is even. We compute
\begin{align*}
	\int_{\TT^3}\abs{v_0(x,t)}^2\,dx&=1-4\lambda t+\int_{1/4-\lambda t}^{1/4+\lambda t}\left[1-F\left(\frac{x_3-1/4}{\lambda t}\right)\right]^3dx_3\\
	&+\int_{3/4-\lambda t}^{3/4+\lambda t}\left[-1+F\left(\frac{x_3-3/4}{\lambda t}\right)\right]^3dx_3=1-2\delta \lambda t,
\end{align*}
where we have used that $v_0=u_0$ except close to $x_3=1/4$ and $x_3=3/4$. Next, we write
\[\int_{\TT^3}\abs{v}^2=\int_{\TT^3}\abs{v_0}^2+\int_{\TT^3}\abs{v-v_0}^2+2\int_{\TT^3}v_0\cdot(v-v_0)^2.\]
Let us choose $\lambda$ sufficiently small so that $\lambda< (32C')^{-1}\delta$, where $C'$ is the constant in (\ref{diferencia vortex}). Hence, for $t\in\mathcal{I}_k$ we have
\[\abs{\int_{\TT^3}\abs{v}^2-\int_{\TT^3}\abs{v_0}^2}\leq [(C'\lambda)^2+2(C'\lambda)]\,\abs{\Omega_k}< [\delta/8](4\cdot2^{-k}\lambda T)= 2^{-(k+1)}\delta\lambda T,\]
where we have used that $v=v_0$ outside $\Omega_k$ for $t\in\mathcal{I}_k$. By definition of $\mathcal{I}_k$, we see that $2^{-(k+1)}T\leq t$ for any $t\in\mathcal{I}_k$. We conclude that
\[1-3\delta\lambda t<\int_{\TT^3}\abs{v(x,t)}^2\,dx<1-\delta\lambda t.\]
Therefore, the weak solution $v$ is admissible. In addition, we see that we obtain a sequence of different solutions $\{v_i\}_{i=i_0}^\infty$ by repeating the construction with $\lambda_i=4^{-i}$ for sufficiently large $i_0$.

\section{Blowup}\label{section blowup}

In this final section we prove Theorem~\ref{teorema blowup} on the existence of H\"older continuous weak solutions of the 3d Euler equations that exhibit a singular set of maximal dimension. The proof makes use of some technical lemmas that are presented in Subsections~\ref{SS.buiblow} and~\ref{SS.lem}.

\subsection{Building blocks}\label{SS.buiblow}
The fundamental element in our construction is the following simple blowup, whose proof is postponed to Subsection~\ref{SS.lem}:
\begin{lemma}
	\label{un blowup}
	Let $0<\beta<1/3$ and let $q>2$. Let $a\in C^\infty(\RR,\RR^3)$ be a bounded map. Given $\varepsilon>0$, there exists a weak solution of the Euler equations $(v_\varepsilon,p_\varepsilon)$ in $\RR^3\times\RR$ such that:
	\begin{itemize}
		\item $(v_\varepsilon,p_\varepsilon)= (a,-\partial_t a\cdot x)$ outside $B(0,\varepsilon)\times(0,\varepsilon)$,
		\item the $q$-singular set of $v_\varepsilon$ is $\mathscr{S}^q_{v_\varepsilon}=\{(0,\varepsilon)\}$,
		\item $v_\varepsilon\in C^\beta_{\text{loc}}(\RR^3\times\RR\backslash \mathscr{S}^q_{v_\varepsilon})$, 
		\item the relative energy $e_\varepsilon(t)\coloneqq \norm{v_\varepsilon(\cdot,t)-a}^2_{L^2(\RR^3)}$ is continuous and so is the map $t\mapsto \int a\cdot(v_\varepsilon-a)dx$.
	\end{itemize}
	Furthermore, there exists a constant $C>0$ depending on $\norm{a}_{L^\infty}<\infty$ but not on $\varepsilon$ such that
	\begin{equation}
		\norm{v_\varepsilon(\cdot,t)-a(t)}^2_{L^2(\RR^3)}+\norm{p_\varepsilon(\cdot,t)+\partial_ta(t)\cdot x}_{L^1(\RR^3)}\leq C\varepsilon^3 \qquad \forall t\in \RR.
		\label{estimación normas varepsilon}
	\end{equation}
\end{lemma}

Once we known how to construct a single blowup, as stated in the previous lemma, we will use the following result to glue many of them together. Its proof is completely independent of Lemma~\ref{un blowup}.
\begin{lemma}
	\label{pegar blowups}
	Let $a\in C^\infty(\RR,\RR^3)$. Let $\{(v_i,p_i)\}_{i=0}^\infty$ be a sequence of weak solutions of the Euler equations in $\RR^3\times\RR$. Suppose that the sets $F_i$ given by the closure of
	\[ \{(x,t)\in \RR^4:(v_i,p_i)(x,t)\neq (a(t),-\partial_t a(t)\cdot x)\}\]
	are pairwise disjoint and that
	\[\sum_{i=0}^\infty\left(\norm{v_i-a}_{L^2(\RR^4)}+\norm{p_i+\partial_t a\cdot x}_{L^1(\RR^4)}\right)<\infty.\]
	Then, the pair $(v,p)$ given by
	\[v\coloneqq a+\sum_{i=0}^\infty(v_i-a), \qquad p\coloneqq -\partial_t a\cdot x+\sum_{i=1}^\infty (p_i+\partial_t a\cdot x)\]
	is a weak solution of the Euler equations.
\end{lemma}
\begin{proof}
	The hypotheses readily yield $v\in L^2_\text{loc}(\RR^4)$ and $p\in L^1_\text{loc}(\RR^4)$. Furthermore, it is easy to see that any partial sum 
	\[\T{v}_k\coloneqq a+\sum_{i=0}^k(v_i-a), \qquad \T{p}_k\coloneqq \partial_t a\cdot x+\sum_{i=0}^k(p_i-\partial_t a\cdot x)\] 
	is a subsolution. Indeed, fix $\chi_i\in C^\infty(\RR^4)$ such that $\chi_i^{-1}(\{0\})=F_i$ and consider $\theta_1\coloneqq \chi_1/(\chi_1+\chi_2)$ and $\theta_2\coloneqq 1-\theta_1$. They are well defined because $F_1\cap F_2=\varnothing$. We see that $\theta_1$ vanishes in $F_2$ and $\theta_2$ vanishes in $\theta_1$, so for any $\phi\in C^\infty_c(\RR^4,\RR^3)$ we have
	\begin{align*}
		\int_{\RR^4}\left(\partial_t\phi\cdot \T{v}_2+\nabla\phi:(\T{v}_2\otimes \T{v}_2+\T{p}_2\Id)\right)&=\int_{\RR^4}\left(\partial_t(\theta_1\phi)\cdot \T{v}_2+\nabla(\theta_1\phi):(\T{v}_2\otimes \T{v}_2+\T{p}_2\Id)\right)\\&\hspace{10pt}+\int_{\RR^4}\left(\partial_t(\theta_2\phi)\cdot \T{v}_2+\nabla(\theta_2\phi):(\T{v}_2\otimes \T{v}_2+\T{p}_2\Id)\right) \\
		&=\int_{\RR^4}\left(\partial_t(\theta_1\phi)\cdot v_1+\nabla(\theta_1\phi):(v_1\otimes v_1+p_1\Id)\right)\\&\hspace{10pt}+\int_{\RR^4}\left(\partial_t(\theta_2\phi)\cdot v_2+\nabla(\theta_2\phi):(v_2\otimes v_2+p_2\Id)\right)\\&=0,
	\end{align*}
	which follows from the definition of $(v_i,p_i)$ being a weak solution of the Euler equations using the test-function $\theta_i\phi$. An analogous argument proves that $\T{v}_2$ is weakly divergence-free. Therefore, $(\T{v}_2,\T{p}_2)$ is a weak solution of the Euler equations. Furthermore, iterating this argument we conclude that $(\T{v}_k,\T{p}_k)$ is a weak solution of the Euler equations for any $k\geq 1$. Since $\T{v}_k$ converges to $v$ in $L^2(K)$ and $\T{p}_k$ converges to $p$ in $L^1(K)$ for any compact subset $K\subset \RR^4$, it is easy to see that $(v,p)$ is a weak solution of the Euler equations.
\end{proof}

\subsection{Proof of \cref{un blowup}}\label{SS.lem}

	Note it suffices to prove the result for $\varepsilon=1$ because for $\varepsilon\neq 1$ we could simply take
	\[(v_\varepsilon,p_\varepsilon)(x,t)\coloneqq (v_1,p_1)(x/\varepsilon,t/\varepsilon)\]
	with $a$ rescaled accordingly. It is clear that this scaling preserves the $q$-singular set, that is,
	\[\mathscr{S}^q_{v_\varepsilon}=\{(x,t)\in \RR^3\times\RR:(x/\varepsilon,t/\varepsilon)\in \mathscr{S}^q_{v_1}\}=\{(0,\varepsilon)\}.\]
	Furthermore, since $p_1+\partial_t a\cdot x$ will be the only compactly supported solution of
	\begin{align*}
		-\Delta (p_1+\partial_ta\cdot x)&=\Div\Div(v_1\otimes v_1)=\Div\Div(v_1\otimes v_1-a\otimes a)\\
		&=\Div\Div[(v_1-a)\otimes(v_1-a)+(v_1-a)\otimes a+a\otimes(v_1-a)],
	\end{align*}
	standard Calderón-Zygmund estimates yield 
	\[\norm{p_1(\cdot,t)+\partial_ta(t)\cdot x}_{L^1(\RR^3)}\leq C\norm{v_1(\cdot,t)-a(t)}_{L^2(\RR^3)}^2.\]
	Thus, to prove (\ref{estimación normas varepsilon}) it suffices to show that $(v_1-a)\in L^\infty_tL^2_x$.
	
	From now on, we will assume $\varepsilon=1$ and drop the subscript $1$ for simplicity. By \cite{Gavrilov} there exists a nontrivial steady solution of the Euler equations with compact support $(u,\pi)\in C^\infty_c(\RR^3)$. By rescaling, we may assume that its support is contained in the ball $B(0,1/4)$. We define
	\begin{align*}
		U(x)&\coloneqq u\left(x-e_1/4\right)-u\left(x+e_1/4\right), \\
		P(x)&\coloneqq \pi\left(x-e_1/4\right)-\pi\left(x+e_1/4\right).
	\end{align*}
	Therefore, $(U,P)\in C^\infty_c(B(0,1))$ is a nontrivial steady solution of the Euler equations such that
	\[\int \xi\cdot U=0 \qquad \forall \xi\in \ker \gradsim.\]
	Hence, by \cref{invertir divergencia matrices} there exists $S_0\in C^\infty_c(B(0,1),\mathcal{S}^3)$ such that $\Div S_0=U$. We introduce a parameter $\alpha\in (-1,0)$ that will be fixed later and we define
	\[S\coloneqq (1+\alpha)(x\otimes U+U\otimes x)-(4+5\alpha)S_0.\]
	Since $\Div(x\otimes U+U\otimes x)=4U+x\cdot\nabla U$ and $U$ is divergence-free, we have
	\begin{equation}
		\label{div div autosimilar}
		\Div\Div S=\Div\left(-\alpha \,U+(1+\alpha)x\cdot\nabla U\right)=0.
	\end{equation}
	Furthermore, we see that the triplet $(U,P,S)$ satisfies
	\begin{equation}
		\label{ec autosimilar}
		-\alpha \,U+(1+\alpha)x\cdot\nabla U+\Div(U\otimes U+P\Id)=\Div S.
	\end{equation}
	Consider the self-similar ansatz 
	\begin{align*}
		\T{v}_0(x,t)&\coloneqq (1-t)^\alpha\, U\left(\frac{x}{(1-t)^{1+\alpha}}\right), \\
		\T{p}_0(x,t)&\coloneqq (1-t)^{2\alpha}P\left(\frac{x}{(1-t)^{1+\alpha}}\right), \\
		\T{R}_0(x,t)&\coloneqq (1-t)^{2\alpha}S\left(\frac{x}{(1-t)^{1+\alpha}}\right).
	\end{align*}
	It follows from \cref{ec autosimilar} that the triplet $(\T{v}_0,\T{p}_0,\T{R}_0)$ is a subsolution of the Euler equations in $\RR^3\times[0,1)$ (in which the Reynolds stress is not normalized to be trace-free, yet). Regarding the scaling, we see that
	\[\norm{\T{v}_0(\cdot,t)}_{L^r}=(1-t)^{\alpha+3(1+\alpha)/r}\norm{U}_{L^r}.\]
	We choose
	\[\alpha\coloneqq -\frac{3}{10}-\frac{3}{2(3+q)},\]
	which ensures that
	\[\lim_{t\to 1}\norm{\T{v}_0(\cdot,t)}_{L^2}=0, \qquad \lim_{t\to 1}\norm{\T{v}_0(\cdot,t)}_{L^q}=+\infty.\]
	
	Next, we fix $\chi\in C^\infty([0,1])$ that vanishes in a neighborhood of $0$ and is identically $1$ in a neighborhood of $1$. Since $(\T{v}_0,\T{p}_0,\T{R}_0)$ is a subsolution, is is easy to see that the triplet $(v_0,p_0,R_0)$ given by
	\begin{align*}
		v_0(x,t)&\coloneqq a(t)+\chi(t)\T{v}_0(x,t),\\
		p_0(x,t)&\coloneqq -\partial_t a\cdot x+\chi(t)^2\T{p}_0(x,t),\\
		R_0(x,t)&\coloneqq\chi(t)\T{R}_0+\chi'(t)(1-t)^{-1}S_0((1-t)^{-1/2}x)+\chi(t)(a\otimes \T{v}_0+\T{v}_0\otimes a)(x,t)
	\end{align*}
	is also a subsolution and that $\Div\Div R_0$ vanishes. Let $\mathcal{I}_k\coloneqq [1-2^{-k},1-2^{-(k+1)}]$ and $\Omega_k\coloneqq B(0,2^{-(1+\alpha)k})$. We see that the support of $R_0(\cdot,t)$ is contained in $\Omega_k$ for $t\in \mathcal{I}_k$. In fact, so is the support of $(v_0,p_0,R_0)(\cdot,t)$. In addition, we see that
	\[\norm{v_0-a}_{0,\mathcal{I}_k}\lesssim2^{-\alpha(k+1)}, \qquad
		\norm{R_0}_{0,\mathcal{I}_k}\lesssim2^{-2\alpha(k+1)},\]
	where the implicit constant depends on $\norm{a}_{L^\infty}$ but not on $k$.
	
	We would like to apply \cref{teorema integración convexa segunda parte}. While our current situation does not exactly meet the hypotheses of \cref{teorema integración convexa segunda parte},  this is not really an issue. In spite of the fact that our subsolution is not compactly supported, what matters is the support of the Reynolds stress, as argued in Subsection~\ref{SS.drop}. Furthermore, although $R_0$ is not trace-free, this can be easily solved in the proof of \cref{teorema integración convexa segunda parte}. After applying \cref{anular R en t0}, one simply has to replace $(\T{v}_0,\T{p}_0,\T{R}_0)$ by \[\left(\T{v}_0,\;\T{p}_0-\frac{1}{3}\tr(\T{R}_0),\; \T{R}_0-\frac{1}{3}\tr(\T{R}_0)\Id\right).\]
	
	Therefore, there exists a weak solution $(v,p)$ of the Euler equations in $\RR^3\times (0,1)$ with $v\in C^\beta_\text{loc}(\RR^3\times(0,1))$ and $(v,p)=(a,0)$ in $(\RR^3\backslash \Omega_k)\times\mathcal{I}_k$ for $k\geq 0$. Furthermore, since $(v_0,p_0,R_0)=(a,0,0)$ for $t$ sufficiently close to $0$, a careful revision of \cref{teorema integración convexa segunda parte} will convince us that $(v,p)=(a,0)$ in a neighborhood of $t=0$. Thus, we may extend $(v,p)$ to the interval $(-\infty,1)$.
	
	On the other hand, for $t\in \mathcal{I}_k$ we have
	\[
	\norm{v-a}_{0,\mathcal{I}_k}\leq \norm{v_0-a}_{0,\mathcal{I}_k}+\norm{v-v_0}_{0,\mathcal{I}_k}\leq  \norm{v_0-a}_{0,\mathcal{I}_k}+C\norm{R_0}_{0,\mathcal{I}_k}^{1/2}\leq C\,2^{-\alpha k},\]
	where the constant changes after each inequality and it is allowed to depend on $\norm{a}_{L^\infty}$ but not on $k$. Hence,
	\begin{equation}
		\label{norm difference v-a}
		\norm{v-a}_{0,\mathcal{I}_k}\abs{\Omega_k}\leq C\,2^{-(3+4\alpha)k}, \qquad \norm{v-a}_{0,\mathcal{I}_k}^2\abs{\Omega_k}\leq C\,2^{-(3+5\alpha)k}.
	\end{equation}
	Our choice of $\alpha$ ensures that both quantities go to $0$ when $k\to\infty$. We conclude that the maps $t\mapsto\norm{v(\cdot,t)-a}^2_{L^2(\RR^3)}$ and $t\mapsto \int a\cdot(v_\varepsilon-a)dx$ can be  extended by 0 to a continuous function in the whole $\RR$.
	
	Next, we show that we may extend our weak solution to $\RR^3\times\RR$ by setting $(v,p)=(a,0)$ for $t\geq 1$. For simplicity, we still denote this extension as $(v,p)$. It will be a weak solution in $\RR^3\times\RR$ if and only if for any solenoidal test-function $\phi\in C^\infty_c(\RR^3\times \RR)$ we have
	\[\int_\RR\int_{\RR^3}[\partial_t\phi\cdot v+\nabla\phi:(v\otimes v)]dx\hspace{0.5pt}dt=0.\]
	We split the integral:
	\begin{align*}
		\int_1^\infty\int_{\RR^3}[\partial_t\phi\cdot v+\nabla\phi:(v\otimes v)]dx\hspace{0.5pt}dt&=\int_1^\infty\int_{\RR^3}[\partial_t\phi\cdot a+\nabla\phi:(a\otimes a)]dx\hspace{0.5pt}dt\\ &=\int_1^\infty\int_{\RR^3}[\partial_t\phi(x,t)\cdot a(t)]dx\hspace{0.5pt}dt=0,
	\end{align*}
	where we have used that, for a fixed time $t$, $a(t)$ is just a constant vector and $\partial_t\phi(\cdot,t)$ is a compactly-supported divergence-free. Thus, the spatial integral vanishes for each time $t$. We conclude that $(v,p)$ will be a weak solution in $\RR^3\times\RR$ if and only if
	\begin{equation}
		\int_{-\infty}^1\int_{\RR^3}[\partial_t\phi\cdot v+\nabla\phi:(v\otimes v)]dx\hspace{0.5pt}dt=0.
		\label{extended field is weak solution}
	\end{equation}
	We choose a cutoff function $\theta\in C^\infty(\RR)$ that equals $0$ if $t\leq 0$ and equals $1$ if $t\geq 1/2$. For $j\geq 1$ consider $\theta_j(t)\coloneqq \theta(1+2^{j}(t-1))$. Since $(1-\theta_j)\phi\in C^\infty_c(\RR^3\times(-\infty,1))$ is solenoidal and $(v,p)$ is a weak solution in $\RR^3\times(-\infty,1)$, we see that
	\begin{equation}
		\int_{-\infty}^1\int_{\RR^3}[\partial_t\phi\cdot v+\nabla\phi:(v\otimes v)]dx\hspace{0.5pt}dt=\int_{1-2^{-j}}^1\int_{\RR^3}[\partial_t(\theta_j\phi)\cdot v+\nabla(\theta_j\phi):(v\otimes v)]dx\hspace{0.5pt}dt.
		\label{reescribir integral}
	\end{equation}
	Let us study the second term on the right-hand side. We fix $t\in \mathcal{I}_k$ and we write
	\[\int_{\RR^3}\nabla(\theta_j\phi):(v\otimes v)\,dx=\int_{\RR^3}\nabla(\theta_j\phi):(v\otimes v-a\otimes a)\,dx.\]
	Taking into account that \[v\otimes v-a\otimes a=(v-a)\otimes (v-a)+a\otimes(v-a)+(v-a)\otimes a\] and (\ref{norm difference v-a}), we surmise that for $t\in \mathcal{I}_k$
	\[\abs{\int_{\RR^3}\nabla(\theta_j\phi):(v\otimes v)\,dx}\leq C\norm{\phi}_1 2^{-(3+5\alpha)k}.\]
	Therefore,
	\[\abs{\int_{1-2^{-j}}^1\int_{\RR^3}\nabla(\theta_j\phi):(v\otimes v)\,dx\hspace{0.5pt}dt}\leq C\norm{\phi}_1\sum_{k=j}^\infty 2^{-(3+5\alpha)k}\abs{\mathcal{I}_k}\leq C\norm{\phi}_1 2^{-(4+5\alpha)j}\]
	because $3+5\alpha>0$ by our choice of $\alpha$. Regarding the first term on the right-hand side of (\ref{reescribir integral}), we split the integral into
	\[\int_{1-2^{-j}}^1\int_{\RR^3}\partial_t(\theta_j\phi)\cdot v\,dx\hspace{0.5pt}dt=\int_{1-2^{-j}}^1\int_{\RR^3}\partial_t(\theta_j\phi)\cdot a\,dx\hspace{0.5pt}dt+\int_{1-2^{-j}}^1\int_{\RR^3}\partial_t(\theta_j\phi)\cdot (v-a)\,dx\hspace{0.5pt}dt.\] 
	Again, the first term vanishes because if we keep $t$ fixed $a(t)$ is just a constant vector and $\partial_t(\theta_j\phi)(\cdot,t)$ is a compactly-supported divergence-free field. Thus, the spatial integral vanishes for each time $t$. Concerning the second term, we estimate
	\begin{align*}
		\abs{\int_{1-2^{-j}}^1\int_{\RR^3}\partial_t(\theta_j\phi)\cdot(v-a)\,dx\hspace{0.5pt}dt}&\leq \sum_{k=j}^\infty \norm{\partial_t(\theta_j\phi)}_0\norm{v-a}_{0,\mathcal{I}_k}\abs{\Omega_k}\abs{\mathcal{I}_k} \\
		&\leq C(\norm{\phi}_0+\norm{\partial_t\phi}_0)2^j\sum_{k=j}^\infty 2^{-(3+4\alpha)k}2^{-k}\\
		&\leq C(\norm{\phi}_0+\norm{\partial_t\phi}_0)2^{-(3+4\alpha)j}.
	\end{align*}
	Hence,
	\[\abs{\int_{-\infty}^1\int_{\RR^3}[\partial_t\phi\cdot v+\nabla\phi:(v\otimes v)]dx\hspace{0.5pt}dt}\leq C(\norm{\phi}_1+\norm{\partial_t\phi}_0)2^{-(3+4\alpha)j},\]
	where the constant may depend on $a$ but it is independent of $j$. Since $3+4\alpha>0$ and $j\geq 1$ is arbitrary, we conclude that (\ref{extended field is weak solution}) holds, so $(v,p)$ is a weak solution in $\RR^3\times\RR$.
	
	Regarding the $q$-singular set, fix a spatial ball $B(0,r)$ and let $t_k\coloneqq 1-2^{-k}$ for $k\geq 1$. For sufficiently large $k$ we have $\chi(t_k)=1$ and, since the velocity is unchanged at the endpoints of the intervals $\mathcal{I}_k$, we have
	\[v(x,t_k)=a+2^{-\alpha k}\hspace{1.5pt}U\left(2^{(1+\alpha)k}x\right).\]
	We compute
	\begin{align*}
		\norm{v(\cdot,t_k)}_{L^q(B(0,r))}&\geq \norm{v(\cdot,t_k)-a}_{L^q(B(0,r))}-\norm{a}_{L^q(B(0,r))}\\&\geq -\norm{a}_{L^\infty}\left(\frac{4}{3}\pi r^3\right)^{1/q}+ 2^{-[\alpha+3(1+\alpha)/q]k}\norm{U}_{L^q(\RR^3)},
	\end{align*}
	where we have used that for sufficiently large $k$ the ball $B(0,2^{(1+\alpha)k}r)$ contains the support of $U$. By our choice of $\alpha$ we have $\alpha+3(1+\alpha)/q<0$, so
	\[\lim_{k\to\infty}\norm{v(\cdot,t_k)}_{L^q(B(0,r))}=+\infty.\]
	Since the ball $B(0,r)$ is arbitrary, we see that $(0,1)\in \mathscr{S}^q_v$.
	
	Finally, since $v\in C^\beta_\text{loc}(\RR^3\times(-\infty,1))$ and $v=a$ in $(\RR^3\backslash B(0,2^{-k}))\times(1-2^{-k},1)$ for any $k\geq 0$, we conclude that $v\in C^\beta_\text{loc}(\RR^3\times\RR\backslash\{(0,1)\})$. In particular, the $q$-singular set reduces to $\{(0,1)\}$. This completes the proof of the lemma.

\subsection{Proof of \cref{teorema blowup}}
After a temporal rescaling, we may assume that $T=1$. After a translation, we may assume that $\overline{B}(0,4\rho)$ is contained in $U$ for sufficiently small $\rho>0$. By subtracting a time dependent constant, we may assume that $p_0(0,t)=0$. Let $a(t)=v_0(0,t)$. We glue $(v_0,p_0)$ and $(a,-\partial_ta\cdot x)$ using \cref{pegado subsoluciones}, obtaining a subsolution $(v_1,p_1,\R{R}_1)$ such that
\[(v_1,p_1,\R{R}_1)(x,t)=\begin{cases} (v_0,p_0,0) \hspace{54pt} \text{if }x\notin B(0,4\rho), \\ (a(t),-\partial_ta(t)\cdot x,0) \quad \text{if }x\in \overline{B}(0,3\rho).\end{cases}\]
It is not difficult to deduce from \cref{estimaciones pegado subsolucones} that by reducing $\rho$ we can obtain $\norm{v_1-v_0}_0$ and $\rho^3\|\R{R}_1\|_0^{1/2}$ arbitrarily small. 

We apply \cref{teorema integración convexa} to obtain a weak solution of the Euler equations $(v_2,p_2)$ that equals $(v_0,p_0)$ outside $B(0,4\rho)\times[0,1]$ and $(a,-\partial_t a\cdot x)$ in $\overline{B}(0,3\rho)\times[0,1]$. Note that we may choose a nonincreasing energy profile that is arbitrarily close to the original one but still satisfies (\ref{condición energía 1}) because $\rho^3\|\R{R}_1\|_0^{1/2}$ is arbitrarily small. Since $\norm{v_1-v_0}_0$ is arbitrarily small, we conclude that $\norm{v_2-v_0}_0$ may be chosen to be arbitrarily small.

Next, we construct the blowup in $B(0,\rho)\times(0,1]$. By~\cite{Besicovitch}, for any $k\geq 1$ there exists a function $f^k\in C^{2^{-k}}(B(0,\rho))$ taking values in $[1-2\cdot 4^{-k},1-4^{-k}]$ and whose graph $G^k$ has Hausdorff dimension $4-2^{-k}$. We want the $q$-singular set $\mathscr{S}^q$ of the final solution to contain all of these $G^k$ so that its Hausdorff dimension is $4$. To do that, we choose a sequence $\{x_i\}_{i=1}^\infty$ dense in $B(0,\rho)$ and we denote $\tau_i^k\coloneqq f^k(x_i)$. We see that $\{(x_i,\tau_i^k)\}_{i=1}^\infty$ is dense in $G^k$. Using \cref{un blowup} we will construct a sequence of blowups converging to each of the points $(x_i,\tau_i^k)$. Thus, these blowups would accumulate at every point in $\bigcup_{k\geq 1}G^k$, which means that this set will be contained in the singular set, as we wanted.

In order to glue the blowups given by \cref{un blowup} using \cref{pegar blowups}, they must have disjoint supports. Hence, we have study the geometry of the situation. Let 
\[t_{ij}^k\coloneqq \tau^k_i-4^{-(k+j)}\]
so that the sequence $\{t_{ij}^k\}_{j=1}^\infty$ is contained in $(1-4^{-(k-1)},1-4^{-k})$ and converges to $\tau_i^k$. We want to apply \cref{un blowup} to construct a blowup in
\[U_{ij}^k\coloneqq B(x_i,\varepsilon_{ij}^k)\times(t_{ij}^k,t_{ij}^k+\varepsilon_{ij}^k).\]
It is clear that choosing $\varepsilon_{ij}^k$ sufficiently small ensures that the sets $\overline{U}{}_{ij}^k$ are disjoint for a fixed $i$ and $k$, but it is not so clear if $i$ is not kept fixed. 

We will try to isolate the sequence corresponding to a fixed $i$. Let $L^k\coloneqq \norm{f^k}_{C^{2^{-k}}}$ and consider the sets
\[\mathcal{C}_{i}^k\coloneqq \left\{(x,t)\in \RR^3\times[1- 4^{-(k-1)},\tau_i^k]:\abs{t-\tau_i^k}\geq(L^k+1)\abs{x-x_i}^{2^{-k}}\right\}.\]
By the definition of $L^k$, for any $i\neq i'\geq 1$ we have $(x_{i},\tau_{i}^k)\notin \mathcal{C}_{i'}^k$. We define $j_0(1)=1$ and for $i>1$ we define $j_0(i)$ to be the minimum $j_0\geq 1$ such that 
\[\left[\mathcal{C}^k_i\cap[\RR^3\times(t^k_{ij_0},\tau^k_i)]\right]\cap \mathcal{C}^k_{i'}=\varnothing \qquad \forall i'<i.\]
As we have mentioned, $(x_{i},\tau_{i}^k)\notin \mathcal{C}_{i'}^k$ for $i'<i$, so there exists a neighborhood of $(x_{i},\tau_{i}^k)$ disjoint from the union of these sets. Since $\mathcal{C}^k_i\cap[\RR^3\times(t^k_{ij_0},\tau^k_i)]$ will be contained in this neighborhood of $(x_{i},\tau_{i}^k)$ for sufficiently large $j_0$, $j_0(i)$ is well defined. The point is that we will only add blowups for $j\geq j_0(i)$.

Next, let us define
\[\varepsilon_{ij}^k\coloneqq 4^{-2^k[i+j+k+\log_4(L^k+1)]}\delta,\] 
where $0<\delta\leq 1$ will be chosen later. Since $\varepsilon_{ij}^k\leq 4^{-(k+j+1)}$, we see that the $\{U^k_{ij}\}_{j=1}^\infty$ are pairwise disjoint for fixed $k,i$. Furthermore, this definition ensures that
\[4^{-(k+j)}-\varepsilon^k_{ij}\geq (L^k+1)(\varepsilon^k_{ij})^{2^{-k}},\]
which means that $U_{ij}^k\subset \mathcal{C}_i^k$. Therefore, it follows from the definition of $j_0(i)$ that the sets
\[\{U^k_{ij}:i,k\geq 1, j\geq j_0(i)\}\]
are pairwise disjoint, as claimed.

Let $(v_{ij}^k,p_{ij}^k)$ be the weak solution of the Euler equations given by \cref{un blowup} using the parameter $\varepsilon^k_{ij}$. After a translation, we may assume that the set where $(v_{ij}^k,p_{ij}^k)\neq (a,-\partial_t a\cdot x)$ is contained in $U^k_{ij}$. We define
\begin{align*}
	v_3&\coloneqq a+\sum_{k,i\geq 1}\sum_{j\geq j_0(i)}(v_{ij}^k-a), \\ \qquad p_3&\coloneqq -\partial_ta\cdot x+\sum_{k,i\geq 1}\sum_{j\geq j_0(i)}(p_{ij}^k+\partial_ta\cdot x).
\end{align*}
By \cref{pegar blowups}, the pair $(v_3,p_3)$ is a weak solution of the Euler equations in $\RR^3\times\RR$. Note that any point $(x,t)$ not in the closure of $\bigcup_{k\geq 1}G^k$ has a neighborhood $V$ that intersects only a finite number of the $U^k_{ij}$. We conclude that the $q$-singular set of the weak solution is the closure of the union of the $q$-singular sets of the $v_{ij}^k$, that is,
\[\mathscr{S}^q=\{(x_i,t_{ij}^k+\varepsilon_{ij}^k):k,i\geq 1, j\geq j_0(i)\}\cup\bigcup_{k\geq 1}G^k\cup \left(\,\overline{B}(0,\rho)\times\{1\}\right)\]
and $v\in C^\beta_\text{loc}(\RR^4\backslash \mathscr{S}^q)$. Regarding the energy, let us consider the partial sums
\begin{align*}
	\T{v}_{ij}^k&\coloneqq a+\sum_{k'=1}^k\sum_{i'=1}^i\sum_{j'= j_0(i')}^j(v_{i'j'}^{k'}-a), \\
	e_{ij}^k(t)&\coloneqq \int_{B(0,2\rho)}\abs{\T{v}_{ij}^k(x,t)}^2dx.
\end{align*}
Taking into account the identity $v^2=(v-a)^2+a^2+2a\cdot(v-a)$, we may write
\[e_{ij}^k(t)=\frac{32}{3}\pi\rho^3a(t)^2+\sum_{k'=1}^k\sum_{i'=1}^i\sum_{j'= j_0(i')}^j\left(\int|v_{i'j'}^{k'}-a|^2dx+2\int a\cdot(v_{i'j'}^{k'}-a)\,dx\right)\]
because the support of the $(v_{i'j'}^{k'}-a)$ are pairwise disjoint. We see that $e_{ij}^k$ is continuous by \cref{un blowup}. Furthermore, by Equation~(\ref{estimación normas varepsilon}) and our choice of $\varepsilon_{ij}^k$ it converges uniformly to $\int_{B(0,2\rho)}\abs{v_4(x,t)}^2dx$ which is, therefore, continuous. In addition, it can get arbitrarily close to $\frac{32}{3}\pi\rho^3a(t)^2$ by reducing $\delta>0$.

Let us glue this blowup to $(v_2,p_2)$. Since $x_i\in B(0,\rho)$ and we may assume that $\varepsilon_{ij}^k\leq \rho$ by further reducing $\delta$, we see that $U_{ij}^k\subset B(0,2\rho)\times(0,1)$, so $(v_4,p_4)=(a,-\partial_t a\cdot x)$ outside $B(0,2\rho)\times(0,1)$. Hence, it glues well with $(v_2,p_2)$. We conclude that there exists a weak solution of the Euler equations $(v_4,p_4)$ that equals $(v_0,p_0)$ outside $B(0,4\rho)\times(0,1)$ and has a $q$-singular set $\mathscr{S}^q\subset \overline{B}(0,\rho)\times(0,1]$ with Hausdorff dimension $4$. In addition, $v\in C^\beta_\text{loc}((\RR^3\times[0,1])\backslash \mathscr{S}^q)$. Furthermore, $\norm{v_4(\cdot,0)-v_0(\cdot,0)}_{C^0}$ is arbitrarily small and $t\mapsto \int \abs{v_4(x,t)}^2dx$ is continuous and arbitrarily close to $t\mapsto \int \abs{v_0(x,t)}^2dx$. 

To complete the proof it suffices to modify $(v_4,p_4)$ in $(B(0,3\rho)\backslash\overline{B}(0,2\rho))\times[0,1]$ to ensure that the energy profile is nonincreasing. Let $\T{e}(t)\coloneqq \int \abs{v_4(x,t)}^2dx$ and fix a nonincreasing function $e(t)>e_4(t)$. It is easy to obtain a sequence of strictly positive smooth functions $\{\delta_k\}_{k=1}^\infty$ whose sum is $e(t)-\T{e}(t)$ and such that $\norm{\delta_k}_{L^\infty}\lesssim 2^{-k}$. Let $r_k\coloneqq r_0\rho \hspace{0.5pt}2^{-k/3}$. By reducing $r_0$ if necessary, we can find a sequence of pairwise disjoint balls $\overline{B}(x_k,r_k)\subset B(0,3\rho)\backslash \overline{B}(0,2\rho)$.

Fix $k\geq 1$ and let $e_k(t)\coloneqq a(2^{-k/3}t)^2\abs{B(0,r_0\rho)}+2^k\delta_k(2^{-k/3}t)$. We use \cref{teorema integración convexa} to construct a weak solution of the Euler equations $(u_k,\pi_k)$ that equals $(a(2^{-k/3}t),-(\partial_t a)(2^{-k/3}t)\cdot x)$ for $x\notin B(0,r_0\rho)$ and such that $\int_{B(0,r_0\rho)}\abs{u_k}^2dx=e_k(t)$. In addition, $v\in C^\beta$.

Finally, we consider
\begin{align*}
	v_5(x,t)&\coloneqq a(t)+\sum_{k=1}^\infty [u_k(2^{k/3}(x-x_k),2^{k/3}t)-a(t)], \\
	p_5(x,t)&\coloneqq -\partial_t a(t)\cdot x+\sum_{k=1}^\infty [\pi_k(2^{k/3}x,2^{k/3}t)+\partial_t a(t)\cdot x].
\end{align*}
The pair $(v_5,p_5)$ is a weak solution of the Euler equations that equals $(a,-\partial_t a\cdot x)$ for $x\notin B(0,3\rho)\backslash \overline{B}(0,2\rho)$ and $v_5\in C^\beta$. Regarding the energy
\begin{align*}
	\int_{B(0,3\rho)\backslash \overline{B}(0,2\rho)}\abs{v_5(x,t)}^2dx&=a(t)^2\abs{B(0,3\rho)\backslash \overline{B}(0,2\rho)}\\&+\sum_{k=1}^\infty\int_{B(x_k,r_k)}\left[|u_5(2^{k/3}(x-x_k),t)|^2-a(t)^2\right]\,dx \\
	&=a(t)^2\abs{B(0,3\rho)\backslash \overline{B}(0,2\rho)}+\sum_{k=1}^\infty \delta_k(t) \\
	&= \int_{B(0,3\rho)\backslash \overline{B}(0,2\rho)}\abs{v_4(x,t)}^2dx+e(t)-\T{e}(t).
\end{align*}
We glue $(v_5,p_5)$ to $(v_4,p_4)$, obtaining the desired weak solution of the Euler equations $(v_6,p_6)$. Note that $\norm{(v_6-v_0)(\cdot,0)}_{C^0}$ can be taken to be arbitrarily small because we can do so with $\norm{v_2-v_0}_0$ and $e-\T{e}$. This completes the proof.

\section*{Acknowledgements}
The authors are very grateful to the reviewers for their suggestions and corrections on a previous version of the manuscript. This work has received funding from the European Research Council (ERC) under the European Union's Horizon 2020 research and innovation programme through the grant agreement~862342 (A.E.). It is partially supported by the grants CEX2023-001347-S, RED2022-134301-T and PID2022-136795NB-I00 (A.E. and D.P.-S.) funded by MCIN/AEI/10.13039 /501100011033, and Ayudas Fundaci\'on BBVA a Proyectos de Investigaci\'on Cient\'ifica 2021 (D.P.-S.). J.P.-T. was partially supported by an FPI grant CEX2019-000904-S-21-4 funded by MICIU/AEI and by FSE+.


\newcommand\cV{\mathcal{V}}

\appendix 

\section{H\"older and Besov spaces}\label{A:Holder}
Let $\Omega$ be an open subset of Euclidean space. We denote by $C^0(\Omega)$ the set of bounded continuous functions on $\Omega$, which we equip with the supremum norm, denoted by $\norm{f}_0\coloneqq \sup_{x\in \Omega}\abs{f(x)}$. More generally, for any $N\geq0$ we define the space $C^N(\Omega)$ as the set of functions that have bounded continuous derivatives of any order $k\leq N$. On this space, we define the following seminorms and norms, respectively:
\[[f]_N\coloneqq \max_{\abs{\beta}=N}\norm{D^\beta f}_0, \hspace{40pt}\norm{f}_N\coloneqq \sum_{j=0}^N [f]_j.\]
Here $\beta$ denotes a multi-index and $\abs{\beta}$ denotes its length. Given $N\geq 0$ and $\alpha\in(0,1)$, we define the Hölder space $C^{N+\alpha}(\Omega)$ as the set of functions $f\in C^N(\Omega)$ such that the following quantity is finite: 
\[[f]_{N+\alpha}\coloneqq \max_{\abs{\beta}=N}\sup_{x\neq y} \frac{\abs{D^\beta f(x)-D^\beta f(y)}}{\abs{x-y}^\alpha}.\]
This set becomes a Banach space when equipped with the following norm:
\[\norm{f}_{N+\alpha}\coloneqq \norm{f}_N+[f]_{N+\alpha}.\]
When we work in a subset $E\subset \Omega$ instead of the whole $\Omega$, we will write $\norm{\cdot}_{N;E}$. When the functions also depend on time, we also take the supremum in $t\in [0,T]$.

The Hölder norms satisfy the following inequalities:
\begin{equation}
	\label{holder desigualdad 1}
	[f]_s\leq C\left(\varepsilon^{t-s}[f]_t+\varepsilon^{-s}\norm{f}_0\right)
\end{equation}
for $0\leq s \leq t$ and all $\varepsilon >0$, and
\begin{equation}
	\label{holder estimación producto}
	[fg]_s\leq [f]_s\norm{g}_0+\norm{f}_0[g]_s
\end{equation}
for $0\leq s\leq 1$. The constant $C$ only depends on the Hölder exponents involved and on the domain $\Omega$. For the applications in this article, since $f$ will be a compactly supported function defined on the whole $\mathbb R^n$, $C$ will just depend on the Hölder exponents. From (\ref{holder desigualdad 1}) with $\varepsilon=\norm{f}_0^{1/r}[f]_r^{-1/r}$ we obtain the following interpolation inequality:
\begin{equation}
[f]_s\leq C\norm{f}_0^{1-s/r}[f]_r^{s/r}.
\label{interpolation inequality}
\end{equation}

Let $\beta$ by a multi-index. By induction on $\abs{\beta}$ and the rule for the first derivative of a product, one easily deduces
\[D^\beta (fg)=\sum_{\abs{\gamma}+\abs{\delta}=\abs{\beta}}C_{\abs{\beta},\gamma,\delta}\,D^\gamma f\,D^\delta g,\]
from which it immediately follows that
\begin{equation}
	\label{holder estimate product N-th}
	[fg]_N\leq C_N\sum_{j=0}^N[f]_j\,[g]_{N-j}.
\end{equation}

The following proposition is standard:

\begin{proposition}
	\label{schauder estimates}
	Let $N\in \mathbb{N}$ and $\alpha\in (0,1)$. Let $f\in C^{\hspace{1pt}N,\alpha}_c(\mathbb{R}^m,\mathbb{R})$ and $F\in C^{\hspace{1pt}N,\alpha}_c(\mathbb{R}^m,\mathbb{R}^m)$. There exists a constant $C=C(N,\alpha)$ such that the potential-theoretic solutions of
	\[\Delta \phi =f, \qquad \Delta \psi =\Div F\]
	satisfy
	\[\norm{\phi}_{N+2+\alpha}\leq C\norm{f}_{N+\alpha}, \qquad \norm{\psi}_{N+1+\alpha}\leq C\norm{F}_{N+\alpha}.\]
\end{proposition} 
Note that in the previous proposition, one does not get any information about the $C^\alpha$ norm of the solution, aside from estimating it by higher-order norms. For this, we will need to introduce negative regularity spaces. Let us consider a Littlewood--Payley decomposition, e.g.\ as in~\cite[Section 2.2]{BCD}. For this, we take smooth radial functions $\chi,\varphi:\RR^3\to[0,1]$, whose supports are contained in the ball $B(0,\frac43)$ and in the annulus $\{\frac34 <|\xi|<\frac83\}$ respectively, with the property that
\[
\chi(\xi)+\sum_{N=0}^\infty\varphi(2^{-N}\xi)=1
\]
for all~$\xi\in\RR^3$. In terms of the Fourier multipliers $P_{<}:=\chi(D)$ and $P_N:=\varphi(2^{-N}D)$, the Besov norm $B^s_{\infty,\infty}$ (which is equivalent to the H\"older norm $C^s$ if $s\in\RR^+\backslash\NN$, and strictly weaker if $s\in\NN$) can be written as
\begin{equation}\label{E.BesovF}
	\|f\|_{B^{s}_{\infty,\infty}}:=\|P_{<}f\|_0+\sup_{N\geq0}2^{Ns}\|P_Nf\|_0\,.
\end{equation}
Here $s\in\RR$ is any real number. Again, when dealing with time-dependent functions, we consider the supremum in time of $\|f(t)\|_{B^{s}_{\infty,\infty}}$.

\section{Some auxiliary estimates}


The first lemma of this appendix shows that we can find a cutoff function with well-behaved bounds on its derivatives:
\begin{lemma}
	\label{cutoff}
	Let $A\subset \mathbb{R}^n$ be a measurable set and let $r>0$. There exists a cutoff function $\varphi_r\in C^\infty(\mathbb{R}^n,[0,1])$ whose support is contained in $A+B(0,r)$ and such that $\varphi_r\equiv 1$ in a neighborhood of $A$. Furthermore, for any $N\geq 0$ we have \[\left\|\varphi_r\right\|_N\leq C(N,n)\,r^{-N}.\]
\end{lemma}
\begin{proof}
	We choose a nonnegative function $\psi\in C^\infty_c(B(0,1))$ such that $\int \psi=1$. For $\varepsilon>0$ we denote
	$$A_\varepsilon\coloneqq A+B(0,\varepsilon), \qquad \psi_\varepsilon (x)=\varepsilon^{-n}\,\psi(x/\varepsilon).$$
	Note that $\int \psi_\varepsilon =1$. The desired cutoff function is:
	$$\varphi_r(x)\coloneqq \left(\mathds{1}_{A_{r/2}}\ast \psi_{r/2}\right)(x)=\int_{A_{r/2}}\psi_{r/2}\left(x-y\right)\;dy.$$
	It is easy to see that its support is contained in $A+B(0,r)$ and that $\varphi_r\equiv 1$ in a neighborhood of $A$. Furthermore, it is smooth and
	$$\begin{aligned}\partial^\alpha \varphi_r(x)&=\left(\mathds{1}_{A_{r/2}}\ast \partial^\alpha\psi_{r/2}\right)(x)=\int_{A_{r/2}}\partial^\alpha_x\left[(r/2)^{-n}\;\psi\left(\frac{x-y}{r/2}\right)\right]\;dy\\&=(r/2)^{-|\alpha|}\int_{A_{r/2}}(r/2)^{-n}(\partial^\alpha\psi)\left(\frac{x-y}{r/2}\right)\;dy.\end{aligned}$$
	Hence,
	\begin{align*}
		\left|\partial^\alpha \varphi_r(x)\right| &\leq (r/2)^{-|\alpha|}\int_{A_{r/2}}(r/2)^{-n}\left|(\partial^\alpha\psi)\left(\frac{x-y}{r/2}\right)\right|\;dy\\ &\leq (r/2)^{-|\alpha|}\int_{\mathbb{R}^n}(r/2)^{-n}\left|(\partial^\alpha\psi)\left(\frac{x-y}{r/2}\right)\right|\;dy\\&\leq(r/2)^{-|\alpha|}\int_{\mathbb{R}^n}\left|\partial^\alpha\psi(y)\right|\,dy, 
	\end{align*}
	from which the result follows.
\end{proof}

The second instrumental lemma provides a bound for the solution to a transport equation. The proof is standard, see e.g.~\cite{Quinto}.

\begin{lemma}
	\label{transporte}
	Let $f\in C^\infty(\RR^3\times\RR)$ be the solution of the transport equation
	\[\begin{cases}\partial_t f+v\cdot\nabla f=g, \\ f(\cdot,0)=f_0\end{cases}\]
	for some vector field $v\in C^\infty(\RR^3\times\RR,\RR^3)$ and $g\in C^\infty(\RR^3\times\RR)$. Then, for $0\leq \alpha\leq 1$ and $\abs{t}\norm{v}_1\leq 1$ we have
	\[\norm{f(\cdot,t)}_\alpha\leq e^\alpha\left(\norm{f_0}_\alpha+\int_0^t\norm{g(\cdot,\tau)}_\alpha d\tau\right).\]
\end{lemma}

In this article we also need and extension theorem for Hölder continuous functions defined on a domain $\Omega$:
\begin{theorem}
	\label{extension holder functions}
	Let $N\geq0$ and $\alpha\in(0,1)$. Let $\Omega\subset\RR^\dim$ be a domain with smooth boundary. Then, there exists a linear map $T:C^{N+\alpha}(\Omega)\to C^{N+\alpha}(\RR^\dim)$ such that
	\begin{itemize}
		\item $Tf=f$ on $\Omega$ for each $f\in C^{N+\alpha}(\Omega)$ and
		\item the norm of $T$ is bounded by a constant depending only on $\Omega$ and $N$.
	\end{itemize}
\end{theorem}
\begin{proof}
The proof is essentially \cite[Lemma 6.37]{GT} and is based on a rectification of the boundary and a reflection. We must, however, make some remarks. In \cite{GT} they assume $N\geq 1$ because they are considering sets with less regular boundaries. In the case of a smooth boundary, the result also holds for $C^\alpha$ functions. We must warn the reader that they are using the notation $C^{N+\alpha}\left(\overline{\Omega}\right)$ to denote Hölder continuous functions because they use $C^{N+\alpha}\left(\Omega\right)$ to denote locally Hölder continuous functions. Finally, it is also interesting to note that the form of operator $T$ itself does not depend on $N$ and $\alpha$, although its norm as an operator $C^{N+\alpha}(\Omega)\to C^{N+\alpha}(\RR^\dim)$ does, of course. This holds provided that we follow the construction of \cite{GT} but use smooth parametrizations of $\Omega$, instead of less regular ones.
\end{proof}

The last result of this appendix is a lemma that establishes a bound for a Besov norm of functions that are compactly supported in a collared neighborhood of an $(n-1)$-dimensional surface. This can be seen as a Poincaré inequality with negative regularity.

\begin{lemma}
	\label{frecuencias bajas de dominio pequeño}
	Let $\Omega\subset \RR^\dim$, $\dim\geq 2$, be a bounded domain with smooth boundary. Let $\alpha\in (0,1)$ and let $r>0$ be sufficiently small (depending on $\Omega$). Consider a function $f\in C_c^\infty(\RR^\dim)$ supported in $\{x\in\RR^\dim:0<\dist(x,\Omega)<r\}$. We have
	\[\norm{f}_{B^{-1+\alpha}_{\infty,\infty}}\leq C(\Omega,\alpha)\,r^{1-\alpha}\norm{f}_0.\]
\end{lemma}
\begin{proof}
	Let us denote $U \coloneqq \{x\in\RR^\dim:0<\dist(x,\Omega)<r\}$.  We begin by computing
	\[\abs{\Proy_N f(x)}\leq 2^{\dim N}\int_U \abs{h(2^N(x-y))}\abs{f(y)}\,dy\lesssim \norm{f}_0\, 2^{\dim N}\int_U \langle2^N(x-y)\rangle^{-8n}\,dy,\]
	where we have used that $|h(x)|\lesssim \langle x\rangle^{-8n}$ because $h$ is in the Schwartz class. Here $P_Nf$ stands for the nonhomogeneous dyadic blocks in the Littlewood-Paley decomposition of $f$, and $h$ is the corresponding convolution kernel.
 
 We claim that the following estimate holds:
	\begin{equation}\label{E.intApp}
 2^{\dim N}\int_U \langle2^N(x-y)\rangle^{-8n}\,dy\lesssim \min\{1,2^Nr\}.
 \end{equation}
	Using this, the bound for $f$ readily follows:
	\begin{align*}
	\norm{f}_{B^{-1+\alpha}_{\infty,\infty}}&=\sup_{N}2^{N(-1+\alpha)}\norm{\Proy_N f}_{L^\infty}\\ 
		&\lesssim \sup_{2^N<r^{-1}}2^{N(-1+\alpha)}(\norm{f}_0 2^Nr)+\sup_{2^N\geq r^{-1}}2^{N(-1+\alpha)}\norm{f}_0\\
		&\lesssim r^{1-\alpha}\norm{f}_0.
	\end{align*}
	Therefore, the problem is reduced to estimating the integral above. 

Let 
\[
U_R \coloneqq \{x\in\RR^\dim:0<\dist(x,\Omega)<R\}.
\]
For each point $x\in\partial\Omega$ and for each small enough~$R>0$, there is a  boundary normal chart
\[
X_{x,R}:Q_R\to U_{2R}\,,
\]
where $Q_R:=[0,R)\times(-R,R)^{n-1}$, such that $\Psi_{x,R}(0)=x$ and
\[
X_{x,R}(Q_R)\supset \{y\in \overline{U_R}: |x-y|<R/2\}\,.
\]
Recall that, by the definition of boundary normal coordinates, the pullback $X_{x,R}^*g_0$ of the Euclidean metric~$g_0:=(\delta_{ij})$ to this chart satisfies $X_{x,R}^*g_0(0)=g_0$. One can also assume that $\nabla X_{x,R}(0)$ is the identity matrix, and that 
\begin{equation}\label{E.distApp}
|x-y|>\frac R{4}    
\end{equation}
for all $y\in X_{x,R}(Q_{R/2})$ and all $x\not \in X_{x,R}(Q_{R/2})\cup\Omega$.

Since $\partial\Omega$ is a smooth compact hypersurface of~$\RR^n$, it is standard that there is some small enough $R>0$ and a finite collection of charts $\{ X^j:=X_{x_j,R}\}_{j=1}^J$ as above such that $\{U^j:=X^j(Q_{R/2})\}_{j=1}^J$ is a cover of the set $\overline{U_{R/2}}$ and which satisfy
\begin{equation}\label{E.C2App}
\|(X^j)^*g_0-g_0\|_{C^2(Q_R)}<\frac1{100}
\end{equation}
and
\begin{equation}\label{E.dist2App}
|X^j(z)-X^j(\tilde z)|\geq \frac12|z-\tilde z|
\end{equation}
for all $z,\tilde z\in Q_R$. Moreover, the distance between the point $X^j(z)$ and the boundary is comparable to its first coordinate, in the sense that $$
\frac{z_1}2\leq \dist(X^j(z),\Omega)\leq 2z_1
$$
for all $z\in Q_R$.

Let us now estimate the integral~\eqref{E.intApp} over the subset
$U^j\cap U$, with $x\in U$. If $x\not\in X^j(Q_{R})$, by~\eqref{E.distApp}, one immediately has
\[
2^{\dim N}\int_{U^j\cap U} \langle2^N(x-y)\rangle^{-8n}\,dy\lesssim 2^{\dim N}\langle 2^NR\rangle^{-8n}|U^j\cap U|\lesssim r\leq \{1,2^N r\}\,,
\]
where we have also used that $|U|\lesssim r$. If $x\in X^j(Q_{R})$, one can write $x=X^j(\tilde z)$ for some $\tilde z\in Q_R$. By~\eqref{E.dist2App}, one then has
\begin{align*}
   2^{\dim N}\int_{U^j\cap U} &\langle2^N(x-y)\rangle^{-8n}\,dy \lesssim 2^{\dim N}\int_{0}^{2r}\int_{(-R,R)^{\dim -1}}\langle2^N(z-\tilde z)\rangle^{-8n}\, J_{X^j}(z)\,dz_1\, dz'\\
   & \lesssim 2^{\dim N}\int_{0}^{2r}\int_{(-R,R)^{\dim -1}}\langle2^N(z_1-\tilde z_1)\rangle^{-4n}\, \langle2^N(z'-\tilde z')\rangle^{-4n}\,dz_1\, dz'\,.
\end{align*}
Here we have used the notation $z=(z_1,z')\in [0,R)\times (-R,R)^{n-1}$ and the fact that the Jacobian $J_{X^j}$ is bounded by~\eqref{E.C2App}. 

We can now carry out the integrations in $z_1$ and $z'$ separately. Since 
\[
2^{(\dim-1) N}\int_{(-R,R)^{\dim -1}} \langle2^N(z'-\tilde z')\rangle^{-4n}\,dz_1\, dz'\leq 2^{(\dim-1) N}\int_{\RR^{\dim -1}} \langle2^N(z'-\tilde z')\rangle^{-4n}\,dz_1\, dz'\lesssim 1\,,
\]
the estimate follows from the fact that
\begin{align*}
2^{ N}\int_{0}^{2r}\langle2^N(z_1-\tilde z_1)\rangle^{-4n}\,dz_1&\leq 2^{ N}\int_{-4r}^{4r}\langle2^N s\rangle^{-4n}\,ds\\
&= \int_{-2^{N+2}r}^{2^{N+2}r}\langle s\rangle^{-4n}\,ds\lesssim\min\{1,2^N r\}\,.
\end{align*}
The estimate~\eqref{E.intApp} then follows by summing over~$j$.
\end{proof}



\end{document}